\documentclass[a4paper,12pt]{amsart}

\usepackage{amssymb,amsmath,amsthm,mathtools,thmtools,xcolor,comment,enumerate, graphicx, ascmac,mathrsfs,framed,url,tikz,tikz-cd, appendix, upquote}
\usetikzlibrary{calc,decorations.pathmorphing}

\mathtoolsset{showonlyrefs}

\usepackage{geometry}
\allowdisplaybreaks[4]

\usepackage[
    bookmarksnumbered,
    colorlinks,
    hypertexnames=false,
    pdfdisplaydoctitle,
    pdfusetitle,
    unicode
]{hyperref}

\usepackage{upref}
\usepackage{zref-clever}
\usepackage{zref-xr}
\zcsetup{
    cap,
    nameinlink=false,
    sort=false
}
\AddToHook{env/table/begin}{\zcsetup{reftype=table}}

\theoremstyle{definition}
\newtheorem{thm}{Theorem}[subsection]
\newtheorem{defi}[thm]{Definition}
\newtheorem{prop}[thm]{Proposition}

\newtheorem{coro}[thm]{Corollary}
\newtheorem{lemma}[thm]{Lemma}
\newtheorem{exam}[thm]{Example}

\newtheorem{rema}[thm]{Remark}

\newtheorem*{prob*}{Problem}

\newcommand{\setzcreftypefor}[2]{
    \AddToHook{env/#1/begin}{\zcsetup{reftype=#2}}
}

\makeatletter
\newcommand{\labelitem}[2]{
    \item[(#1)]
    \phantomsection
    \def\@currentlabel{(#1)}
    \label{#2}
}
\makeatother

\setzcreftypefor{thm}{theorem}
\setzcreftypefor{defi}{definition}
\setzcreftypefor{prop}{proposition}
\setzcreftypefor{coro}{corollary}
\setzcreftypefor{lemma}{lemma}
\setzcreftypefor{exam}{example}
\setzcreftypefor{exams}{example}
\setzcreftypefor{rema}{remark}
\setzcreftypefor{claim}{claim}

\newcounter{thmAref}
\newcounter{thmBref}
\newcounter{thmCref}
\newcounter{thmDref}

\newtheorem{thmA*}[thmAref]{Theorem}
\newtheorem{thmB*}[thmBref]{Theorem}
\newtheorem{thmC*}[thmCref]{Theorem}
\newtheorem{thmD*}[thmDref]{Theorem}

\AddToHook{env/thmA*/begin}{\zcsetup{reftype=theorem}}
\AddToHook{env/thmB*/begin}{\zcsetup{reftype=theorem}}
\AddToHook{env/thmC*/begin}{\zcsetup{reftype=theorem}}
\AddToHook{env/thmD*/begin}{\zcsetup{reftype=theorem}}

\numberwithin{equation}{section}

\DeclareMathOperator{\rank}{\mathrm{rank}}
\DeclareMathOperator{\diam}{\mathrm{diam}}
\DeclareMathOperator{\Mag}{\mathrm{Mag}}
\DeclareMathOperator{\MH}{\mathrm{MH}}
\DeclareMathOperator{\MC}{\mathrm{MC}}

\let\Im\relax
\DeclareMathOperator{\Im}{\mathrm{Im}}

\DeclareMathOperator{\id}{\mathrm{id}}

\newcommand{\ones}{\mathbf{1}}

\newcommand{\midd}{\mathrel{}\middle |\mathrel{}}
\newcommand{\abs}[1]{\left\lvert#1\right\rvert}

\title{Combinatorial criteria for diagonality in magnitude homology}

\author{Ryutaro Misawa}
\address[Ryutaro Misawa]{Division of Mathematics\\
Graduate School of Information Sciences\\
Tohoku University\\
6-3-09 Aramaki-Aza-Aoba, Aoba-ku, Sendai, 980-8579\\
Japan}
\email{misawa.ryutaro.q2@dc.tohoku.ac.jp}

\author{Yusaku Nishimura}
\address[Yusaku Nishimura]{Graduate School of Fundamental Science and Engineering\\ 
Waseda University\\
3-4-1 Ookubo, Shinjuku-ku, Tokyo, 169-8555\\
Japan}
\email{n2357y@ruri.waseda.jp}

\author{Takuya Saito}
\address[Takuya Saito]{Institute for Chemical Reaction Design and Discovery\\
Hokkaido University\\
Kita 21 Nishi 10, Kita-ku, Sapporo, Hokkaido, 001-0021\\
Japan}
\email{saito@icredd.hokudai.ac.jp}

\author{Yuuho Tanaka}
\address[Yuuho Tanaka]{Mathematical Sciences Program\\
Department of Science and Technology\\
Faculty of Science and Technology\\
Oita University\\
700 Dannoharu, Oita, 870-1192\\
Japan}
\email{tanaka-yuuho@oita-u.ac.jp}

\date{\today}

\begin{document}

\begin{abstract}
A graph is diagonal if its magnitude homology vanishes whenever the homological degree differs from the length degree.
Ivanov and Mukoseev showed that, for finite connected graphs, diagonality is equivalent to Koszulity of the distance algebra over every field.
We introduce a combinatorial condition called reducivity and prove that it implies diagonality by a recursive reduction of magnitude chain complexes using chain homotopy equivalences.
Our criteria for reducivity establish diagonality for weakly modular graphs and basis graphs of even $\Delta$-matroids and $q$-matroids.
Moreover, the diagonality of weakly modular graphs implies that, for finite connected locally weakly modular graphs, diagonality is equivalent to simple connectivity of the triangle-square complex and to off-diagonal vanishing in homological degree two.
These conditions are also equivalent to Koszulity and to quadraticity of the distance algebra over every field.
This extends the local-to-global characterization of weakly modular graphs, a discrete analogue of the Cartan--Hadamard theorem.
Furthermore, using reducivity, we classify diagonal strongly regular graphs with smallest eigenvalue $-2$ and diagonal Taylor graphs.
We also show that the Hilbert series of a distance algebra does not determine its quadraticity, even for distance-regular graphs.
\vspace{2mm}


\end{abstract}

\maketitle

\makeatletter
\def\l@subsection{\@tocline{2}{0pt}{2pc}{5pc}{}}
\def\l@subsubsection{\@tocline{3}{0pt}{3pc}{7pc}{}}
\def\l@paragraph{\@tocline{4}{0pt}{1pc}{7pc}{}}
\makeatother
\setcounter{tocdepth}{3}
\tableofcontents

\section{Introduction}\label{sec:intro}
In Lawvere's framework \cite{Law73}, metric spaces can be regarded as enriched categories.
Magnitude is an invariant of enriched categories that extends the Euler characteristic of finite categories \cite{Lei08} and applies in particular to metric spaces \cite{Lei13}.
For a finite connected graph, magnitude is a rational function defined by shortest-path distances, and its formal power series expansion has integer coefficients \cite{Lei19}.
Hepworth and Willerton \cite{HW17} introduced magnitude homology, whose graded Euler characteristic recovers this formal power series.

Diagonal graphs are of particular interest because their magnitude determines their magnitude homology up to isomorphism.
For a graph $G$, we write $\MH_{n,\ell}(G)$ for its magnitude homology with integer coefficients in homological degree $n$ and length degree $\ell$.
A graph $G$ is called \emph{diagonal} if $\MH_{n,\ell}(G)=0$ whenever $n\ne\ell$ (see \zcref{subsec:maghom}).
In addition, diagonal graphs have several interesting algebraic properties.
The magnitude cohomology ring of a finite metric space, introduced by Hepworth \cite{Hep22}, determines the space up to isometry.
For a diagonal graph, this ring is a quotient of the path algebra of the graph by relations involving shortest paths of length $2$ \cite{Hep22}.
Asao and Ivanov \cite{AI24} expressed the magnitude homology and cohomology of finite quasimetric spaces in terms of Tor and Ext, respectively, over certain graded algebras called distance algebras.
Ivanov and Mukoseev \cite{IM26} showed that a graph $G$ is diagonal if and only if its distance algebra is Koszul over every field.
In this case, the distance algebra and the path cochain algebra are Koszul and are quadratic duals of each other \cite{IM26}.
Thus, finding combinatorial conditions for diagonality is also important for the study of the Koszulity of distance algebras (see \zcref{subsec:maghom}).
Hepworth and Willerton asked the following question.

\begin{prob*}[{\cite[Section~1.5, p.~36]{HW17}}]
\label{prob:diagonal-characterization}
Can diagonal graphs be characterized in graph-theoretic terms?
\end{prob*}
By the result of Ivanov and Mukoseev, this problem is equivalent to asking for a graph-theoretic characterization of graphs whose distance algebras are Koszul over every field.

Hepworth and Willerton proved that complete graphs and trees are diagonal, that diagonality is preserved under Cartesian products, and that the join of any two nonempty graphs is diagonal \cite{HW17}.
Gu \cite{Gu18} used algebraic Morse theory to prove the diagonality of pawful graphs and the icosahedral graph.
Bottinelli and Kaiser \cite{BK21} proved that median graphs are diagonal, and Asao \cite{Asao2026} proved that quasi-median graphs are diagonal.

For fixed endpoints and a fixed length, Asao and Izumihara \cite{AI21} related magnitude homology to the relative homology of a pair of simplicial complexes, with a shift in homological degree.
Tajima and Yoshinaga \cite{TY23} applied discrete Morse theory to the associated complexes and gave examples of diagonal graphs of diameter $2$ that are not pawful.
Asao, Hiraoka, and Kanazawa \cite{AHK24} studied the relation between girth and off-diagonal vanishing and showed that the probability that a random graph is diagonal exhibits a phase transition.

Asao and Wakatsuki \cite{AW24} gave a criterion for diagonality of finite geodetic metric spaces in terms of $4$-cuts.
In particular, for finite geodetic graphs, off-diagonal vanishing in homological degree $2$ implies diagonality.
Ivanov and Mukoseev \cite{IM26} gave a combinatorial characterization of off-diagonal vanishing in homological degree $2$ in terms of relations between paths.

In this paper, we develop a recursive reduction of magnitude chain complexes using partitions of the vertex set defined by shortest paths.
The main step is to establish a chain homotopy equivalence that reduces both the homological degree and the length degree by one.
This leads to a combinatorial sufficient condition for diagonality, which we call \emph{reducivity}.
Using this condition, we prove the diagonality of weakly modular graphs, basis graphs of even $\Delta$-matroids, and basis graphs of $q$-matroids.
In particular, the diagonality of weakly modular graphs gives homological and algebraic characterizations of weak modularity within the class of locally weakly modular graphs, in addition to the topological characterization given by the discrete Cartan--Hadamard theorem.
We also apply our criteria to classify diagonal graphs in two classes of distance-regular graphs.

\smallskip
\paragraph{\textit{Reducive graphs and homological reduction}}
Our condition for diagonality is based on the following reduction of magnitude chain complexes.
Let $G=(V,E)$ be a finite connected graph, and write $d$ for its shortest-path distance and $\Gamma(u)$ for the neighborhood of a vertex $u$.
Fix $u\in V$ and $S\subseteq\Gamma(u)$.
The target set $T_u(S)$ consists of the vertices that can be reached from $u$ by a shortest path whose first step is to a vertex of $S$, i.e.,
\begin{equation}
    T_u(S)=\{x\in V\mid d(u,x)=1+d(w,x)\text{ for some }w\in S\}.
\end{equation}
For adjacent vertices $u,w$ and subsets $S\subseteq\Gamma(u)\setminus\{w\}$ and $R\subseteq\Gamma(w)$, put $T_u^S(w)\coloneq T_u(\{w\})\setminus T_u(S)$.
We consider the following partition of the vertex set:
\begin{equation}
    V=T_u^S(w)\sqcup T_w(R).
\end{equation}
Fix endpoints $u,u'$ and consider the subcomplex $\MC_{*,\ell}(G;uT_u(S),u')$ generated by vertex sequences with these endpoints whose second entry belongs to $T_u(S)$.
If the above partition holds, then for $\ell\ge2$ we have a chain homotopy equivalence
\begin{equation}
    \frac{\MC_{*,\ell}(G;uT_u(S\cup\{w\}),u')}{\MC_{*,\ell}(G;uT_u(S),u')} \simeq \MC_{*-1,\ell-1}(G;wT_w(R),u'),
\end{equation}
where the complex on the left is a quotient complex and the complex on the right retains the original differential (see \zcref{lem:chain-isom_2,lem:chain-isom_2_rho}).
This reduction decreases both the homological and length degrees by one and hence preserves the difference between them.
A graph is called \emph{reducive} if we can choose a family $\mathcal{R}_u\subseteq 2^{\Gamma(u)}$ containing $\Gamma(u)$ for each vertex $u$ so that the following condition holds at every vertex.
For each nonempty $S\in\mathcal{R}_u$, there exist $w\in S$ and $R\in\mathcal{R}_w$ such that $S\setminus\{w\}\in\mathcal{R}_u$ and $V=T_u^{S\setminus\{w\}}(w)\sqcup T_w(R)$ (see \zcref{def:family}).
Using the long exact sequence in homology, we prove the following theorem by induction on the length degree and the cardinality of $S$.

\begin{thmA*}[\zcref{thm:reducible}]\label{thm:main-reducive}
Every reducive graph is diagonal.
\end{thmA*}

By this theorem, to prove that a graph is diagonal, it suffices to construct families of subsets of vertex neighborhoods satisfying the above conditions.
Thus, reducivity also gives a combinatorial sufficient condition for the distance algebra of a graph to be Koszul over every field.
Moreover, the class of reducive graphs is closed under Cartesian products (see \zcref{prop:reducive-cartesian-product}).

\smallskip
\paragraph{\textit{Combinatorial criteria and applications}}
We characterize reducivity in terms of distances and common neighbors (see \zcref{thm:tryy}).
We also give a sufficient condition for reducivity using suitable families of cliques that cover vertex neighborhoods (see \zcref{lem:basis-reduction}).
Applying these criteria, we obtain the following result.

\begin{thmB*}\label{thm:main-classes}
Every graph in the following three classes is reducive and hence diagonal.
\begin{enumerate}[(i)]
    \item Weakly modular graphs (\zcref{prop:wmg-reducible}).
    \item Basis graphs of even $\Delta$-matroids (\zcref{prop:even-delta-reducible}).
    \item Basis graphs of $q$-matroids (\zcref{thm:q-reducive}).
\end{enumerate}
\end{thmB*}

Weakly modular graphs are defined by the triangle and quadrangle conditions on shortest-path distances (see \zcref{subsec:wm-mesh}).
For these graphs, the family at each vertex in the definition of reducivity can be taken to consist of all subsets of its neighborhood.
Weakly modular graphs include median and quasi-median graphs, and weak modularity and pawfulness are equivalent for graphs of diameter at most $2$ (see \zcref{prop:wmg-paw-med}).
Thus, \zcref{thm:main-classes}~(i) gives a common proof of the known diagonality results for these classes.

For basis graphs of even $\Delta$-matroids and $q$-matroids, we use the exchange axioms for sets and subspaces, respectively, to construct the families of cliques required by our criterion (see \cite{Chepoi2007,JurPel2018,CeJu2024}).
Basis graphs of even $\Delta$-matroids can be realized as the $1$-skeletons of the corresponding polytopes \cite{Chepoi2007}.
The diagonality of basis graphs of ordinary matroids, pointed out by 
Koizumi \cite{Koizumi2026}, is also covered by \zcref{thm:main-classes}~(ii).
Bilinear forms graphs are graphs on matrix spaces over finite fields, in which two matrices are adjacent if and only if their difference has rank one.
They can be realized as basis graphs of $q$-matroids.
Bilinear forms graphs of diameter at least $3$ are not weakly modular, but they are diagonal by \zcref{thm:main-classes}~(iii) (see \zcref{prop:bilinear-q-basis,exam:bilinear-not-wmg}).

\smallskip
\paragraph{\textit{Nonpositive curvature and local-to-global properties}}
Weakly modular graphs form a common generalization of many classes of graphs arising in the study of combinatorial nonpositive curvature.
They include the $1$-skeletons of $\mathrm{CAT}(0)$ cubical complexes, systolic complexes, and bucolic complexes \cite[Section~1]{CCHO20}.
In particular, median graphs are exactly the $1$-skeletons of $\mathrm{CAT}(0)$ cubical complexes \cite[Theorem~2.7]{CCHO20}.
Asao \cite{Asa21} showed that magnitude homology vanishes in all positive homological degrees for $\mathrm{CAT}(0)$ geodesic metric spaces.
On the other hand, Bottinelli and Kaiser \cite{BK21} proved the diagonality of median graphs equipped with the graph metric.
Our result in \zcref{thm:main-classes}~(i) extends the latter result and gives off-diagonal vanishing for the $1$-skeletons of the above classes of complexes.

The triangle-square complex of a graph $G$ is the CW complex obtained from $G$ by attaching a $2$-cell along the boundary of each triangle and each chordless square.
In their discrete analogue of the Cartan--Hadamard theorem, Chalopin, Chepoi, Hirai, and Osajda \cite[Theorem~3.1]{CCHO20} characterized weak modularity of locally weakly modular graphs in terms of simple connectivity of their triangle-square complexes.
Combining this characterization with the results of Ivanov and Mukoseev \cite{IM26} and \zcref{thm:main-classes}~(i), we obtain the following equivalences.

\begin{thmC*}[\zcref{thm:triangle-square}]\label{thmC}
Let $G$ be a finite connected locally weakly modular graph.
The following conditions are equivalent:
\begin{enumerate}[(i)]
    \item The triangle-square complex is simply connected.
    \item The graph $G$ is weakly modular.
    \item The graph $G$ is diagonal.
    \item For every integer $\ell>2$, we have $\MH_{2,\ell}(G)=0$.
    \item The distance algebra of $G$ is Koszul over every field $\mathbb{F}$.
    \item The distance algebra of $G$ is quadratic over every field $\mathbb{F}$.
\end{enumerate}
\end{thmC*}

\zcref{thmC} extends this local-to-global characterization of weak modularity by adding homological and algebraic conditions.
In particular, within the class of finite connected locally weakly modular graphs, quadraticity of the distance algebra over every field implies Koszulity over every field.
Equivalently, under the hypothesis of local weak modularity, off-diagonal vanishing in homological degree $2$ implies off-diagonal vanishing in all homological degrees.

\begin{figure}[htbp]
\centering
\begin{tikzpicture}[
    line width=0.6pt,
    every node/.style={font=\small,align=center},
    class/.style={draw,thick,inner sep=4.5pt},
    drg/.style={draw,rounded corners=2pt,inner sep=6pt,font=\itshape},
    yscale=1.8
]
    \node[class] (sup) at (8.5,4.5) {Superclass};
    \node[class] (sub) at (4.5,4.5) {Subclass};
    \node[class] (diagonal) at (2,3.5) {Diagonal};
    \node[class] (reducive) at (2,2.7) {Reducive};
    \node[class] (wmg) at (0.5,1) {Weakly modular};
    \node[class] (meshed) at (6,2.7) {Meshed};
    \node[class] (median) at (-0.25,0) {Median};
    \node[class] (grid) at (-0.5,-1) {Grid};
    \node[class] (pawful) at (-2,0) {Pawful};
    \node[class] (mat) at (6.5,0) {Matroid basis};
    \node[class] (dmat) at (4.25,1) {Even $\Delta$-matroid basis};
    \node[class] (qmat) at (8,1) {$q$-matroid basis};
    \node[class] (tree) at (-1.75,-1) {Tree};
    \node[class] (star) at (-3,-2) {Star};
    \node[class] (path) at (-1.75,-2) {Path};
    \node[drg] (ham) at (1.25,-1) {Hamming};
    \node[drg] (hypcube) at (1.25,-2) {Hypercube};
    \node[drg] (john) at (5.5,-1) {Johnson};
    \node[drg] (gr) at (8,-1) {Grassmann};
    \node[drg] (bili) at (9,0) {Bilinear};
    \node[drg] (conf) at (-4.5,-1) {Conference${}^{\ddagger}$};
    \node[class] (dualpolar) at (-4.5,0) {Dual polar${}^{\dagger}$};
    \node[drg] (halvedcube) at (4,-2) {Halved cube};
    \node (c4) at (-4.25,-2.5) {\textit{$4$-cycle $C_4$}};
    \node (house) at (-5,1.25) {House};
    \node (icos) at (-4,0.75) {\textit{Icosahedral}};
    \node (leg) at (2.75,4.5) {Legend:};
    \draw[->,thick] (reducive) --node[auto=left] {\scriptsize\zcref{thm:reducible}}(diagonal);
    \draw[->,thick] (wmg)--node[auto=left, pos=0.4, rotate=15, sloped] {\scriptsize\zcref{prop:wmg-reducible}}(reducive);
    \draw[->,dashed] (wmg)--(meshed);
    \draw[->,dashed] (median)--(wmg);
    \draw[->,thick] (pawful)--node[auto=left, pos=0.2] {\scriptsize \zcref{prop:wmg-paw-med}} (wmg);
    \draw[->,dotted, thick] (mat)--(dmat);
    \draw[->,dashed] (dmat)-- (meshed);
    \draw[->,thick] (qmat)--node[auto=right] {\scriptsize \zcref{prop:q-basis-meshed}} (meshed);
    \draw[->,thick] (dmat)--node[auto=left, rotate=-15, sloped] {\scriptsize \zcref{prop:even-delta-reducible}} (reducive);
    \draw[->,thick] (qmat)--node[auto=left, rotate=-10, sloped] {\scriptsize \zcref{thm:q-reducive}} (reducive);
    \draw[->,dashed] (ham)--(wmg);
    \draw[->,dotted, thick] (hypcube)--(ham);
    \draw[->,dotted, thick] (star)--(tree);
    \draw[->,dotted, thick] (tree)--(median);
    \draw[->,dotted, thick] (star)--(pawful);
    \draw[->,thick] (conf)--node[auto=left, pos=0.3] {\scriptsize \zcref{thm:pawful-item2}} (pawful);
    \draw[->,dashed] (john)--(wmg);
    \draw[->,dotted, thick] (john)-- (mat);
    \draw[->,thick] (gr)--node[auto=left, pos=0.6] {\scriptsize \zcref{prop:Grassmann-wmg}} (wmg);
    \draw[->,dotted, thick] (gr)--(qmat);
    \draw[->,thick] (bili)--node[auto=left, pos=0.7] {\scriptsize \zcref{prop:bilinear-q-basis}} (qmat);
    \draw[->,dashed] (dualpolar)--(wmg);
    \draw[->,dashed] (halvedcube)--(wmg);
    \draw[->,dotted, thick] (halvedcube)--(dmat);
    \draw[->,dotted, thick] (path)--(tree);
    \draw[->,dotted, thick] (path)--(grid);
    \draw[->,dotted, thick] (grid)--(median);
    \draw[->,dotted, thick] (hypcube)--(grid);
    \draw[->,dotted, thick] (c4)--(hypcube);
    \draw[->,dotted, thick] (c4)--(pawful);
    \draw[->,thick] (house)--node[auto=left, pos=0.2, rotate=5, sloped] {\scriptsize \zcref{exam:house-tryy}}(reducive);
    \draw[->,thick] (icos)--node[auto=left, pos=0.2, rotate=10, sloped] {\scriptsize \zcref{exam:icosahedral-tryy}}(reducive);
    \draw[->,thick] (sub)--node[auto=left] {\scriptsize proved here}(sup);
    \draw[->,dashed] (4.5,3.5)--node[auto=left] {\scriptsize shown in \cite{CCHO20}}(8.5,3.5);
    \draw[->,dotted] (4.5,4)--node[auto=left] {\scriptsize standard inclusions}(8.5,4);
\end{tikzpicture}
\caption{Inclusions among graph classes. An arrow means that the class at its tail is contained in the class at its head. Names in italics denote distance-regular graphs. Solid arrows show inclusions proved in this paper, with references to the relevant results. Dashed arrows show inclusions from \cite{CCHO20}. Dotted arrows show standard inclusions. (${}^{\dagger}$ Dual polar graphs of classical polar spaces are distance-regular. ${}^{\ddagger}$ Excluding the $5$-cycle.)}
\label{fig:graph-class-basic-inclusions}
\end{figure}
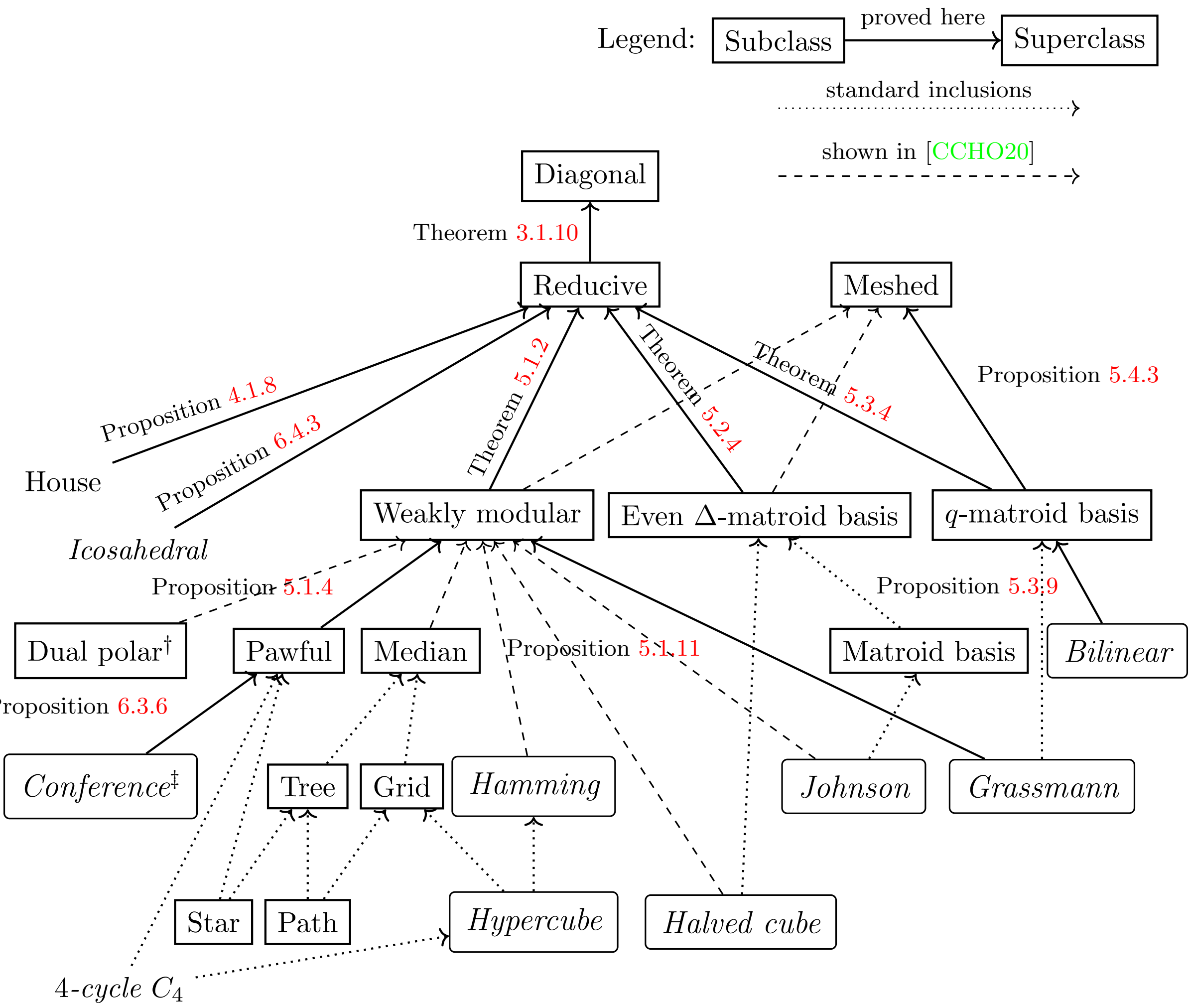

Meshedness is a condition on shortest-path distances that is weaker than weak modularity (see \zcref{subsec:nonred-mesh}).
Every graph in the three classes in \zcref{thm:main-classes} is meshed (see \zcref{prop:wmg-matroid-meshed,prop:q-basis-meshed}).
It is therefore natural to ask whether the hypothesis of local weak modularity in \zcref{thmC} can be replaced by meshedness.
However, we construct a meshed graph $G$ on $10$ vertices whose triangle-square complex is simply connected but for which $\MH_{3,4}(G)\ne0$ (see \zcref{exam:meshed-not-diagonal}).
Thus, the hypothesis cannot be replaced by meshedness.
The inclusions among the main graph classes are shown in \zcref{fig:graph-class-basic-inclusions}.

\smallskip
\paragraph{\textit{Distance-regular graphs}}
Finally, we apply our criteria to distance-regular graphs.
These graphs have high symmetry with respect to distance and hence form another natural class of graphs for the study of magnitude homology.
This combinatorial symmetry is described by the intersection array \cite{vDKT16}.
Strongly regular graphs form the class of distance-regular graphs of diameter $2$.
Every graph in this class other than the $5$-cycle has smallest eigenvalue at most $-2$ \cite[Section~1.1.10]{BvM22}.
The strongly regular graphs with smallest eigenvalue $-2$ have been classified (see \cite{S1968,BvM22}).
We also consider Taylor graphs, a class of distance-regular graphs of diameter $3$.
This class includes the icosahedral graph, whose diagonality was established by Gu \cite{Gu18}, as well as graphs arising from the Higman--Sims group and the Conway group $\mathrm{Co}_3$, both sporadic finite simple groups (see \cite{BGK19,IM24}).
Using our criteria together with structural results on these graphs, we obtain the following classifications.

\begin{thmD*}\label{thm:main-classification}
The following statements hold.
\begin{enumerate}[(i)]
    \item A connected strongly regular graph whose adjacency matrix has smallest eigenvalue $-2$ is diagonal if and only if it is isomorphic to neither the Petersen graph nor the Shrikhande graph (see \zcref{thm:srg-minus2-wmg}).
    \item A Taylor graph is diagonal if and only if it is not isomorphic to the $6$-cycle $C_6$ (see \zcref{cor:taylor-wmg-classification}).
\end{enumerate}
\end{thmD*}

The Chang graphs provide further examples of reducive graphs that are not weakly modular (see \zcref{exam:chang-tryy}).

The intersection array of a distance-regular graph determines its magnitude \cite{AW24}.
Moreover, the intersection array and the distance between two endpoints determine the ranks of the magnitude chain groups with these endpoints and their Euler characteristic in each length degree (see \zcref{prop:drg-zeta-inverse}).
Hence, for a diagonal distance-regular graph, these data determine the magnitude homologies with fixed endpoints up to isomorphism.
It follows from \cite{IM26} that, for such a graph, these data also determine the dimensions of the graded Tor groups between the corresponding right and left simple modules over the distance algebra for any choice of the coefficient field.

In general, distance-regular graphs with the same intersection array can have different magnitude homology.
Gu \cite{Gu18} gave two graphs with the same intersection array but with different ranks of $\MH_{3,4}$.
Since $\MH_{2,3}$ vanishes for both graphs, this example does not settle whether the rank of $\MH_{2,3}$ is determined by the intersection array.
We give examples of graphs with the same intersection array for which the rank of $\MH_{2,3}$ differs (see \zcref{exam:96-19-2-4}).
Moreover, although the distance algebras of these graphs have the same Hilbert series, some are quadratic over every field, whereas others are not quadratic over any field.
Thus, even for distance-regular graphs, the Hilbert series of a distance algebra does not determine its quadraticity.

\smallskip
\paragraph{\textit{Organization of the paper}}
In \zcref{sec:pre}, we review the basics of magnitude, magnitude homology, and distance algebras.
In \zcref{sec:reducive}, we construct a reduction of magnitude chain complexes and prove that reducive graphs are diagonal.
In \zcref{sec:reducive-criteria}, we give two criteria for reducivity.
In \zcref{sec:class}, we apply these criteria to three classes of graphs and discuss the relation between diagonality and meshedness.
In \zcref{sec:drg}, we study diagonality of distance-regular graphs and the invariants determined by their intersection arrays.

\smallskip
\paragraph{\textit{Acknowledgments}}
The authors thank Junnosuke Koizumi for pointing out the diagonality of matroid basis graphs and for helpful comments.
The authors also thank Yoh Kitajima for helpful comments on the relationship between weakly modular graphs and combinatorial nonpositive curvature.
The research of Ryutaro Misawa was supported by a JSPS Grant-in-Aid for JSPS Fellows (26KJ0570).
The research of Yusaku Nishimura was supported by the Yoshida Scholarship Foundation through the Doctor 21 program and by a Grant-in-Aid for Young Scientists (Early Bird) from the Waseda Research Institute for Science and Engineering.
The research of Takuya Saito was supported by JST CREST Grant Number JPMJCR24Q2.
The research of Yuuho Tanaka was supported by JSPS KAKENHI Grant Number JP25K23339.

\smallskip
\paragraph{\textit{Human–AI collaboration}}
This work was developed through extensive interaction between the authors and OpenAI’s GPT-6.2 Astra and GPT-6.2 Sol, and Google Gemini 4 Neon.
AI assistance included developing mathematical arguments, performing computations, and drafting and revising the manuscript.
The authors critically evaluated the AI-generated material, independently verified the mathematical arguments and references, and take full responsibility for the paper.

\section{Preliminaries}\label{sec:pre}
\subsection{Notation}
For an integer $m\ge1$, put $[m]=\{1,\ldots,m\}$, and set $[0]=\emptyset$.
For $0\le i\le m$, let $\binom{[m]}{i}$ be the family of $i$-element subsets of $[m]$.
We write a singleton $\{u\}$ as $u$ when there is no ambiguity.
We write $\abs{S}$ for the cardinality of a set $S$.
For two sets $X$ and $Y$, the symmetric difference is $X\triangle Y=(X\cup Y)\setminus (X\cap Y)$.
We write $\id$ for the identity map and $\Im(\phi)$ for the image of a map $\phi$.

We write $\mathbb{Q}, \mathbb{R}, \mathbb{C}, \mathbb{F}_q, \mathbb{Z},\mathbb{Z}_{\ge 0}$, and $\mathbb{Z}_{> 0}$ for the rational, real, and complex fields, the field with $q$ elements, the ring of integers, and the sets of nonnegative and positive integers, respectively.
For a ring $\mathbb{F}$, we write $\mathbb{F}[t]$ and $\mathbb{F}[\![t]\!]$ for the polynomial ring and the formal power series ring in $t$ over $\mathbb{F}$, respectively.
If $W$ is a subspace of a vector space $U$, we write $W\le U$.
For sets $X$ and $Y$ and a ring $\mathbb{F}$, let $\mathbb{F}^X$ be the set of tuples of elements of $\mathbb{F}$ indexed by $X$, and let $\mathbb{F}^{X\times Y}$ be the set of $\abs{X}\times \abs{Y}$ matrices with entries in $\mathbb{F}$ whose rows and columns are indexed by $X$ and $Y$, respectively.
We write $\Xi_{ij}$ for the $(i,j)$-entry of a matrix $\Xi$.
We write $I$ for the identity matrix, $J$ for the all-ones matrix, $O$ for the zero matrix, and $\ones$ for the all-ones column vector.
For a set $X$, we write $\mathbb{Z}X$ for the free abelian group with basis $X$.
For a polynomial or formal power series $\alpha(t)$ in $t$, we write $[t^\ell]\alpha(t)$ for the coefficient of $t^\ell$.

\subsection{Graphs}\label{subsec:graph}

Throughout this paper, a graph $G=(V,E)=(V(G),E(G))$ is finite, simple, and undirected.
Thus $E\subseteq \binom{V}{2}$. We write $uw$ for an edge $\{u,w\}$ when there is no ambiguity.
Two vertices $u,w\in V$ are \emph{adjacent} if $uw\in E$.
Two graphs $G,G'$ are \emph{isomorphic} if there is a bijection $\phi:V(G)\to V(G')$ that induces a bijection $E(G)\to E(G')$ given by $uw\mapsto \phi(u)\phi(w)$.
Given a graph $G=(V,E)$, a graph $G'=(V',E')$ is a \emph{subgraph} of $G$ if $V'\subseteq V$ and $E'\subseteq E$.
The \emph{induced subgraph} of $G$ on $S\subseteq V$ has vertex set $S$ and contains every edge of $G$ with both endpoints in $S$.

For an integer $m\in \mathbb{Z}_{>0}$, a graph is \emph{complete} if it is isomorphic to $K_m=([m],\binom{[m]}{2})$.
A subset $C\subseteq V$ of the vertex set of a graph $G=(V,E)$ is a \emph{clique} if the induced subgraph on $C$ is complete.
We also regard the empty set as a clique.
For an integer $m\in \mathbb{Z}_{>0}$, a \emph{path} of length $m-1$ is a graph isomorphic to $([m],\{12,23,\dots,(m-1)m\})$.
For an integer $m\geq3$, an $m$-\emph{cycle} is a graph isomorphic to $C_m=([m],\{12,23,\dots,(m-1)m,m1\})$.
A \emph{tree} is a connected graph with no cycle as a subgraph.
The \emph{complete multipartite graph} $K_{m_1,\ldots,m_s}$ is formed by partitioning the vertex set into $s\ge2$ nonempty sets of sizes $m_1,\ldots,m_s$ and joining any two vertices in different sets.
When there are $m$ sets, each of size $r$, we write $K_{m\times r}$.
In particular, $K_{m,n}$ is a complete bipartite graph, and $K_{1,n}$ is a \emph{star}.

The \emph{Cartesian product} $G\square G'$ of graphs $G,G'$ has vertex set $V(G)\times V(G')$. Two vertices $(u,x)$ and $(w,y)$ are adjacent if
\begin{equation}
\bigl(u=w\ \text{and}\ xy\in E(G')\bigr)
\quad\text{or}\quad
\bigl(uw\in E(G)\ \text{and}\ x=y\bigr).
\end{equation}
The \emph{join} of two nonempty graphs $G,G'$ is obtained from their disjoint union by adding all edges between $V(G)$ and $V(G')$.

A graph $G$ is \emph{connected} if $V(G)$ is nonempty and $G$ contains a path with endpoints $x,y$ for every pair $x,y \in V(G)$.
Throughout the paper, graphs are assumed to be connected unless stated otherwise.
For any two vertices $x,y \in V$, the \emph{distance} $d_G(x,y)=d(x,y)$ in $G$ is the length of a shortest path joining $x,y$.
The distance function is thus $d: V \times V \to \mathbb{Z}_{\geq 0}$.
The \emph{diameter} of $G$ is
\begin{equation}
D \coloneq \max_{x,y \in V}d(x,y).
\end{equation}
For a vertex $u\in V$ and an integer $i\in\mathbb{Z}_{\ge 0}$, the set $\Gamma_i(u)$ of vertices at distance $i$ from $u$ is defined by
\begin{equation}
\Gamma_i(u) \coloneq \{x\in V\mid d(x,u)=i\}.
\end{equation}
For $i=1$, we write $\Gamma(u)$ for $\Gamma_1(u)$ and call it the \emph{neighborhood} of $u$.
By definition, $u \notin \Gamma(u)$.
For $S\subseteq V$, put $\Gamma(S)\coloneq\bigcup_{s\in S}\Gamma(s)$, with $\Gamma(\emptyset)=\emptyset$.
The \emph{degree} of $u$ is $\abs{\Gamma(u)}$, written $\deg(u)$.
When $\abs{\Gamma_r(u)}$ is independent of $u\in V$, we denote this common value by $k_r$, and write $k = k_1$.
A graph is \emph{regular} if it has constant degree $k$.

For each distance $i$ in a graph $G$, define the \emph{distance matrix} $A_i \in \{0,1\}^{V\times V}$ by 
\begin{equation}
    A_i = (\delta_{i,d(x,y)})_{(x,y)\in V^2}.
\end{equation}
Here $\delta_{i,d(x,y)}$ is the Kronecker delta.
By definition, $A_0$ is the identity matrix $I$, and $A_1$ is the \emph{adjacency matrix} of $G$.
We write $A$ for the adjacency matrix $A_1$.

\subsection{Magnitude of graphs}\label{subsec:mag}
For a formal variable $t$, define the \emph{zeta matrix} of $G$ by $Z_G\coloneq (t^{d(x,y)})_{(x,y)\in V\times V}\in\mathbb{Q}[t]^{V\times V}$.
The constant term of $Z_G$ is $I$, so $Z_G$ is invertible over $\mathbb{Q}[\![t]\!]$. The unique \emph{magnitude weighting} $\omega\in\mathbb{Q}[\![t]\!]^V$ is defined as $\omega\coloneq Z_G^{-1}\ones$.
The \emph{magnitude} of $G$ is $\Mag(G;t)=\sum_{x\in V}\omega_x=\ones^\top Z_G^{-1}\ones$.

\begin{prop}\label{prop:valency-constant}
Suppose that $k_r=\abs{\Gamma_r(u)}$ is independent of $u$ for each $0\le r\le D$.
Then
\begin{equation}
    \Mag(G;t)=\frac{\abs{V}}{\sum_{r=0}^D k_rt^r}
    =\frac{\abs{V}^2}{\sum_{x,y\in V}t^{d(x,y)}}.
\end{equation}
\end{prop}

\begin{proof}
Since $Z_G\ones=(\sum_{r=0}^D k_rt^r)\ones$, the magnitude weighting is $\omega=({1}/{\sum_{r=0}^D k_rt^r})\ones$.
The two equalities follow from $\sum_{x,y\in V}t^{d(x,y)}=\abs{V}\sum_{r=0}^D k_rt^r$.
\end{proof}

For a vertex-transitive graph, this recovers Speyer's formula \cite[Lemma~3.2]{Lei19}.
If a regular graph has $v$ vertices, degree $k$, and diameter $2$, then $k_2=v-k-1$. Hence
\begin{equation}\label{eq:valency-constant-diam2}
    \Mag(G;t)=\frac{v}{1+kt+(v-k-1)t^2}.
\end{equation}

In terms of the distance matrices, $Z_G=\sum_{i=0}^Dt^iA_i=I+\sum_{i=1}^Dt^iA_i$. Expanding $Z_G^{-1}$ gives
\begin{align}
    Z_G^{-1}
    &=\sum_{n=0}^{\infty}(-1)^n(Z_G-I)^n \notag\\
    &=I+t(-A)+t^2(A^2-A_2)
      +t^3(-A^3+AA_2+A_2A-A_3)+\cdots
      \label{eq:mobius}
\end{align}
where $A_i$ is the zero matrix $O$ for $i>D$.

Next, we define sets of vertex sequences of a given total length.

\begin{defi}
For $n,\ell\ge0$ and $u,u'\in V$, define
\begin{equation}
    P_{n,\ell}(G;u,u')
    \coloneq
    \left\{(x_0,\ldots,x_n)\in V^{n+1}
    \;\middle|\;
    \begin{aligned}
        &u=x_0\ne x_1\ne\cdots\ne x_n=u',\\
        &\sum_{i=1}^n d(x_{i-1},x_i)=\ell
    \end{aligned}
    \right\}.
\end{equation}
Only consecutive vertices must be distinct.
We set $P_{0,\ell}(G;u,u')$ equal to $\{(u)\}$ if $\ell=0$ and $u=u'$, and to the empty set otherwise. We also set $P_{-1,\ell}(G;u,u')=\emptyset$.
Put
\begin{equation}\label{eq:definition-of-Pnl}
    P_{n,\ell}(G)\coloneq \bigsqcup_{u,u'\in V}P_{n,\ell}(G;u,u').
\end{equation}
\end{defi}

Since $((Z_G-I)^n)_{uu'}=\sum_{\ell=n}^{\infty}\abs{P_{n,\ell}(G;u,u')}t^\ell$ for each $n\ge0$, \eqref{eq:mobius} gives
\begin{equation}\label{eq:entries-mobius}
    (Z_G^{-1})_{uu'}
    =\sum_{\ell=0}^{\infty}
    \left(\sum_{n=0}^{\ell}(-1)^n\abs{P_{n,\ell}(G;u,u')}\right)t^\ell.
\end{equation}
Summing over all entries of $Z_G^{-1}$ gives
\begin{equation}\label{eq:magnitude-series}
    \Mag(G;t)
    =\sum_{u,u'\in V}(Z_G^{-1})_{uu'}
    =\sum_{\ell=0}^{\infty}
    \left(\sum_{n=0}^{\ell}(-1)^n\abs{P_{n,\ell}(G)}\right)t^\ell.
\end{equation}
Thus, every coefficient of the magnitude is an integer.

\subsection{Magnitude homology and distance algebras}\label{subsec:maghom}
\begin{defi}[\cite{HW17}]
The \emph{magnitude chain complex} $\MC_{*,*}(G)$ has chain groups $\MC_{n,\ell}(G)=\mathbb{Z}P_{n,\ell}(G)$ and differential $\partial=\sum_{i=1}^{n-1}(-1)^i\partial_i\colon\MC_{n,\ell}(G)\to\MC_{n-1,\ell}(G)$, where
\begin{equation}
    \partial_i(x_0,\ldots,x_n)=
    \begin{cases}
        (x_0,\ldots,x_{i-1},x_{i+1},\ldots,x_n)
        &\text{if }d(x_{i-1},x_{i+1})
        =d(x_{i-1},x_i)+d(x_i,x_{i+1}),\\
        0 &\text{otherwise}.
    \end{cases}
\end{equation}
We set $\partial=0$ for $n\le1$.
The \emph{magnitude homology} of $G$ is $\MH_{n,\ell}(G)=H_n(\MC_{*,\ell}(G))$.
\end{defi}

The groups $\MC_{n,\ell}(G;u,u')=\mathbb{Z}P_{n,\ell}(G;u,u')$ with fixed endpoints form a subcomplex. We write its homology as $\MH_{n,\ell}(G;u,u')$.
By \eqref{eq:definition-of-Pnl},
\begin{equation}\label{eq:chain-hom-direct-sum}
    \MC_{n,\ell}(G)=\bigoplus_{u,u'\in V}\MC_{n,\ell}(G;u,u'),
    \qquad
    \MH_{n,\ell}(G)=\bigoplus_{u,u'\in V}\MH_{n,\ell}(G;u,u').
\end{equation}
The chain groups vanish in degrees $n>\ell$, so \eqref{eq:entries-mobius} gives the Euler characteristic
\begin{align}
    \chi(\MH_{*,\ell}(G;u,u'))
    =\chi(\MC_{*,\ell}(G;u,u'))
      &=\sum_{n=0}^{\ell}(-1)^n\abs{P_{n,\ell}(G;u,u')}\\
      &=[t^\ell](Z_G^{-1})_{uu'}.
      \label{eq:euler-char-component}
\end{align}
Summing over all pairs of endpoints and using \eqref{eq:magnitude-series} gives
\begin{equation}
    \Mag(G;t)
    =\sum_{\ell=0}^{\infty}t^\ell\chi(\MH_{*,\ell}(G))
    =\sum_{\ell=0}^{\infty}t^\ell\chi(\MC_{*,\ell}(G)).
\end{equation}
In this sense, magnitude homology categorifies magnitude.

The next proposition follows from $\MC_{\ell+1,\ell}(G)=0$.
\begin{prop}\label{prop:freeness-diagonal-part}
The group $\MH_{\ell,\ell}(G)$ is free abelian.
\end{prop}

The off-diagonal groups vanish for $\ell\le2$ \cite[Theorem 4.1, Example 4.6]{LS21}. 
Therefore we obtain the following from \eqref{eq:mobius}, \eqref{eq:euler-char-component}, and \zcref{prop:freeness-diagonal-part}.
\begin{prop}[{\cite[Proposition~2.9]{HW17}, \cite[Corollary~4.6]{IM26}}]\label{prop:low-degree-magnitude-homology}
For any $u,u'\in V$,
\begin{equation}
    \MH_{0,0}(G;u,u')\cong\mathbb{Z}^{\,I_{uu'}},\quad
    \MH_{1,1}(G;u,u')\cong\mathbb{Z}^{\,A_{uu'}},\quad
    \MH_{2,2}(G;u,u')\cong\mathbb{Z}^{\,(A^2-A_2)_{uu'}}.
\end{equation}
\end{prop}

Although the off-diagonal magnitude homology does not vanish in general for $\ell \geq 3$, there often exist graphs for which it vanishes even for $\ell \geq 3$.

\begin{defi}[{\cite[Definition~7.1]{HW17}}]
A graph $G$ is \emph{diagonal} if $\MH_{n,\ell}(G)=0$ for all $n,\ell\in\mathbb{Z}_{\ge0}$ with $n\ne\ell$.
\end{defi}

Suppose that $G$ is diagonal. By \zcref{prop:freeness-diagonal-part} and the Euler characteristic formula, we have $\rank\MH_{\ell,\ell}(G)=(-1)^\ell[t^\ell]\Mag(G;t)$.
By~\eqref{eq:euler-char-component}, we also have
\begin{equation}\label{eq:rank-of-diagonal}
    (-1)^\ell\rank\MH_{\ell,\ell}(G;u,u')
    =\sum_{n=0}^{\ell}(-1)^n\abs{P_{n,\ell}(G;u,u')}
    =[t^\ell](Z_G^{-1})_{uu'}.
\end{equation}
The groups $\MH_{\ell,\ell}(G;u,u')$ are free abelian, so the ranks in \eqref{eq:rank-of-diagonal} determine their isomorphism types.
In particular, the magnitude of a diagonal graph determines its magnitude homology.
If $k_r=\abs{\Gamma_r(u)}$ is also independent of $u$, then \zcref{prop:valency-constant} shows that $\MH_{n,\ell}(G)$ is determined by $k_0,\ldots,k_D$.

\smallskip
\paragraph{\textit{Distance algebras}}
Let $\mathbb{F}$ be a field and let $G=(V,E)$ be a finite nonempty connected graph.
The \emph{distance algebra} $\mathsf{DA}_{\mathbb{F}}(G)$ of $G$ is the $\mathbb{F}$-vector space with basis $V\times V$, equipped with the multiplication
\begin{equation}
    (u,w)(x,y)=
    \begin{cases}
        (u,y) &\text{if } w=x \text{ and } d(u,w)+d(w,y)=d(u,y),\\
        0 &\text{otherwise}.
    \end{cases}
\end{equation}
It is a graded algebra in which $(u,v)$ is homogeneous of degree $d(u,v)$, with unit $\sum_{u\in V}(u,u)$.
For brevity, write $\mathsf{DA}=\mathsf{DA}_{\mathbb{F}}(G)$ and let $\mathsf{DA}_i$ denote its degree-$i$ part.
Then $\mathsf{DA}_0\cong\mathbb{F}^{V}$ is semisimple, and $\mathsf{DA}$ is generated over $\mathsf{DA}_0$ by $\mathsf{DA}_1$.
The algebra $\mathsf{DA}$ is called \emph{quadratic} if the kernel of the canonical surjection from the tensor algebra of the $\mathsf{DA}_0$-bimodule $\mathsf{DA}_1$ onto $\mathsf{DA}$ is generated as a two-sided ideal by its degree-two part.
It is called \emph{Koszul} if $\mathsf{DA}_0$, regarded as a graded left $\mathsf{DA}$-module via the projection $\mathsf{DA}\to\mathsf{DA}_0$, admits a graded projective resolution whose $n$-th term is generated in degree $n$.
See \cite{AI24} and \cite[Sections~3.3 and~6.1]{IM26}.

By \cite[Proposition~6.4 and Theorem~6.5]{IM26}, the distance algebra $\mathsf{DA}_{\mathbb{F}}(G)$ is quadratic over every field $\mathbb{F}$ if and only if $\MH_{2,\ell}(G)=0$ for all $\ell>2$, and it is Koszul over every field $\mathbb{F}$ if and only if $G$ is diagonal.

\smallskip
\paragraph{\textit{Examples of diagonal graphs}}

Hepworth and Willerton showed that complete graphs and trees are diagonal~\cite[Example~2.5 and Corollary~6.8]{HW17}, and they also showed that diagonality is preserved under Cartesian products~\cite[Proposition~7.3]{HW17}.
For example, Hamming graphs are obtained as Cartesian products of complete graphs (and are therefore diagonal).

\begin{defi}[{\cite[Section~12.4.1]{BH2012}}]
Let $D\ge1$ and $m\ge2$ be integers.
The \emph{Hamming graph} $H(D,m)$ has vertex set $[m]^D$, with two vertices adjacent if they differ in exactly one coordinate.
The graph $H(D,2)$ is the \emph{hypercube} $Q_D$.
\end{defi}

The graph $H(D,m)$ is diagonal because it is the Cartesian product of $D$ copies of $K_m$. A \emph{grid graph}, defined as a Cartesian product of finitely many paths, is diagonal for the same reason.

The join of any two nonempty graphs is diagonal as well \cite[Theorem~7.5]{HW17}.
Hence, complete multipartite graphs \cite[Example~7.6]{HW17}, including complete bipartite graphs $K_{m,n}$ and \emph{hyperoctahedral graphs} $K_{m\times2}$, are diagonal.
A \emph{wheel graph} is obtained by adding one vertex to a cycle and joining the new vertex to every vertex of the cycle. The resulting graph is a join and is therefore diagonal.

A \emph{block graph} is obtained by successively gluing complete graphs along single vertices. The homology decomposition for such gluings \cite[Corollary~6.9]{HW17} shows that block graphs are diagonal (see also \cite[Appendix~B]{Gu18}).

In addition, Gu showed that so-called pawful graphs are diagonal.

\begin{defi}[{\cite[Definition~4.2]{Gu18}}]
A connected graph $G=(V,E)$ is \emph{pawful} if $\diam G\le2$ and $\Gamma(x)\cap\Gamma(y)\cap\Gamma(z)\ne\emptyset$ whenever $d(x,y)=d(y,z)=2$ and $d(x,z)=1$.
\end{defi}

For example, complete graphs and complete multipartite graphs are pawful.
Gu used algebraic Morse theory to show that pawful graphs are diagonal \cite[Theorem~4.4]{Gu18}.

\begin{defi}[{\cite[Section~2.2.2.1]{CCHO20}}]
A vertex $w$ is a \emph{median} of $x_1,x_2,x_3\in V$ if $d(x_i,x_j)=d(x_i,w)+d(w,x_j)$ for any distinct $i,j\in [3]$.
A connected graph is \emph{median} if every triple of vertices has a unique median.
\end{defi}

Bottinelli and Kaiser proved that median graphs are diagonal \cite{BK21}.
Trees are median, and the Cartesian product of median graphs is median. In particular, hypercubes and grid graphs are median.

The \emph{icosahedral graph}, the $1$-skeleton of the regular icosahedron, is a diagonal graph that is neither pawful nor median \cite[Theorem~4.5]{Gu18}.

\section{Reducive graphs}\label{sec:reducive}

In this section, we introduce reducivity as a combinatorial condition for the diagonality of graphs. Using this condition, we can derive the diagonality of various graphs.

\subsection{The main theorem}
\begin{defi}\label{def:target}
For a vertex $u\in V$ of a graph $G=(V,E)$ and a set $S\subseteq\Gamma(u)$, define the target set from $u$ through $S$, denoted by $T_u(S)$, as follows:
\begin{equation}
        T_u(S)\coloneq \left\{x\in V\mathrel{}\middle|\mathrel{} d(u,x)=1+d(w,x)\text{ for some }w\in S\right\}.
\end{equation}
For $R\subseteq\Gamma(u)$, put
\begin{equation}
        T^R_u(S)\coloneq T_u(S)\setminus T_u(R).
\end{equation}
\end{defi}

The target set $T_u(S)$ consists of the vertices that can be reached from $u$ by a shortest path whose first step is to a vertex of $S$.
Note that we have $T_u(\emptyset)=\emptyset$.
The definition gives the following properties.

\begin{prop}\label{prop:target-union}
    Let $u\in V$ be a vertex of a graph $G=(V,E)$. For $S,S'\subseteq\Gamma(u)$, we have $T_u(S\cup S')=T_u(S)\cup T_u(S')$.
\end{prop}

\begin{prop}\label{prop:target-neighbor}
    For any $u\in V$ and $S\subseteq \Gamma(u)$, $S=T_u(S)\cap \Gamma(u)$ holds.
\end{prop}

\begin{defi}[Reducible condition]\label{df:reducible-cond}
Let $u \in V$, $w \in \Gamma(u)$, $S \subseteq \Gamma(u) \setminus w$, and $R \subseteq \Gamma(w)$.
We say that $u$ is $(S,R)$-\emph{reducible} to $w$ if
\begin{equation}\label{eq:reducible-cond}
    V=T_u^S(w) \sqcup T_w(R).
\end{equation}
\end{defi}

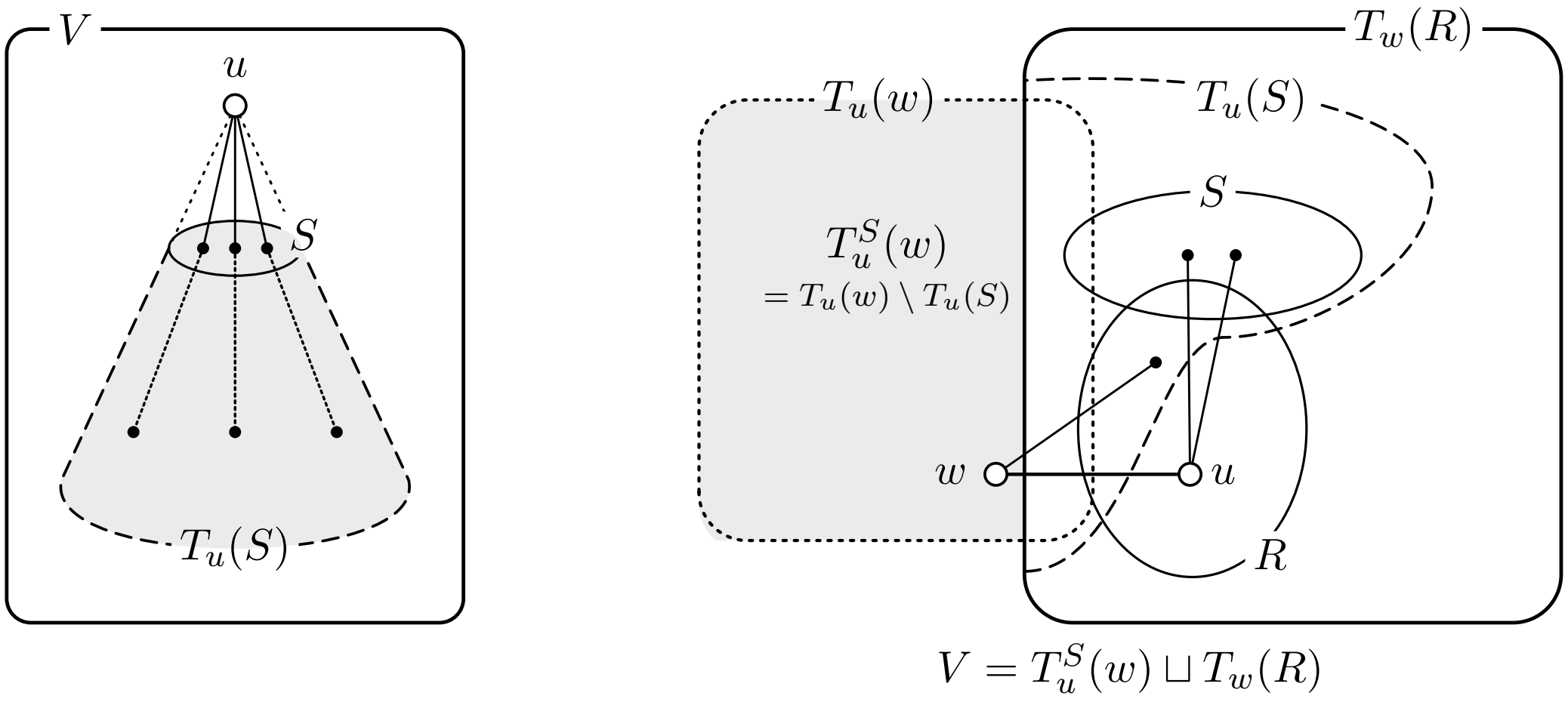
\begin{figure}[htbp]
\centering
\tikzset{
    redvertex/.style={circle,draw,fill=white,inner sep=0pt,minimum size=5.5pt,line width=.7pt},
    redsample/.style={circle,fill=black,inner sep=0pt,minimum size=3pt},
    redlocal/.style={line width=.6pt},
    redtargetS/.style={line width=.7pt,dash pattern=on 4pt off 2.5pt},
    redtargetw/.style={line width=.8pt,dotted},
    redpartition/.style={line width=.9pt},
    rededge/.style={line width=.55pt},
    redextension/.style={line width=.65pt,densely dotted},
    redlabel/.style={inner sep=2pt,outer sep=0pt}
}
\begin{minipage}[t]{.31\linewidth}
    \centering
    \begin{tikzpicture}[x=1cm,y=1cm,line cap=round,line join=round,
        every node/.style={font=\small},baseline=(current bounding box.north)]
        \coordinate (a) at (2,4.53);
        \coordinate (sl) at (1.42,3.28);
        \coordinate (sr) at (2.58,3.28);
        \coordinate (a1) at (1.72,3.28);
        \coordinate (a2) at (2,3.28);
        \coordinate (a3) at (2.28,3.28);
        \coordinate (b1) at (1.11,1.67);
        \coordinate (b2) at (2,1.67);
        \coordinate (b3) at (2.89,1.67);
        \path[fill=black!8]
            (sl) arc[start angle=180,end angle=0,x radius=.58cm,y radius=.24cm]
            -- (3.52,1.25)
            .. controls (3.70,.45) and (.30,.45) .. (.48,1.25) -- cycle;
        \node[redlabel] (labelV) at (.60,5.20) {$V$};
        \node[redlabel] (labelS) at (2.60,3.4) {$S$};
        \node[redlabel] (labelTuS) at (2,.65) {$T_u(S)$};
        \begin{scope}
            \pgfseteorule
            \clip (-.25,-.25) rectangle (4.25,5.55)
                (labelV.south west) rectangle (labelV.north east);
            \draw[redpartition,rounded corners=6pt] (0,0) rectangle (4,5.2);
        \end{scope}
        \begin{scope}
            \pgfseteorule
            \clip (-.25,-.25) rectangle (4.25,5.55)
                (labelS.south west) rectangle (labelS.north east)
                (labelTuS.south west) rectangle (labelTuS.north east);
            \draw[redlocal,dotted] (a)--(sl) (a)--(sr);
            \draw[redtargetS] (sl)--(.48,1.25)
                .. controls (.30,.45) and (3.70,.45) .. (3.52,1.25)--(sr);
            \draw[redlocal] (2,3.28) ellipse[x radius=.58,y radius=.24];
        \end{scope}
        \draw[rededge] (a)--(a1) (a)--(a2) (a)--(a3);
        \draw[redextension] (a1)--(b1) (a2)--(b2) (a3)--(b3);
        \node[redvertex,label=above:$u$] at (a) {};
        \foreach \p in {a1,a2,a3,b1,b2,b3} \node[redsample] at (\p) {};
        \path[use as bounding box] (2,-.75)--(2,5.55);
    \end{tikzpicture}\par
    \smallskip
    {\small\textup{(a)} Target set\par}
\end{minipage}\hspace{.025\linewidth}
\begin{minipage}[t]{.59\linewidth}
    \centering
    \begin{tikzpicture}[x=1cm,y=1cm,line cap=round,line join=round,
        every node/.style={font=\small},baseline=(current bounding box.north)]
        \begin{scope}
            \clip[rounded corners=12pt] (.20,.62) rectangle (3.65,4.58);
            \fill[black!8] (.20,.72) rectangle (3.05,4.58);
        \end{scope}
        \coordinate (u) at (4.5,1.3);
        \coordinate (w) at (2.8,1.3);
        \coordinate (s1) at (4.90,3.22);
        \coordinate (s2) at (4.48,3.22);
        \coordinate (r1) at (4.2,2.28);
        \node[redlabel] (labelTuw) at (1.78,4.58) {$T_u(w)$};
        \node[redlabel] (labelTwR) at (6.46,5.20) {$T_w(R)$};
        \node[redlabel] (labelTuS) at (5.04,4.58) {$T_u(S)$};
        \node[redlabel] (labelS) at (4.70,3.78) {$S$};
        \node[redlabel] (labelR) at (5.2,.6) {$R$};
        \begin{scope}
            \pgfseteorule
            \clip (-.25,-.25) rectangle (8,5.55)
                (labelTuw.south west) rectangle (labelTuw.north east);
            \draw[redtargetw,rounded corners=12pt] (.20,.72) rectangle (3.65,4.58);
        \end{scope}
        \begin{scope}
            \pgfseteorule
            \clip (-.25,-.25) rectangle (8,5.55)
                (labelTwR.south west) rectangle (labelTwR.north east);
            \draw[redpartition,rounded corners=12pt] (3.05,0) rectangle (7.75,5.2);
        \end{scope}
        \begin{scope}
            \pgfseteorule
            \clip (-.25,-.25) rectangle (8,5.55)
                (labelTuS.south west) rectangle (labelTuS.north east);
            \draw[redtargetS] (3.05,.45)
                .. controls (4.0,0.45) and (4.2,2.5) .. (4.8,2.5)
                .. controls (5.5,2.5) and (6.62,3.1) .. (6.62,3.82)
                .. controls (6.62,4.50) and (4.78,4.85) .. (3.05,4.75);
        \end{scope}
        \begin{scope}
            \pgfseteorule
            \clip (-.25,-.25) rectangle (8,5.55)
                (labelS.south west) rectangle (labelS.north east);
            \draw[redlocal] (4.70,3.22) ellipse[x radius=1.30,y radius=.56];
        \end{scope}
        \begin{scope}
            \pgfseteorule
            \clip (-.25,-.25) rectangle (8,5.55)
                (labelR.south west) rectangle (labelR.north east);
            \draw[redlocal] (4.52,1.70) ellipse[x radius=1.0,y radius=1.30];
        \end{scope}
        \draw[rededge] (u)--(s1) (u)--(s2) (w)--(r1);
        \draw[line width=.85pt] (u)--(w);
        \node[align=center] at (1.85,3.12)
            {$T_u^S(w)$\\[-1pt]{\tiny $=T_u(w)\setminus T_u(S)$}};
        \node[redvertex] at (u) {};
        \node[right=1pt] at (u) {$u$};
        \node[redvertex,label=left:$w$] at (w) {};
        \foreach \p in {s1,s2,r1} \node[redsample] at (\p) {};
        \node[font=\small] at (3.975,-.40)
            {$V=T_u^S(w)\sqcup T_w(R)$};
        \path[use as bounding box] (2,-.75)--(2,5.55);
    \end{tikzpicture}\par
    \smallskip
    {\small\textup{(b)} Reducible condition\par}
\end{minipage}\par

\smallskip
\begin{tikzpicture}[baseline=(legend.base),every node/.style={font=\small}]
    \node[draw,fill=black!8,minimum width=7mm,minimum height=3mm,inner sep=0pt,line width=.5pt] (swatch) at (0,0) {};
    \node[anchor=west] (legend) at (.48,0)
        {Shaded regions: (a) $T_u(S)$; \quad (b) $T_u^S(w)$.};
\end{tikzpicture}
\caption{A target set and the reducible condition.}
\label{fig:Reducible}
\end{figure}

From \zcref{prop:target-neighbor}, we see that $R$ is subject to several conditions, such as $u \in R$ and $S \cap \Gamma(w) \subset R$. A discussion of these and related conditions on $R$ will be given in \zcref{sec:reducive-criteria}. 
\zcref{fig:Reducible} illustrates a target set and the reducible condition.

\begin{defi}\label{def:generalized-TRYY}\label{def:family}
A graph $G=(V,E)$ is \emph{reducive} if there is a family $(\mathcal{R}_u)_{u\in V}$, with $\mathcal{R}_u\subseteq2^{\Gamma(u)}$ for every $u\in V$, satisfying both of the following conditions at each vertex:
{\setlength{\leftmargini}{44pt}
\begin{itemize}
\labelitem{R1}{item:canonical-condition}   $\Gamma(u)\in\mathcal{R}_u$.
\labelitem{R2}{item:reducible-condition} For every nonempty $S\in\mathcal{R}_u$, there exist $w\in S$ and $R\in\mathcal{R}_w$ such that $S\setminus w\in \mathcal{R}_u$ and $u$ is $(S\setminus w,R)$-reducible to $w$.
\end{itemize}}
\end{defi}

Condition \ref{item:reducible-condition} is designed for an induction on the length degree and the cardinality of $S$. 
More precisely, it allows us to reduce the homology of the quotient complex generated by $\{(u,x_1,\ldots, x_n)\mid x_1\in T_u^{S\setminus w}(w)\}$ to that of the subcomplex generated by $\{(w,x_2,\ldots, x_n)\mid x_2\in T_w(R)\}$.

We first illustrate the reducive condition with some elementary examples.
Complete graphs and complete bipartite graphs are reducive, whereas the $5$-cycle is not.
These examples illustrate how the partition condition can be verified or fail.

\begin{exam}[Complete graphs]\label{exam:complete-graph-reducive}
The complete graph $K_m$ is reducive for every $m\ge1$.
Write $V=V(K_m)$ and take $\mathcal R_u=2^{\Gamma(u)}$ for every $u\in V$.
Then \ref{item:canonical-condition} holds. Every shortest path between distinct vertices
has length $1$, so $T_u(S)=S$ for every $S\subseteq\Gamma(u)$.
If $S\ne\emptyset$, choose any $w\in S$ and set $R=\Gamma(w)$. Then
$T_u^{S\setminus w}(w)=\{w\}$ and $T_w(R)=V\setminus w$.
The two sets partition $V$, giving \ref{item:reducible-condition}.
For $m=1$ there is no nonempty $S$ to check.
\end{exam}

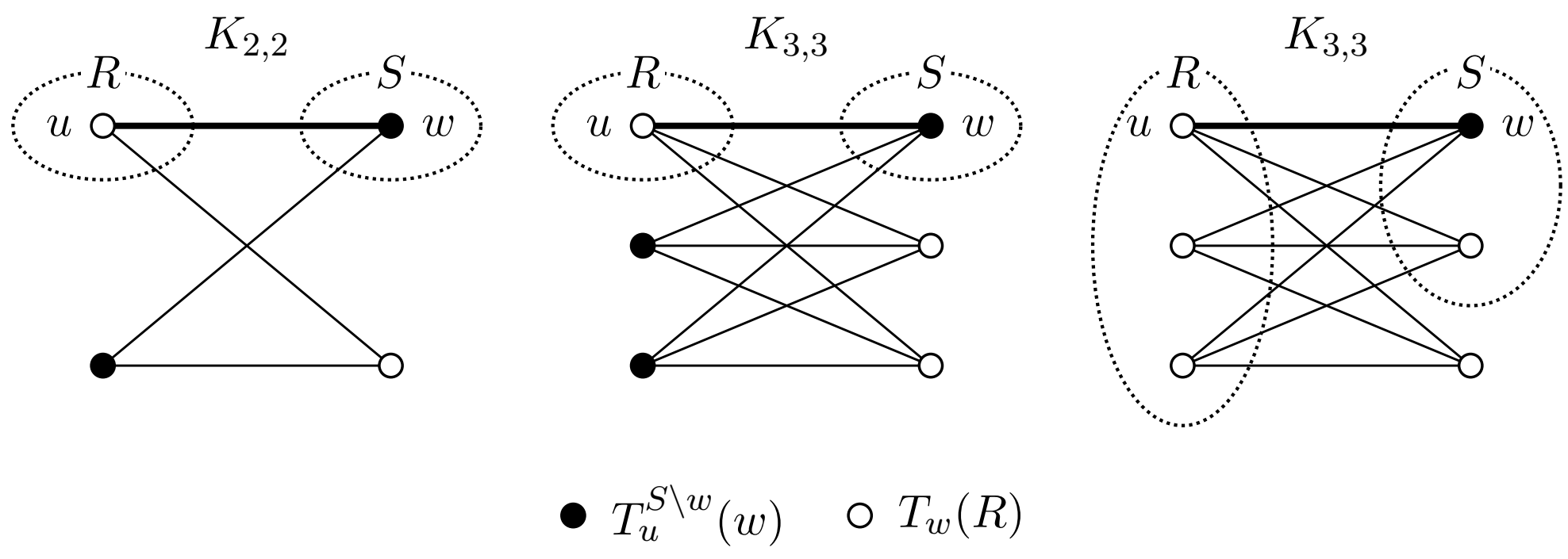
\begin{figure}[htbp]
\centering
\begin{tikzpicture}[
    font=\small, line width=.55pt,
    vertex/.style={circle,draw=black,fill=white,inner sep=0pt,
        minimum size=5.5pt,line width=.65pt},
    filled vertex/.style={vertex,fill=black},
    set boundary/.style={densely dotted,line width=.8pt},
    set label/.style={fill=white,inner sep=1pt}
]
    \begin{scope}[xshift=-4.5cm]
        \node at (0,1.72) {$K_{2,2}$};
        \coordinate (left-u) at (-1.20,1.0);
        \coordinate (left-w) at (1.20,1.0);
        \coordinate (left-u1) at (-1.20,-1.0);
        \coordinate (left-w1) at (1.20,-1.0);

        \draw[line width=1.5pt] (left-u)--(left-w);
        \draw (left-u)--(left-w1);
        \draw (left-u1)--(left-w);
        \draw (left-u1)--(left-w1);

        \draw[set boundary] (-1.20,1.0)
            ellipse [x radius=.75cm,y radius=.45cm];
        \draw[set boundary] (1.20,1.0)
            ellipse [x radius=.75cm,y radius=.45cm];
        \node[set label] at (-1.2,1.45) {$R$};
        \node[set label] at (1.2,1.45) {$S$};

        \node[vertex,label=left:{$u$}] at (left-u) {};
        \node[filled vertex,label=right:{$w$}] at (left-w) {};
        \node[filled vertex] at (left-u1) {};
        \node[vertex] at (left-w1) {};
    \end{scope}
    
    \begin{scope}[xshift=0cm]
        \node at (0,1.72) {$K_{3,3}$};
        \coordinate (right-u) at (-1.20,1.0);
        \coordinate (right-w) at (1.20,1.0);
        \coordinate (right-u2) at (-1.20,0.0);
        \coordinate (right-w2) at (1.20,0.0);
        \coordinate (right-u3) at (-1.20,-1.0);
        \coordinate (right-w3) at (1.20,-1.0);

        \draw[line width=1.5pt] (right-u)--(right-w);
        \draw (right-u)--(right-w2);
        \draw (right-u)--(right-w3);
        \draw (right-u2)--(right-w);
        \draw (right-u2)--(right-w2);
        \draw (right-u2)--(right-w3);
        \draw (right-u3)--(right-w);
        \draw (right-u3)--(right-w2);
        \draw (right-u3)--(right-w3);

        \draw[set boundary] (-1.20,1.0)
            ellipse [x radius=.75cm,y radius=.45cm];
        \draw[set boundary] (1.20,1.0)
            ellipse [x radius=.75cm,y radius=.45cm];
        \node[set label] at (-1.2,1.45) {$R$};
        \node[set label] at (1.2,1.45) {$S$};

        \node[vertex,label=left:{$u$}] at (right-u) {};
        \node[filled vertex,label=right:{$w$}] at (right-w) {};
        \node[filled vertex] at (right-u2) {};
        \node[vertex] at (right-w2) {};
        \node[filled vertex] at (right-u3) {};
        \node[vertex] at (right-w3) {};
    \end{scope}
    
    \begin{scope}[xshift=4.5cm]
        \node at (0,1.72) {$K_{3,3}$};
        \coordinate (right-u) at (-1.20,1.0);
        \coordinate (right-w) at (1.20,1.0);
        \coordinate (right-u2) at (-1.20,0.0);
        \coordinate (right-w2) at (1.20,0.0);
        \coordinate (right-u3) at (-1.20,-1.0);
        \coordinate (right-w3) at (1.20,-1.0);

        \draw[line width=1.5pt] (right-u)--(right-w);
        \draw (right-u)--(right-w2);
        \draw (right-u)--(right-w3);
        \draw (right-u2)--(right-w);
        \draw (right-u2)--(right-w2);
        \draw (right-u2)--(right-w3);
        \draw (right-u3)--(right-w);
        \draw (right-u3)--(right-w2);
        \draw (right-u3)--(right-w3);

        \draw[set boundary] (-1.20,0)
            ellipse [x radius=.75cm,y radius=1.50cm];
        \draw[set boundary] (1.20,.50)
            ellipse [x radius=.75cm,y radius=1.00cm];
        \node[set label] at (-1.2,1.45) {$R$};
        \node[set label] at (1.2,1.45) {$S$};

        \node[vertex,label=left:{$u$}] at (right-u) {};
        \node[filled vertex,label=right:{$w$}] at (right-w) {};
        \node[vertex] at (right-u2) {};
        \node[vertex] at (right-w2) {};
        \node[vertex] at (right-u3) {};
        \node[vertex] at (right-w3) {};
    \end{scope}

    \node[filled vertex] at (-1.78,-2.25) {};
    \node[anchor=west] at (-1.61,-2.25)
        {$T_u^{S\setminus w}(w)$};
    \node[vertex] at (.61,-2.25) {};
    \node[anchor=west] at (.78,-2.25) {$T_w(R)$};
\end{tikzpicture}
\caption{Reductions in bipartite graphs.}
\label{fig:bipartite-reduction-examples}
\end{figure}

\begin{exam}[Complete bipartite graphs]\label{exam:complete-bipartite-reducive}
The complete bipartite graph $K_{m_1,m_2}$ is reducive for all $m_1,m_2\ge1$.
Write $V=V(K_{m_1,m_2})$ and take $\mathcal R_u=2^{\Gamma(u)}$ for every $u\in V$.
This gives \ref{item:canonical-condition}. For fixed $u$, the two parts are
$\{u\}\sqcup\Gamma_2(u)$ and $\Gamma(u)$. For every nonempty $S\subseteq\Gamma(u)$,
$T_u(S)=S\sqcup\Gamma_2(u)$.
Indeed, each vertex of $\Gamma_2(u)$ can be reached from $u$ in two steps
through any vertex of $S$. 

Choose any $w\in S$. If $S=\{w\}$, take $R=\{u\}$. Then
$T_u^{S\setminus w}(w)=\{w\}\sqcup\Gamma_2(u)$ and $T_w(R)=\{u\}\sqcup(\Gamma(u)\setminus w)$.
If $\abs{S}\ge2$, then $T_u(S\setminus w)$ contains $\Gamma_2(u)$,
so $T_u^{S\setminus w}(w)=\{w\}$. Take $R=\Gamma(w)$. Then
$T_w(R)=V\setminus w$. Both cases give the partition required by \ref{item:reducible-condition}
(see Figure~\ref{fig:bipartite-reduction-examples}).
\end{exam}

\begin{figure}[htbp]
\centering
\begin{tikzpicture}[
  vertex/.style={circle,draw,fill=white,minimum size=2.5mm,inner sep=0pt},
  marked/.style={vertex,fill=black},
  line width=.55pt,font=\small]
  \foreach \panel/\shift in {1/0,2/4.9}{
    \begin{scope}[xshift=\shift cm]
      \node[vertex,label={90:$u$}] (u) at (90:1.18) {};
      \node[vertex,label={162:$z$}] (z) at (162:1.18) {};
      \node[marked,label={-54:$x$}] (x) at (-54:1.18) {};
      \ifnum\panel=1
        \node[marked,label={18:$w$}] (w) at (18:1.18) {};
        \node[vertex,label={-126:$y$}] (y) at (-126:1.18) {};
      \else
        \node[vertex,label={18:$w$}] (w) at (18:1.18) {};
        \node[marked,label={-126:$y$}] (y) at (-126:1.18) {};
      \fi
      \draw (u)--(w)--(x)--(y)--(z)--(u);
      \ifnum\panel=1
        \node at (0,-1.83) {$T_u^{\{z\}}(w)=\{w,x\}$};
      \else
        \draw[line width=1.45pt] (w)--(x)--(y);
        \node at (0,-1.83) {$T_w(x)=\{x,y\}$};
      \fi
    \end{scope}
  }
\end{tikzpicture}
\caption{Black vertices form the indicated sets. The bold path is the unique shortest path from $w$ to $y$.}
\label{fig:c5-nonreducive}
\end{figure}
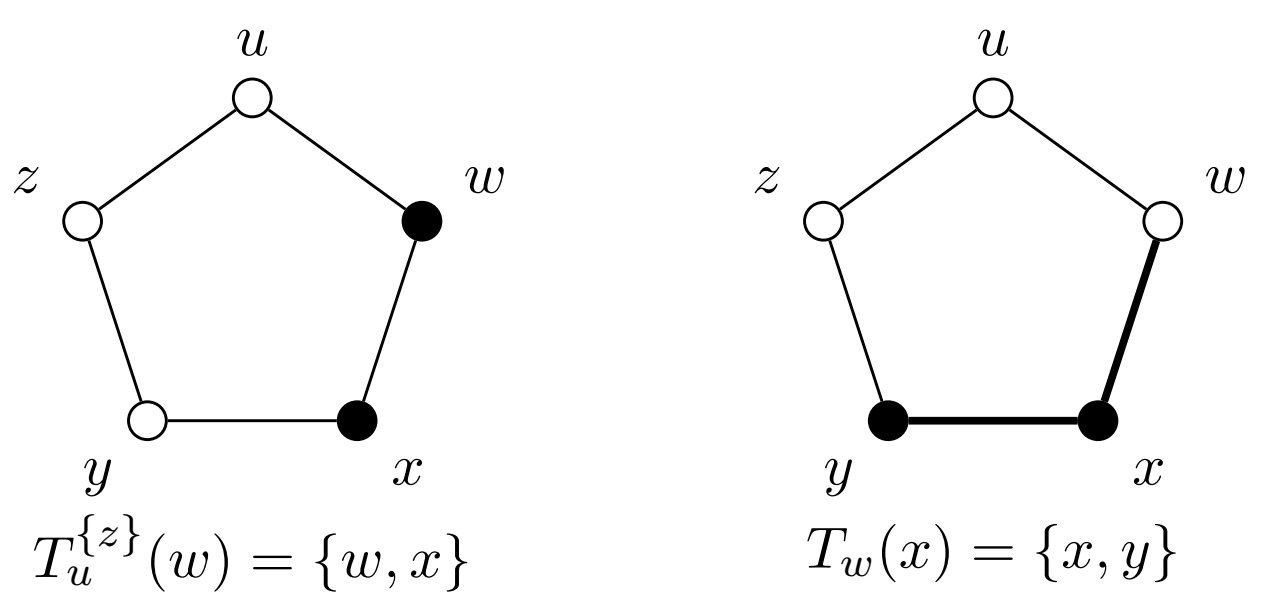

\begin{exam}[The $5$-cycle]\label{exam:c5-nonreducive}
The $5$-cycle $C_5$ is not reducive.
Label its vertices $u,w,x,y,z$ in cyclic order, as in Figure~\ref{fig:c5-nonreducive}.
By \ref{item:canonical-condition}, $S=\Gamma(u)=\{w,z\}$ must belong to $\mathcal R_u$.
Reflection fixing $u$ interchanges $w$ and $z$, so it suffices to rule out a reduction from $u$ to $w$.
We have $T_u(w)=\{w,x\}$ and $T_u(z)=\{z,y\}$, hence $T_u^{\{z\}}(w)=\{w,x\}$.
The required partition would give $T_w(R)=\{u,y,z\}$ for some $R\subseteq\Gamma(w)$.
The unique shortest path from $w$ to $y$ is $(w,x,y)$, so $y\in T_w(R)$ forces $x\in R$.
Then $x\in T_w(R)$ as well, contradicting the partition.
\end{exam}

Like diagonal graphs and median graphs, reducive graphs are closed under Cartesian products.

\begin{prop}\label{prop:reducive-cartesian-product}
If $G$ and $G'$ are reducive, then $G\square G'$ is reducive.
\end{prop}
\begin{proof}
Write $V=V(G)$ and $V'=V(G')$. Choose families $(\mathcal R_u)_{u\in V}$
and $(\mathcal R'_{u'})_{u'\in V'}$ satisfying \ref{item:canonical-condition}
and \ref{item:reducible-condition} in $G$ and $G'$, respectively.
Use $\Gamma'$ and $T'$ for the neighborhoods and target sets in $G'$.
For each $(u,u')\in V\times V'$, set
\begin{equation}
    \mathcal{R}^\square_{(u,u')}
       =\bigl\{(S\times\{u'\})\cup(\{u\}\times S')
          \mid S\in\mathcal{R}_u,\ S'\in\mathcal{R}'_{u'}\bigr\}.
\end{equation}
Taking $S=\Gamma(u)$
and $S'=\Gamma'(u')$ gives the full neighborhood of $(u,u')$,
so \ref{item:canonical-condition} holds.

By definition, the distance between $(u,u')$ and $(v,v')$ is given by $d_G(u,v)+d_{G'}(u',v')$.
For $S\subseteq\Gamma(u)$ and $S'\subseteq\Gamma'(u')$, write
$S\square S'=(S\times\{u'\})\cup(\{u\}\times S')$. Then we have
\begin{equation}\label{eq:cartesian-reduction-target}
    T_{(u,u')}(S\square S')
       =(T_u(S)\times V')\cup(V\times T'_{u'}(S')).
\end{equation}

Fix a nonempty $S\square S'\in\mathcal R^\square_{(u,u')}$.
Interchanging the factors if necessary, assume $S\ne\emptyset$.
By \ref{item:reducible-condition} in $G$, choose $w\in S$ and $R\in\mathcal R_w$
such that $S\setminus w\in\mathcal R_u$ and
$T_u^{S\setminus w}(w)=V\setminus T_w(R)$.
Put $w^\square=(w,u')\in S\square S'$. At $w^\square$ take
$R^\square=(R\times\{u'\})\cup(\{w\}\times S') \in\mathcal R^\square_{(w,u')}$.
We also have
$S\square S'\setminus w^\square\in\mathcal R^\square_{(u,u')}$, since deleting $w^\square$ replaces $S$ by $S\setminus w$ in the first coordinate.
Applying \eqref{eq:cartesian-reduction-target} at $(u,u')$ and at $w^\square$ gives
\begin{align}
    T_{(u,u')}^{S\square S'\setminus w^\square}(w^\square)
       &=T_u^{S\setminus w}(w)\times\bigl(V'\setminus T'_{u'}(S')\bigr)\\
       &=(V\setminus T_w(R))\times\bigl(V'\setminus T'_{u'}(S')\bigr)
       =(V\times V')\setminus T_{w^\square}(R^\square).
\end{align}
This proves \ref{item:reducible-condition}.
\end{proof}

We can now state the main theorem.

\begin{thm}\label{thm:reducible}
Reducive graphs are diagonal.
\end{thm}

\subsection{Reduction of chain complexes}\label{subsec:chain-isom}
The proof of \zcref{thm:reducible} is based on a reduction of magnitude chain complexes.
For a fixed initial vertex $u$, we consider subcomplexes generated by vertex sequences whose second entry belongs to a target set $T_u(S)$.
For a neighbor $w\notin S$ of $u$, the inclusion $T_u(S)\subseteq T_u(S\cup\{w\})$ gives a short exact sequence of chain complexes.
When the reducible condition holds, the quotient complex is chain homotopy equivalent to a corresponding subcomplex based at $w$, with both the homological and length degrees decreased by one.
The resulting long exact sequence in homology allows us to prove diagonality by induction on the length degree and the cardinality of $S$.

\begin{defi}
For $n\ge1$ and $X\subseteq V$, set
\begin{equation}
P_{n,\ell}(G;uX,u')=\{(u=x_0,x_1,\ldots,x_n=u')\in P_{n,\ell}(G;u,u')\mid x_1\in X\}.
\end{equation}
For $n\le0$, set $P_{n,\ell}(G;uX,u')=\emptyset$.
Also put $\MC_{n,\ell}(G;uX,u')\coloneq\mathbb{Z}P_{n,\ell}(G;uX,u')$.
\end{defi}

\begin{prop}
    For $S\subseteq\Gamma(u)$, $\MC_{*,\ell}(G;uT_u(S),u')$ is a subcomplex of $\MC_{*,\ell}(G;u,u')$.
\end{prop}
\begin{proof}
    It suffices to show that $\partial_i(x_0,x_1,\dots,x_n)\in\MC_{n-1,\ell}(G;uT_u(S),u')$ for every basis element $p=(x_0,x_1,\dots,x_n)\in\MC_{n,\ell}(G;uT_u(S),u')$ and every $i\in\{1,\ldots,n-1\}$.

    If $i\neq1$, then $\partial_i p\in\MC_{n-1,\ell}(G;uT_u(S),u')$, since the second entry of a nonzero face is unchanged.

    It remains to consider $i=1$.
    If $\partial_1p=0$, then $\partial_1p\in \MC_{n-1,\ell}(G;uT_u(S),u')$.
    Suppose $\partial_1p\neq0$.
    Then $\partial_1p=(x_0,x_2,\ldots,x_n)$ and $d(x_0,x_2)=d(x_0,x_1)+d(x_1,x_2)$.
    Since $x_1\in T_u(S)$, there is a vertex $y\in S$ such that $d(x_0,x_1)=d(x_0,y)+d(y,x_1)$.
    The triangle inequality gives
    \begin{equation}
        d(x_0,x_2)=d(x_0,y)+d(y,x_1)+d(x_1,x_2)\ge d(x_0,y)+d(y,x_2)\ge d(x_0,x_2)
    \end{equation}
    and hence $x_2\in T_u(S)$.
    This proves $\partial_1p\in\MC_{n-1,\ell}(G;uT_u(S),u')$.
\end{proof}
Write $\MH_{n,\ell}(G;uT_u(S),u')$ for the homology in degree $n$ of the subcomplex $\MC_{*,\ell}(G;uT_u(S),u')$.
Fix $u,u'\in V$, $S\subseteq\Gamma(u)$, $w\in\Gamma(u)\setminus S$, and $\ell\ge0$. Put
\begin{equation}
    \mathcal{Q}_{n,\ell}(G;uT_u^S(w),u')\coloneq\mathbb{Z}P_{n,\ell}(G;uT^S_u(w),u').
\end{equation}

This graded group is contained in $\MC_{*,\ell}(G;uT_u(S\cup\{w\}),u')$, but it need not form a subcomplex.
We define its differential as follows.
The maps $\overline{\partial}_i$ defined below determine the differential on this graded group.

\begin{prop}\label{prop:pre}
    For $n\ge 1$, define $\overline{\partial}_i$ by
      \begin{equation}
    \overline{\partial}_i(u,x_1,\ldots,x_n)\coloneq
        \begin{cases*}
            0 & if $\partial_i(u,x_1,\ldots,x_n)\in \MC_{*,\ell}(G;uT_u(S),u')$,\\
            \partial_i(u,x_1,\ldots,x_n) & otherwise.
        \end{cases*}
      \end{equation}
  Then the map
    \begin{align}\label{eq:diff-of-Q}
        \overline{\partial}\coloneq \sum_{i=1}^{n-1}(-1)^{i}\overline{\partial}_i:\mathcal{Q}_{n,\ell}(G;uT_u^S(w),u')\to \mathcal{Q}_{n-1,\ell}(G;uT_u^S(w),u')
    \end{align}
  is well-defined. In particular, $(\mathcal{Q}_{*,\ell}(G;uT_u^S(w),u'),\overline{\partial})$ is a chain complex.
\end{prop}

\begin{proof}
    Take a basis element $p=(u,x_1,\ldots,x_n)\in \mathcal{Q}_{n,\ell}(G;uT_u^S(w),u')$.
    If $i\neq1$ and $\partial_i p\neq0$, the second entry of $\partial_i p$ remains $x_1\in T_u^S(w)$, so $\overline{\partial}_i p=\partial_i p\in\mathcal{Q}_{n-1,\ell}(G;uT_u^S(w),u')$.

    Suppose $i=1$ and $\partial_1 p\neq0$.
    Since $x_1\in T_u^S(w)\subseteq T_u(w)$, we have $d(u,x_1)=1+d(w,x_1)$. Also, $\partial_1 p\neq0$ gives $d(u,x_2)=d(u,x_1)+d(x_1,x_2)$.
    The triangle inequality therefore gives
    \begin{equation}
        d(u,x_2)=1+d(w,x_1)+d(x_1,x_2)\geq1+d(w,x_2)\geq d(u,x_2),
    \end{equation}
    and hence $x_2\in T_u(w)$.
    If $x_2\in T_u(S)$, then $\overline{\partial}_1 p=0$. Otherwise, $x_2\in T_u^S(w)$ and hence $\overline{\partial}_1 p=\partial_1 p\in\mathcal{Q}_{n-1,\ell}(G;uT_u^S(w),u')$.
    Thus $\overline{\partial}$ is well-defined.

    Furthermore, it is clear that $\overline{\partial}^2=0$, and thus $(\mathcal{Q}_{*,\ell}(G;uT_u^S(w),u'),\overline{\partial})$ is a chain complex.
\end{proof}

\begin{prop}\label{prop:chain-isom}
The inclusion of $\MC_{*,\ell}(G;uT_u(S),u')$ into $\MC_{*,\ell}(G;uT_u(S\cup w),u')$ is a chain map.
Moreover, there is an isomorphism
\begin{equation}\label{eq:quotient-complex}
    \mathcal{Q}_{*,\ell}(G;uT_u^S(w),u')\cong\frac{\MC_{*,\ell}(G;uT_u(S\cup w),u')}{\MC_{*,\ell}(G;uT_u(S),u')}.
\end{equation}
\end{prop}

\begin{proof}
    Since $T_u(S)\subseteq T_u(S\cup w)$, the complex $\MC_{*,\ell}(G;uT_u(S),u')$ is a subcomplex of $\MC_{*,\ell}(G;uT_u(S\cup w),u')$.
    The groups $\MC_{n,\ell}(G;uT_u(S),u')$ and $\MC_{n,\ell}(G;uT_u(S\cup w),u')$ have bases $P_{n,\ell}(G;uT_u(S),u')$ and $P_{n,\ell}(G;uT_u(S\cup w),u')$, respectively. The images of the elements of $P_{n,\ell}(G;uT^S_u(w),u')$ form a basis of $\displaystyle \frac{\MC_{n,\ell}(G;uT_u(S\cup w),u')}{\MC_{n,\ell}(G;uT_u(S),u')}$.
    Hence $\displaystyle \frac{\MC_{n,\ell}(G;uT_u(S\cup w),u')}{\MC_{n,\ell}(G;uT_u(S),u')}$ and $\mathcal{Q}_{n,\ell}(G;uT_u^S(w),u')$ are isomorphic as modules.
    Let $\pi$ denote the quotient map from $\MC_{*,\ell}(G;uT_u(S\cup w),u')$ to $\displaystyle \frac{\MC_{*,\ell}(G;uT_u(S\cup w),u')}{\MC_{*,\ell}(G;uT_u(S),u')}$.
    The boundary map $\overline{\partial}$ of the quotient complex is defined by $\overline{\partial}(\pi(z))=\pi(\partial(z))$, and the map in \eqref{eq:diff-of-Q} satisfies the defining equation for $\overline{\partial}$.
\end{proof}

The next proposition follows from \zcref{prop:chain-isom} and the definition of the quotient complex.

\begin{prop}\label{lem:chain-short-exact}
For every $\ell$, the following is a short exact sequence of chain complexes:
\begin{equation}\label{eq:chain-short-exact}
0\longrightarrow
 \MC_{*,\ell}(G;uT_u(S),u')
\xrightarrow{\ \iota\ }
 \MC_{*,\ell}(G;uT_u(S\cup w),u')
\xrightarrow{\ \pi\ }
 \mathcal{Q}_{*,\ell}(G;uT_u^S(w),u')
 \longrightarrow 0.
\end{equation}
For $n\ge1$, the map $\pi$ is given by
\begin{equation}
 \pi(u,x_1,\ldots,x_n)=
 \begin{cases}
 (u,x_1,\ldots,x_n) & \text{if } x_1\notin T_u(S),\\
 0 & \text{if } x_1\in T_u(S).
 \end{cases}
\end{equation}
\end{prop}

The zigzag lemma gives the following corollary.

\begin{coro}\label{coro:chain-long-exact}
There is a long exact sequence
\begin{align}\label{eq:chain-long-exact}
\cdots\longrightarrow\MH_{n,\ell}(G;uT_u(S),u')&\longrightarrow\MH_{n,\ell}(G;uT_u(S\cup w),u')\longrightarrow H_n(\mathcal{Q}_{*,\ell}(G;uT_u^S(w),u'))\\
&\longrightarrow\MH_{n-1,\ell}(G;uT_u(S),u')\longrightarrow\cdots .
\end{align}
\end{coro}

The rest of this subsection proves the following lemma.
\begin{lemma}\label{lem:chain-isom_2}
Fix $u,u'\in V$, $w\in\Gamma(u)$, and $\ell\ge 2$. Suppose $u$ is reducible to $w$ with respect to $S\subseteq\Gamma(u)$ and $R\subseteq\Gamma(w)$. Regard $\MC_{*-1,\ell-1}(G;wT_w(R),u')$ as a chain complex whose group in degree $n$ is $\MC_{n-1,\ell-1}(G;wT_w(R),u')$, with the original differential $\partial$. Then, for every $n\ge1$, the map
\begin{align}
 \phi:\MC_{n-1,\ell-1}(G;wT_w(R),u')&\longrightarrow\mathcal{Q}_{n,\ell}(G;uT^S_u(w),u'),\\
 (w,x_2,\ldots,x_n)&\longmapsto (-1)^n (u,w,x_2,\ldots,x_n) \label{eq:phi}
\end{align}
is a chain map and induces an isomorphism in homology
\begin{equation}\label{eq:wmd-replacement-iso}
 \phi_*:\MH_{n-1,\ell-1}(G;wT_w(R),u') \xrightarrow{\ \cong\ }H_n(\mathcal{Q}_{*,\ell}(G;uT^S_u(w),u')).
\end{equation}
\end{lemma}

\begin{proof}
Write $\phi_n$ for the map in \eqref{eq:phi} with target in homological degree $n$. We prove the lemma in three steps.

\noindent
\textbf{Step 1: $\phi$ is a chain map}
\begin{lemma}\label{lem:chain-isom_2_1}
Suppose that $u\in V$ and $w\in\Gamma(u)$ satisfy the reducible condition with respect to $S\subseteq\Gamma(u)$ and $R\subseteq\Gamma(w)$.
For every $\xi\in\MC_{n-1,\ell-1}(G;wT_w(R),u')$, we have
\begin{equation}
    \phi_{n-1}\partial(\xi)=\overline{\partial}\phi_n(\xi).
    \label{eq:signed-prefix-chain-map}
\end{equation}
In particular, $\phi$ is an injective chain map.
\end{lemma}

\begin{proof}
It suffices to check \eqref{eq:signed-prefix-chain-map} on a basis element $p=(w,x_2,\ldots,x_n)\in P_{n-1,\ell-1}(G;wT_w(R),u')$.
By \eqref{eq:phi},
\begin{align}
    \phi_{n-1}\partial(p)
    &=\phi_{n-1}\left(\sum_{i=1}^{n-2}(-1)^i\partial_i(p)\right)
    =\sum_{i=2}^{n-1}(-1)^{n+i}\partial_i(u,w,x_2,\ldots,x_n).
\end{align}
On the other hand,
\begin{equation}
    \overline{\partial}\phi_n(p)
    =\sum_{i=1}^{n-1}(-1)^{n+i}\overline{\partial}_i(u,w,x_2,\ldots,x_n).
\end{equation}
Since $w\in T_u^S(w)$, for $i\geq2$ we have $\overline{\partial}_i(u,w,x_2,\ldots,x_n)=\partial_i(u,w,x_2,\ldots,x_n)$.
Also, $x_2\in T_w(R)$, and the reducible condition \eqref{eq:reducible-cond} gives $x_2\notin T_u^S(w)$.
Thus either $\partial_1(u,w,x_2,\ldots,x_n)=0$ or $\partial_1(u,w,x_2,\ldots,x_n)\in\MC_{n-1,\ell}(G;uT_u(S),u')$. In either case, $\overline{\partial}_1(u,w,x_2,\ldots,x_n)=0$.
This proves \eqref{eq:signed-prefix-chain-map}.
The map $\phi_n$ is injective because it sends distinct basis elements to distinct signed basis elements.
\end{proof}

\noindent\textbf{Step 2: Construction of $\eta,\psi$; proof of $\psi\phi=\phi$}

For $n\ge1$, define $\eta_n:\mathcal{Q}_{n,\ell}(G;uT_u^S(w),u')\to\mathcal{Q}_{n+1,\ell}(G;uT_u^S(w),u')$ on basis elements by
\begin{equation}
    \eta_n(u,x_1,\ldots,x_n) =
    \begin{cases}
        (u,w,x_1,\ldots,x_n) &\text{if }x_1\neq w,\\
        0 &\text{if }x_1=w,
    \end{cases}\label{eq:eta}
\end{equation}
and extend it linearly. For $n\le0$, we set $\eta_n=0$.

\begin{lemma}\label{lem:chain-isom_2_eta}
The map $\eta_n$ is well-defined, and $\eta_n(\mathcal{Q}_{n,\ell}(G;uT_u^S(w),u'))\subseteq\MC_{n+1,\ell}(G;u\{w\},u')$.
Moreover, for every $\xi\in\mathcal{Q}_{n,\ell}(G;uT_u^S(w),u')$, we have $\xi+\overline{\partial}\eta_n(\xi)\in\MC_{n,\ell}(G;u\{w\},u')$.
\end{lemma}

\begin{proof}
Let $p=(u,x_1,\ldots,x_n)$ be a basis element.
If $x_1\neq w$, then $x_1\in T_u^S(w)$ gives $d(u,x_1)=d(u,w)+d(w,x_1)$, and hence $(u,w,x_1,\ldots,x_n)$ has length $\ell$.
Since $w\in T_u^S(w)$, the element $\eta_n(p)$ is a basis element of $\mathcal{Q}_{n+1,\ell}(G;uT_u^S(w),u')$ whose second entry is $w$.
If $x_1=w$, then $\eta_n(p)=0$.
In both cases, $\eta_n$ is well-defined and its image is contained in $\MC_{n+1,\ell}(G;u\{w\},u')$.

Next, we show that $p+\overline{\partial}\eta_n(p)\in\MC_{n,\ell}(G;u\{w\},u')$. First, suppose $x_1\neq w$.
Since $x_1\in T_u^S(w)$, we have $\overline{\partial}_1(u,w,x_1,\ldots,x_n)=(u,x_1,\ldots,x_n)=p$.
For $i\geq2$, a nonzero face $\overline{\partial}_i(u,w,x_1,\ldots,x_n)$ has second entry $w$.
Hence,
\begin{equation}
    \overline{\partial}\eta_n(p)=-p+\zeta_{\eta(p)}
    \label{eq:insertion-boundary}
\end{equation}
for some $\zeta_{\eta(p)}\in\MC_{n,\ell}(G;u\{w\},u')$.
Thus, $p+\overline{\partial}\eta_n(p)\in\MC_{n,\ell}(G;u\{w\},u')$.
If $x_1=w$, then $\eta_n(p)=0$ and $p$ itself belongs to $\MC_{n,\ell}(G;u\{w\},u')$.
The claim follows by linearity.
\end{proof}

Define $\psi:\mathcal{Q}_{*,\ell}(G;uT_u^S(w),u')\to\mathcal{Q}_{*,\ell}(G;uT_u^S(w),u')$ by
\begin{equation}\label{eq:prefix-retraction}
    \psi=\id+\overline{\partial}\eta+\eta\overline{\partial}.
\end{equation}

\begin{lemma}\label{lem:chain-isom_2_tau}
The map $\psi$ is a chain map, $\psi\simeq\id$, and $\psi\phi=\phi$.
Here $\simeq$ denotes chain homotopy.
\end{lemma}

\begin{proof}
Since $\eta$ is a module homomorphism in each degree, so is $\psi$.
Furthermore, because $\overline{\partial}^{\,2}=0$, we have
\begin{equation}
    \overline{\partial}\psi_n(\xi)
    =\overline{\partial}\xi+\overline{\partial}\eta_{n-1}(\overline{\partial}\xi)
    =\psi_{n-1}(\overline{\partial}\xi).
    \label{eq:normalization-commutes}
\end{equation}
Thus, $\psi$ is a chain map.
By \eqref{eq:prefix-retraction}, we have $\psi-\id=\overline{\partial}\eta+\eta\overline{\partial}$, so $\eta$ is a chain homotopy between $\psi$ and $\id$.
Hence, $\psi\simeq\id$.

For any basis element $p=(w,x_2,\ldots,x_n)$, \eqref{eq:phi} gives $\phi_n(p)=(-1)^n(u,w,x_2,\ldots,x_n)$.
Consequently, \eqref{eq:eta} gives $\eta_n\phi_n(p)=0$, and extending this by linearity results in $\eta\phi=0$.
By \zcref{lem:chain-isom_2_1}, we have $\overline{\partial}\phi=\phi\partial$, so 
\begin{equation}
\psi\phi=\phi+\overline{\partial}\eta\phi+\eta\overline{\partial}\phi=\phi+\overline{\partial}\eta\phi+\eta\phi\partial=\phi.
\end{equation}
\end{proof}

\noindent\textbf{Step 3: Construction of $\rho$ and showing $\phi\rho=\psi\simeq\id$}

\begin{lemma}\label{lem:chain-isom_2_2}
For every $\zeta\in\MC_{n,\ell}(G;u\{w\},u')$, if $\overline{\partial}\zeta\in\MC_{n-1,\ell}(G;u\{w\},u')$, then $\zeta\in\Im(\phi_n)$.
\end{lemma}

\begin{proof}
If $n=1$, then $\MC_{1,\ell}(G;u\{w\},u')=0$ since $\ell\ge2$.

Suppose $n\geq2$. By the reducible condition \eqref{eq:reducible-cond}, the basis expansion of $\zeta$ splits uniquely into two sums. The third entries lie in $T_w(R)$ in the first sum and in $T_u^S(w)$ in the second.
The basis elements whose third entry lies in $T_w(R)$ belong to $\Im(\phi_n)$ by \eqref{eq:phi}.
It therefore suffices to show that all coefficients of basis elements whose third entries lie in $T_u^S(w)$ are zero.

If $p=(u,w,x_2,\ldots,x_n)$ and $x_2\in T_u^S(w)$, then $d(u,x_2)=1+d(w,x_2)$ gives $\overline{\partial}_1(p)=(u,x_2,\ldots,x_n)$.
For $i\geq2$, a nonzero face $\overline{\partial}_i(p)$ has second entry $w$. If $x_2\in T_w(R)$, then the reducible condition gives $\overline{\partial}_1(p)=0$.
Since $\overline{\partial}\zeta\in\MC_{n-1,\ell}(G;u\{w\},u')$, each basis element with second entry in $T_u^S(w)\setminus\{w\}$ has coefficient zero in $\overline{\partial}\zeta$.
Distinct basis elements $(u,w,x_2,\ldots,x_n)$ give distinct faces $(u,x_2,\ldots,x_n)$, so these terms cannot cancel.
Thus, the coefficients of basis elements whose third entry lies in $T_u^S(w)$ are all zero, and hence $\zeta\in\Im(\phi_n)$.
\end{proof}

\begin{lemma}\label{lem:chain-isom_2_tau-image}
We have $\Im(\psi)\subseteq\Im(\phi)$.
\end{lemma}

\begin{proof}
Let $\xi\in\mathcal{Q}_{n,\ell}(G;uT_u^S(w),u')$. By \zcref{lem:chain-isom_2_eta}, we have $\xi+\overline{\partial}\eta_n(\xi)\in\MC_{n,\ell}(G;u\{w\},u')$ and $\eta_{n-1}(\overline{\partial}\xi)\in\MC_{n,\ell}(G;u\{w\},u')$.
Thus, $\psi_n(\xi)\in\MC_{n,\ell}(G;u\{w\},u')$.
Also, \eqref{eq:normalization-commutes} gives $\overline{\partial}\psi_n(\xi)=\overline{\partial}\xi+\overline{\partial}\eta_{n-1}(\overline{\partial}\xi)$.
Applying \zcref{lem:chain-isom_2_eta} to $\overline{\partial}\xi$ gives $\overline{\partial}\psi_n(\xi)\in\MC_{n-1,\ell}(G;u\{w\},u')$.
Thus, \zcref{lem:chain-isom_2_2} gives $\psi_n(\xi)\in\Im(\phi_n)$.
\end{proof}

By \zcref{lem:chain-isom_2_1, lem:chain-isom_2_tau-image}, for every $\xi\in\mathcal{Q}_{n,\ell}(G;uT_u^S(w),u')$, there exists a unique $\rho_n(\xi)\in\MC_{n-1,\ell-1}(G;wT_w(R),u')$ such that
\begin{equation}
    \phi_n(\rho_n(\xi))=\psi_n(\xi).
    \label{eq:retraction-definition}
\end{equation}
The maps $\phi_n,\psi_n$ are module homomorphisms, and the injectivity of $\phi_n$ ensures that $\rho_n$ is a module homomorphism as well.
The maps $\rho_n$ together define
\begin{equation}
    \rho:\mathcal{Q}_{*,\ell}(G;uT_u^S(w),u')\to\MC_{*-1,\ell-1}(G;wT_w(R),u').
\end{equation}

\begin{lemma}\label{lem:chain-isom_2_rho}
The map $\rho$ is a chain map.
Moreover, $\rho\phi=\id$ and $\phi\rho\simeq\id$.
\end{lemma}

\begin{proof}
Since $\phi,\psi$ are chain maps, we have
\begin{align}
    \phi_{n-1}(\partial\rho_n(\xi))
    =\overline{\partial}\phi_n(\rho_n(\xi))
    =\overline{\partial}\psi_n(\xi)
    =\psi_{n-1}(\overline{\partial}\xi)
    =\phi_{n-1}(\rho_{n-1}(\overline{\partial}\xi)).
\end{align}
The injectivity of $\phi_{n-1}$ gives $\partial\rho_n(\xi)=\rho_{n-1}(\overline{\partial}\xi)$.
Thus, $\rho$ is a chain map.

By \zcref{lem:chain-isom_2_tau}, we have $\psi\phi=\phi$.
Also, \eqref{eq:retraction-definition} gives $\phi\rho=\psi$, so $\phi\rho\phi=\psi\phi=\phi$.
The injectivity of $\phi$ then implies $\rho\phi=\id$.
Furthermore, since $\psi\simeq\id$ by \zcref{lem:chain-isom_2_tau}, it follows that $\phi\rho=\psi\simeq\id$.
\end{proof}

We now return to the proof of \zcref{lem:chain-isom_2}.
By \zcref{lem:chain-isom_2_rho}, the maps $\phi$ and $\rho$ are inverse up to chain homotopy. Hence $\MC_{*-1,\ell-1}(G;wT_w(R),u')$ and $\mathcal{Q}_{*,\ell}(G;uT_u^S(w),u')$ are chain homotopy equivalent.
In particular, $\phi_*$ is the isomorphism in \eqref{eq:wmd-replacement-iso}.
\end{proof}

\subsection{Proof of \zcref{thm:reducible}}
Let $G=(V,E)$ be a reducive graph, and fix a family $(\mathcal{R}_u)_{u\in V}$ of subsets of vertex neighborhoods satisfying conditions \ref{item:canonical-condition} and \ref{item:reducible-condition}.
Since $\MC_{n,\ell}(G)=0$ for $n>\ell$, it suffices to consider $n<\ell$.
We prove $\MH_{n,\ell}(G;uT_u(S),u')=0$ for all $u,u'\in V$, $S\in\mathcal{R}_u$, and $0\le n<\ell$, by induction on $(\ell,\abs{S})$ in lexicographic order, with the length degree $\ell$ first.

If $n=0$ or $S=\emptyset$, then $\MC_{n,\ell}(G;uT_u(S),u')=0$.
We may therefore assume $\ell\ge2$, $1\le n<\ell$, and $S\ne\emptyset$.
By \ref{item:reducible-condition}, we can choose $w\in S$ and $R\in\mathcal{R}_w$ such that $S\setminus w\in\mathcal{R}_u$ and $u$ is $(S\setminus w,R)$-reducible to $w$. Under the isomorphism in \zcref{lem:chain-isom_2}, the relevant part of the long exact sequence in \zcref{coro:chain-long-exact} becomes
\begin{equation}
    \MH_{n,\ell}(G;uT_u(S\setminus w),u')\xrightarrow{\iota_*}\MH_{n,\ell}(G;uT_u(S),u')\xrightarrow{\phi_*^{-1}\pi_*} \MH_{n-1,\ell-1}(G;wT_w(R),u').
\end{equation}
The left term vanishes by induction because $\abs{S\setminus w}<\abs{S}$. The right term also vanishes by induction because the length degree decreases and $n-1<\ell-1$.
By exactness, we obtain $\MH_{n,\ell}(G;uT_u(S),u')=0$.

Finally, \ref{item:canonical-condition} gives $\Gamma(u)\in\mathcal{R}_u$. Since $T_u(\Gamma(u))=V\setminus u$, we have $\MH_{n,\ell}(G;u,u')=\MH_{n,\ell}(G;uT_u(\Gamma(u)),u')=0$ for $0\le n<\ell$.
By \eqref{eq:chain-hom-direct-sum}, we conclude that $\MH_{n,\ell}(G)=0$. Thus $G$ is diagonal.\qed

\section{Criteria for reducivity}\label{sec:reducive-criteria}

We give two criteria for reducivity.
The first characterizes reducivity in terms of distances and common neighbors.
The second gives a sufficient condition using families of cliques that cover vertex neighborhoods.
These criteria will be used in the next section to prove reducivity for weakly modular graphs and basis graphs of even $\Delta$-matroids and $q$-matroids.

\subsection{Metric criteria}\label{subsec:tryy}

\begin{thm}\label{thm:tryy-0}
Let $G=(V,E)$ be a graph, and let $(\mathcal{R}_u)_{u\in V}$ be a collection of families with $\mathcal{R}_u\subseteq 2^{\Gamma(u)}$ for every $u\in V$.
Then $(\mathcal{R}_u)_{u\in V}$ satisfies condition~\ref{item:reducible-condition} if and only if it satisfies the following condition at every vertex:
{\setlength{\leftmargini}{55pt}
\begin{itemize}
    \labelitem{R2'}{item:wmc} For every nonempty $S\in\mathcal{R}_u$, there exists $w\in S$ satisfying the following three conditions:
    {
    \begin{itemize}
        \labelitem{R2'-1}{item:fc} We have $S\setminus w\in\mathcal{R}_u$ and $\Gamma(w)\cap\bigl(\{u\}\cup\Gamma(u)\cup\Gamma(S\setminus w)\bigr)\in\mathcal{R}_w$.
        \labelitem{R2'-2}{item:tc-uws} For every $x\in V$ such that $d(u,x)=d(w,x)=r\ge1$, we have
        \begin{equation}
            \Gamma(w)\cap\bigl(\Gamma(u)\cup\Gamma(S\setminus w)\bigr)\cap\Gamma_{r-1}(x)\ne\emptyset.
        \end{equation}
        \labelitem{R2'-3}{item:qc-uws} For every $x\in V$ such that $d(u,x)=r+1$, $d(w,x)=r\ge1$, and $(S\setminus w)\cap\Gamma_r(x)\ne\emptyset$, there exists $z\in(S\setminus w)\cap\Gamma_r(x)$ such that
        \begin{equation}
            \Gamma(w)\cap\Gamma(z)\cap\Gamma_{r-1}(x)\ne\emptyset.
        \end{equation}
    \end{itemize}}
\end{itemize}}
\end{thm}
Conditions \ref{item:tc-uws} and \ref{item:qc-uws} are closely related to the triangle and quadrangle conditions for weakly modular graphs.

\begin{proof}
We use two facts to show that \ref{item:reducible-condition} in \zcref{def:family} is equivalent to \ref{item:wmc}.

\begin{itemize}
    \item If $u$ is $(S\setminus w,R)$-reducible to $w$, then $R$ is uniquely determined by $u,w,S$ (see \zcref{lem:ur}).
    \item For $R$ given by \eqref{eq:tryy-R}, the vertex $u$ is $(S\setminus w,R)$-reducible to $w$ if and only if $(u,w,S)$ satisfies both \ref{item:tc-uws} and \ref{item:qc-uws} (see \zcref{lem:reducive-tryy}).
\end{itemize}

We first prove these two facts.

\begin{lemma}\label{lem:ur}
Let $u\in V$, $S\subseteq\Gamma(u)$, $w\in\Gamma(u)\setminus S$, and $R\subseteq\Gamma(w)$. If $u$ is $(S,R)$-reducible to $w$, then
\begin{equation}
    R=\Gamma(w)\cap\bigl(\{u\}\cup\Gamma(u)\cup\Gamma(S)\bigr).
\end{equation}
In particular, such an $R$ is unique.
\end{lemma}

\begin{proof}
If $u$ is $(S,R)$-reducible to $w$, then \eqref{eq:reducible-cond} gives $T_w(R)=V\setminus T_u^S(w)$.
\zcref{prop:target-neighbor} shows $T_w(R)\cap \Gamma(w)=R$.
Intersecting both sides with $\Gamma(w)$ gives $R=\Gamma(w)\setminus T_u^S(w)$.

For $y\in\Gamma(w)$, we have $d(u,y)\le2$. By the definition of the target set, $y\in T_u^S(w)$ if and only if $d(u,y)=2$ and $y\notin\Gamma(S)$.
Thus $\Gamma(w)\cap T_u^S(w)=\bigl(\Gamma(w)\cap\Gamma_2(u)\bigr)\setminus\Gamma(S)$.
It follows that $R=\Gamma(w)\setminus T_u^S(w)=\Gamma(w)\setminus \left(\bigl(\Gamma(w)\cap\Gamma_2(u)\bigr)\setminus\Gamma(S)\right)=\Gamma(w)\cap\bigl(\{u\}\cup\Gamma(u)\cup\Gamma(S)\bigr)$.
\end{proof}

In what follows, take $u\in V$, $\emptyset\ne S\subseteq\Gamma(u)$, and $w\in S$, and put
\begin{equation}\label{eq:tryy-R}
    R\coloneq\Gamma(w)\cap\bigl(\{u\}\cup\Gamma(u)\cup\Gamma(S\setminus w)\bigr).
\end{equation}
To prove \zcref{lem:reducive-tryy}, we first express \ref{item:tc-uws} and \ref{item:qc-uws} in terms of target sets.

\begin{lemma}\label{lem:TC-equi}
For $u,w,S,R$ as above, the following are equivalent.
\begin{itemize}
    \item The triple $(u,w,S)$ satisfies \ref{item:tc-uws}.
    \item For every $x\in V$, if $d(u,x)=d(w,x)=r\ge1$, then $x\in T_w(R)$.
\end{itemize}
\end{lemma}

\begin{proof}
Suppose $d(u,x)=d(w,x)=r\ge1$.
By the definition of a target set, $x\in T_w(R)$ is equivalent to $R\cap\Gamma_{r-1}(x)\ne\emptyset$.
Since $u\notin\Gamma_{r-1}(x)$, the formula in \eqref{eq:tryy-R} shows that $R\cap\Gamma_{r-1}(x)\ne\emptyset$ is exactly the common-neighbor requirement in \ref{item:tc-uws}.
\end{proof}

\begin{lemma}\label{lem:QC-equi}
For $u,w,S,R$ as above, the following two conditions are equivalent:
\begin{itemize}
    \item The triple $(u,w,S)$ satisfies \ref{item:qc-uws}.
    \item We have $T_u(w)\cap T_u(S\setminus w)\subseteq T_w(R)$.
\end{itemize}
\end{lemma}

\begin{proof}
Assume \ref{item:qc-uws} and let $x\in T_u(w)\cap T_u(S\setminus w)$. Put $r=d(w,x)$. Then $d(u,x)=r+1$ and $(S\setminus w)\cap\Gamma_r(x)\ne\emptyset$.
If $r=0$, then $x=w\notin T_u(S\setminus w)$, a contradiction. Thus $r\ge1$.
By \ref{item:qc-uws}, there exist $z\in(S\setminus w)\cap\Gamma_r(x)$ and $y\in\Gamma(w)\cap\Gamma(z)$ such that $d(y,x)=r-1$.
Since $\Gamma(z)\subseteq \Gamma(S\setminus w)$, we have $\Gamma(w)\cap\Gamma(z)\subseteq R$, so $x\in T_w(R)$.

Conversely, assume $T_u(w)\cap T_u(S\setminus w)\subseteq T_w(R)$ and take $x$ as in \ref{item:qc-uws}. Since $x\in T_u(w)\cap T_u(S\setminus w)$, there exists $y\in R$ with $d(y,x)=r-1$. The triangle inequality gives $d(u,y)\ge d(u,x)-d(y,x)=2$, and $y\in\Gamma(w)$ gives $d(u,y)\le2$. Thus $d(u,y)=2$. By \eqref{eq:tryy-R}, there is some $z\in S\setminus w$ such that $y\in\Gamma(z)$. The bounds
\begin{equation}
    r=d(u,x)-1\le d(z,x)\le1+d(y,x)=r
\end{equation}
give $z\in\Gamma_r(x)$. Hence $y\in\Gamma(w)\cap\Gamma(z)\cap\Gamma_{r-1}(x)$, as required in \ref{item:qc-uws}.
\end{proof}

We can now prove the following lemma.

\begin{lemma}\label{lem:reducive-tryy}
For $u,w,S,R$ as above, the following two conditions are equivalent:
\begin{itemize}
    \item The vertex $u$ is $(S\setminus w,R)$-reducible to $w$.
    \item The triple $(u,w,S)$ satisfies \ref{item:tc-uws} and \ref{item:qc-uws}.
\end{itemize}
\end{lemma}

\begin{proof}
First suppose that $u$ is $(S\setminus w,R)$-reducible to $w$. Then
\begin{equation}\label{eq:tryy-partition}
    V=T_u^{S\setminus w}(w)\sqcup T_w(R).
\end{equation}
Let $x\in V$. If $d(u,x)=d(w,x)$, then $x\notin T_u(w)$, so \eqref{eq:tryy-partition} gives $x\in T_w(R)$.
If $x\in T_u(w)\cap T_u(S\setminus w)$, then $x\notin T_u^{S\setminus w}(w)$ and again $x\in T_w(R)$. By \zcref{lem:TC-equi,lem:QC-equi}, this proves \ref{item:tc-uws} and \ref{item:qc-uws}.

Conversely, suppose that $(u,w,S)$ satisfies \ref{item:tc-uws} and \ref{item:qc-uws}. For each $x\in V$, consider three cases.
\begin{itemize}
    \item Suppose $x\in T_u(w)$ and $x\notin T_u(S\setminus w)$, that is, $x\in T_u^{S\setminus w}(w)$.
    If $x\in T_w(R)$, there exists $y\in R$ such that $d(y,x)=d(w,x)-1=d(u,x)-2$. The triangle inequality gives $y\notin\{u\}\cup\Gamma(u)$, so by \eqref{eq:tryy-R} there exists $z\in S\setminus w$ with $y\in\Gamma(z)$. Then
    \begin{equation}
        d(u,x)-1\le d(z,x)\le1+d(y,x)=d(u,x)-1
    \end{equation}
    gives $x\in T_u(z)\subseteq T_u(S\setminus w)$, a contradiction. Thus $x\notin T_w(R)$.
    \item Suppose $x\in T_u(w)$ and $x\in T_u(S\setminus w)$.
    By \zcref{lem:QC-equi}, we have $x\in T_w(R)$.
    \item Suppose $x\notin T_u(w)$.
    Since $u$ and $w$ are adjacent, either $d(u,x)=d(w,x)-1$ or $d(u,x)=d(w,x)$. In the first case, $u\in R$ gives $x\in T_w(u)\subseteq T_w(R)$. In the second case, $d(u,x)=d(w,x)\ge1$, so \zcref{lem:TC-equi} gives $x\in T_w(R)$.
\end{itemize}
This proves \eqref{eq:tryy-partition}, so $u$ is $(S\setminus w,R)$-reducible to $w$.
\end{proof}
We now complete the proof of \zcref{thm:tryy-0}.
Fix the families $(\mathcal{R}_u)_{u\in V}$.

For each nonempty $S\in\mathcal{R}_u$, condition \ref{item:reducible-condition} gives $w\in S$ and $R\in\mathcal{R}_w$ with $S\setminus w\in\mathcal{R}_u$ such that $u$ is $(S\setminus w,R)$-reducible to $w$.
By \zcref{lem:ur}, the set $R$ is given by \eqref{eq:tryy-R}, so \ref{item:fc} holds.
By \zcref{lem:reducive-tryy}, both \ref{item:tc-uws} and \ref{item:qc-uws} hold. This proves \ref{item:wmc}.

Conversely, suppose families $\mathcal{R}_u$ satisfy \ref{item:wmc}.
For a nonempty $S\in\mathcal{R}_u$, choose $w\in S$ as in \ref{item:wmc} and define $R$ by \eqref{eq:tryy-R}. By \ref{item:fc}, we have $S\setminus w\in\mathcal{R}_u$ and $R\in\mathcal{R}_w$. By \zcref{lem:reducive-tryy}, the vertex $u$ is $(S\setminus w,R)$-reducible to $w$. Thus \ref{item:reducible-condition} holds.
\end{proof}

\begin{coro}\label{thm:tryy}
A graph $G=(V,E)$ is reducive if and only if there exists a family $(\mathcal{R}_u)_{u\in V}$ with $\mathcal{R}_u\subseteq2^{\Gamma(u)}$ for each $u\in V$, satisfying \ref{item:canonical-condition} and \ref{item:wmc}. 
\end{coro}

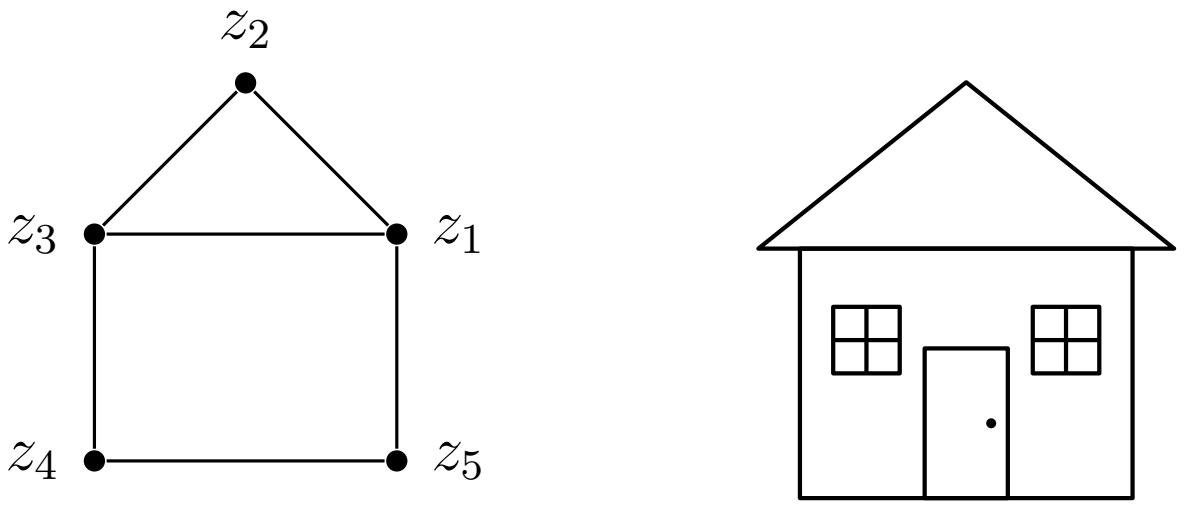
\begin{figure}[htbp]
\centering
\begin{tikzpicture}[vertex/.style={circle,fill=black,minimum size=4pt,inner sep=0pt}, line width=0.6pt]
    \node[vertex,label=left:$z_3$] (v1) at (-1,1.5) {};
    \node[vertex,label=right:$z_1$] (v2) at (1,1.5) {};
    \node[vertex,label=right:$z_5$] (v3) at (1,0) {};
    \node[vertex,label=left:$z_4$] (v4) at (-1,0) {};
    \node[vertex,label=above:$z_2$] (v5) at (0,2.5) {};
    \draw (v1)--(v2)--(v3)--(v4)--(v1);
    \draw (v1)--(v5)--(v2);
\end{tikzpicture}\qquad\qquad
\begin{tikzpicture}[thick, line join=round, scale=0.55]
  \draw (0,0) rectangle (4,3);
  \draw (-0.5,3) -- (2,5) -- (4.5,3) -- cycle;
  \draw (1.5,0) rectangle (2.5,1.8);
  \fill (2.3,0.9) circle (0.06);
  \draw (0.4,1.5) rectangle (1.2,2.3);
  \draw (0.8,1.5) -- (0.8,2.3);
  \draw (0.4,1.9) -- (1.2,1.9);
  \draw (2.8,1.5) rectangle (3.6,2.3);
  \draw (3.2,1.5) -- (3.2,2.3);
  \draw (2.8,1.9) -- (3.6,1.9);
\end{tikzpicture}

\caption{The House graph (left) and a house (right).}
\label{fig:house}
\end{figure}

\begin{exam}[House graph]\label{def:house}
The \emph{House graph} is obtained by adding a chord $z_1z_3$ to a $5$-cycle $(\{z_i\}_{i\in [5]}, \{z_iz_{i+1}\mid i\in [5]\})$ where $z_6=z_1$ (see \zcref{fig:house}).
Tajima and Yoshinaga showed that the House graph is diagonal \cite[Example~5.5]{TY23}.
\end{exam}

\begin{prop}\label{exam:house-tryy}
The House graph $G$ is reducive.
\end{prop}

\begin{proof}
Using the vertex labels in \zcref{def:house}, define
\begin{align}
    \mathcal{R}_{z_j}&\coloneq2^{\Gamma(z_j)}\qquad(j\in\{1,2,3\}),\\
    \mathcal{R}_{z_5}&\coloneq\{\emptyset,\{z_1\},\Gamma(z_5)\},\qquad
    \mathcal{R}_{z_4}\coloneq\{\emptyset,\{z_3\},\Gamma(z_4)\}.
\end{align}
For each $u\in V$, we have $\Gamma(u)\in\mathcal{R}_u$. For a nonempty $S\in\mathcal{R}_u$, we can choose $w\in S$ such that $S\setminus w\in\mathcal{R}_u$. Fix such a $w$.

Define $R$ by \eqref{eq:tryy-R}. The set $R$ always contains $u$, so if $uw\ne z_5z_4$, then $R\in\mathcal{R}_w$.
If $uw=z_5z_4$, then $S\setminus w\in\mathcal{R}_u$ forces $S=\Gamma(u)$, and $R=\Gamma(w)\in\mathcal{R}_w$.
Thus \ref{item:fc} holds.

For $r\ge2$ in \ref{item:tc-uws}, the common neighbor required by \ref{item:triangle-condition} exists except when $x=z_2$ and $uw=z_5z_4$. In this case $S=\Gamma(u)$. We take $y=z_3$ if $(u,w)=(z_5,z_4)$, and $y=z_1$ if $(u,w)=(z_4,z_5)$. In either case, $y\in\Gamma(w)\cap\Gamma(S\setminus w)\cap\Gamma(x)$. Thus \ref{item:tc-uws} holds.
Since the diameter of $G$ is $2$, we only need to check \ref{item:qc-uws} for $r=1$.
For $r=1$, take $x$ as the common neighbor.
Therefore \zcref{thm:tryy} shows that $G$ is reducive.
\end{proof}

\subsection{Clique families}\label{subsec:basis-reducive}
We give a criterion for reducivity using families of cliques that cover vertex neighborhoods.
Recall that the empty set is a clique.

\begin{thm}\label{lem:basis-reduction}
Let $G=(V,E)$ be a graph. Suppose there exists a family $(\mathcal{C}_u)_{u\in V}$ with $\mathcal{C}_u\subseteq2^{\Gamma(u)}$ for each $u\in V$.
Here $\mathcal{C}_u$ is a family of cliques in the neighborhood of $u$, and we assume for each vertex $u$ that $\emptyset\in\mathcal C_u$.
If the following two conditions hold, then $G$ is reducive:
\begin{itemize}
    \labelitem{C1}{item:basis-reduction-clique} For each $u\in V$, we have $\Gamma(u)=\bigcup_{C\in\mathcal{C}_u}C$.
    \labelitem{C2}{item:basis-reduction-partition} For each $u\in V$, $C\in\mathcal{C}_u$, and $w\in\Gamma(u)\setminus C$, there exists a union $R$ of members of $\mathcal{C}_w$ such that $V=T_u^C(w)\sqcup T_w(R)$.
\end{itemize}
\end{thm}

\begin{proof}
We show that the families $\mathcal{R}_u$ defined below satisfy \ref{item:canonical-condition} and \ref{item:reducible-condition}:
\begin{equation}
    \mathcal{R}_u\coloneq\left\{S_0\cup\bigcup_{C\in\mathcal{N}}C\midd\mathcal{N}\subseteq\mathcal{C}_u,\ S_0\subseteq\Gamma(u)\text{ is a clique}\right\}.
\end{equation}

By \ref{item:basis-reduction-clique}, we have $\Gamma(u)\in\mathcal{R}_u$, so \ref{item:canonical-condition} holds.

Write a nonempty $S\in\mathcal{R}_u$ as $S=C_1\cup\cdots\cup C_m$, where $C_1,\ldots,C_{m-1}\in\mathcal{C}_u$ and $C_m$ is a clique. Choose a representation with the smallest $m$.
Put $N=\bigcup_{j=1}^{m-1}C_j$ and $S_0=C_m\setminus N$. When $m=1$, set $N=\emptyset$.
If $S_0=\emptyset$, then $S=N$, which contradicts $S\ne\emptyset$ or the minimality of $m$.
Thus $S_0$ is a nonempty clique and $S=N\sqcup S_0$.
Choose $w\in S_0$. Since $S_0\setminus w$ is also a clique, we have $S\setminus w=N\sqcup(S_0\setminus w)\in\mathcal{R}_u$.

Put $C_0=\emptyset$. Since $w\notin C_j$, condition \ref{item:basis-reduction-partition} gives, for each $j=0,\ldots,m-1$, a union $R_j$ of members of $\mathcal{C}_w$ such that $V=T_u^{C_j}(w)\sqcup T_w(R_j)$.
Suppose first that $S_0\setminus w\neq\emptyset$.
For each $y\in S_0\setminus w$, we have $w\in\Gamma(y)$. Applying the same condition to $\emptyset\in\mathcal{C}_y$ gives a union $R_y'$ of members of $\mathcal{C}_w$ such that $V=T^\emptyset_y(w)\sqcup T_w(R_y')$.
Put $R=\bigcup_{j=0}^{m-1}R_j\cup\bigcup_{y\in S_0\setminus w}R_y'$. Since $R$ is a union of members of $\mathcal{C}_w$, we have $R\in\mathcal{R}_w$.
Target sets preserve unions by \zcref{prop:target-union}, so $T_w(R)=\bigcup_{j=0}^{m-1}T_w(R_j)\cup\bigcup_{y\in S_0\setminus w}T_w(R_y')$. Hence
\begin{align}
    V\setminus T_w(R)
    &=\bigcap_{j=0}^{m-1}\bigl(V\setminus T_w(R_j)\bigr)\cap\bigcap_{y\in S_0\setminus w}\bigl(V\setminus T_w(R_y')\bigr)\\
    &=\bigcap_{j=0}^{m-1}T_u^{C_j}(w)\cap\bigcap_{y\in S_0\setminus w}T_y(w)
    =T_u^N(w)\cap\bigcap_{y\in S_0\setminus w}T_y(w).\label{eq:basis-reduction1}
\end{align}
If $S_0\setminus w=\emptyset$, put $R=\bigcup_{j=0}^{m-1}R_j$. The same argument gives
\begin{align}
    V\setminus T_w(R)&=\bigcap_{j=0}^{m-1}T_u^{C_j}(w)=T_u^{S\setminus w}(w).\label{eq:basis-reductionES}
\end{align}
This proves \ref{item:reducible-condition} in this case. For $S_0\setminus w\neq\emptyset$, it remains to verify \ref{item:reducible-condition}.

Take $x\in T_u(w)$ and $y\in S_0\setminus w$.
Since $y\in\Gamma(u)\cap\Gamma(w)$, the identity $d(u,x)=d(w,x)+1$ and the triangle inequality give $d(w,x)\le d(y,x)\le d(w,x)+1$.
The equality $d(y,x)=d(w,x)$ holds exactly when $x\in T_u(y)$, and $d(y,x)=d(w,x)+1$ holds exactly when $x\in T_y(w)$.
Thus $T_u(w)\cap T_y(w)=T_u(w)\setminus T_u(y)$.
Since $T_u^N(w)\subseteq T_u(w)$, applying this identity to each $y\in S_0\setminus w$ gives
\begin{align}
    T_u^N(w)\cap\bigcap_{y\in S_0\setminus w}T_y(w)
    =T_u^N(w)\setminus\bigcup_{y\in S_0\setminus w}T_u(y)
    =T_u^{N\sqcup(S_0\setminus w)}(w)=T_u^{S\setminus w}(w).\label{eq:basis-reduction2}
\end{align}
Combining \eqref{eq:basis-reduction1} and \eqref{eq:basis-reduction2} gives $V=T_u^{S\setminus w}(w)\sqcup T_w(R)$, so \ref{item:reducible-condition} holds.
\end{proof}

\begin{rema}
    The converse of \zcref{lem:basis-reduction} does not hold in general.
    A reducive graph need not have clique families satisfying the hypotheses of \zcref{lem:basis-reduction}.
    The House graph in \zcref{exam:house-tryy} is an example.
    See \zcref{sec:appendix-house} for details.
\end{rema}

\section{Applications and limitations}\label{sec:class}

We consider three subclasses of reducive graphs: weakly modular graphs, basis graphs of even $\Delta$-matroids, and basis graphs of $q$-matroids.
We also give a meshed graph that is not diagonal. Meshed graphs include all three classes.

\subsection{Weakly modular graphs}\label{subsec:wm-mesh}
Weakly modular graphs include both pawful and median graphs.
\subsubsection{Weakly modular graphs are reducive}

\begin{defi}[{\cite[Definition~2.1]{CCHO20}}]
A connected graph $G=(V,E)$ is \emph{weakly modular} if, for every $x\in V$, the following conditions hold:
\begin{itemize}
\labelitem{TC}{item:triangle-condition}
If $uw\in E$ and $d(x,u)=d(x,w)=r$, then $\Gamma(u)\cap\Gamma(w)\cap\Gamma_{r-1}(x)\ne\emptyset$.
\labelitem{QC}{item:quadrangle-condition}
If $d(u,w)=2$ and there exists $z\in\Gamma(u)\cap\Gamma(w)$ with $d(x,u)=d(x,w)=d(x,z)-1=r$, then $\Gamma(u)\cap\Gamma(w)\cap\Gamma_{r-1}(x)\ne\emptyset$.
\end{itemize}
\end{defi}
These conditions are called the {triangle condition} and the {quadrangle condition}, respectively.

\begin{thm}\label{prop:wmg-reducible}
Weakly modular graphs are reducive and hence diagonal.
\end{thm}

\begin{proof}
    Let $G$ be a weakly modular graph with vertex set $V$.
    For each $u\in V$, define $\mathcal{R}_u\coloneq 2^{\Gamma(u)}$.
    By \zcref{thm:tryy}, it suffices to show that $(u,w,S)$ satisfies \ref{item:wmc}.
    By construction, the families $\mathcal{R}_u$ satisfy \ref{item:fc}.
    
    First suppose $d(u,x)=d(w,x)=r\ge1$.
    Since $uw\in E$, condition \ref{item:triangle-condition} gives
    \begin{equation}\label{eq:tcdef}
    \Gamma(w)\cap\bigl(\Gamma(u)\cup\Gamma(S\setminus w)\bigr)\cap\Gamma_{r-1}(x)
    \supseteq\Gamma(u)\cap\Gamma(w)\cap\Gamma_{r-1}(x)\ne\emptyset.
    \end{equation}
    Hence \ref{item:tc-uws} holds.
    
    Next suppose $d(u,x)=r+1$, $d(w,x)=r\ge1$, and $(S\setminus w)\cap\Gamma_r(x)\ne\emptyset$. Choose $z\in(S\setminus w)\cap\Gamma_r(x)$.
    Since $w,z\in\Gamma(u)$ and $w\ne z$, the distance $d(w,z)$ is either $1$ or $2$.
    If $d(w,z)=1$, then $d(w,x)=d(z,x)=r$ allows us to apply \ref{item:triangle-condition}.
    If $d(w,z)=2$, then $u\in\Gamma(w)\cap\Gamma(z)$ and $d(u,x)=r+1$ allow us to apply \ref{item:quadrangle-condition}.
    In either case,
    \begin{equation}\label{eq:qcdef}
    \Gamma(w)\cap\Gamma(z)\cap\Gamma_{r-1}(x)\ne\emptyset.
    \end{equation}
    Hence \ref{item:qc-uws} also holds.
    Thus, $G$ is reducive, and its diagonality follows from \zcref{thm:reducible}.
\end{proof}

\begin{rema}\label{rema:house-not-wmg}
The House graph is not weakly modular.
Indeed, with the labels in \zcref{def:house}, we have
\begin{equation}
d(z_2,z_5)=d(z_2,z_4)=2,\qquad z_5z_4\in E,\qquad \Gamma(z_5)\cap\Gamma(z_4)=\emptyset,
\end{equation}
so \ref{item:triangle-condition} fails for the vertex $z_2$ and the edge $z_5z_4$.
Hence the converse of \zcref{prop:wmg-reducible} does not hold.
\end{rema}

\subsubsection{Examples of weakly modular graphs}\label{subsubsec:emg-ex}

\begin{prop}\label{prop:wmg-paw-med}
A graph of diameter at most $2$ is weakly modular if and only if it is pawful.
\end{prop}

\begin{proof}
Suppose that $G$ is weakly modular and has diameter at most $2$. Given three vertices with $d(x,y)=d(y,z)=2$ and $d(x,z)=1$, apply \ref{item:triangle-condition} to the vertex $y$ and the edge $xz$ to obtain a common neighbor of the three vertices. Hence $G$ is pawful.

Conversely, suppose that $G$ is pawful. For \ref{item:triangle-condition}, take $x$ as the common neighbor when $r=1$, and use the pawful condition when $r=2$. Under the hypotheses of \ref{item:quadrangle-condition}, we have $1\le r<D\le2$, hence $r=1$, and we can again take $x$ as the common neighbor.
\end{proof}

Complete graphs and complete multipartite graphs are pawful, and hence they are weakly modular by \zcref{prop:wmg-paw-med}.
Weak modularity is also preserved by gluing along a single vertex \cite[Proposition~2.18]{CCHO20}.
For Cartesian products, we have the following equivalence.

\begin{prop}\label{prop:cartesian-wmg}
For nonempty connected graphs $G$ and $G'$, the product $G\square G'$ is weakly modular if and only if both $G$ and $G'$ are weakly modular.
\end{prop}

\begin{proof}
If $G$ and $G'$ are both weakly modular, then $G\square G'$ is weakly modular \cite[Proposition~2.18]{CCHO20}.
Conversely, suppose that $G\square G'$ is weakly modular.
Fix $z\in V(G')$. The induced subgraph on $V(G)\times\{z\}$ is isomorphic to $G$, and every shortest path between two of its vertices stays in this subgraph.
Now consider three vertices $(u,z),(w,z),(x,z)$ of $G\square G'$.
A common neighbor given by \ref{item:triangle-condition} or \ref{item:quadrangle-condition} lies on a shortest path from $(x,z)$ to $(u,z)$ or $(w,z)$. The common neighbor therefore belongs to $V(G)\times\{z\}$.
Hence $G$ is weakly modular. Interchanging the factors gives the weak modularity of $G'$.
\end{proof}

By \zcref{prop:cartesian-wmg} and the preservation of weak modularity under gluing along single vertices, we obtain the following.
\begin{prop}\label{prop:ham-wmg}
Hamming graphs and block graphs are weakly modular.
\end{prop}

Median graphs are also weakly modular.
\begin{prop}[{\cite[Section~6.1]{BC08}}]
    Every median graph is weakly modular.
\end{prop}

\begin{defi}[{\cite[Section~12.4.2]{BH2012}}]
The \emph{Johnson graph} $J(m,\kappa)$ has vertex set $\binom{[m]}\kappa$, with two vertices $u,w$ adjacent if $\abs{u\cap w}=\kappa-1$.
\end{defi}

\begin{prop}[{\cite[Proposition~4.4(1)]{CCHO20}}]\label{prop:johnson-wmg}
The Johnson graph $J(m,\kappa)$ is weakly modular.
\end{prop}

\begin{defi}[{\cite[Section~12.4.3]{BH2012}}]
Let $q$ be a prime power, and let $0\le \kappa\le m$ be integers.
The \emph{Grassmann graph} $J_q(m,\kappa)$ has the $\kappa$-dimensional subspaces of $\mathbb{F}_q^m$ as vertices, with two vertices $u,w$ adjacent if $\dim(u\cap w)=\kappa-1$.
\end{defi}

\begin{prop}\label{prop:Grassmann-wmg}
The Grassmann graph $J_q(m,\kappa)$ is weakly modular.
\end{prop}
\begin{proof}
First, we have $d(u,w)=\kappa-\dim(u\cap w)$.
Along one edge, the dimension of the intersection with a fixed subspace changes by at most $1$, so the right-hand side is a lower bound on the distance.
Conversely, extend a basis of $u\cap w$ to bases of $u,w$, respectively. Replacing the added basis vectors one at a time gives a path of the required length.

Suppose $d(x,u)=d(x,w)=r$ and put $s=\kappa-r$. For $r=1$, the vertex $x$ itself is a common neighbor as required by \ref{item:triangle-condition} and \ref{item:quadrangle-condition}. Thus \ref{item:triangle-condition} and \ref{item:quadrangle-condition} hold in this case.
We assume $r\ge2$ below and verify \ref{item:triangle-condition} and \ref{item:quadrangle-condition}.

To prove \ref{item:triangle-condition}, suppose $uw\in E$ and put $U=u\cap w$. Since $\dim U=\kappa-1$, the dimension $\dim(x\cap U)$ is either $s$ or $s-1$.

If $\dim(x\cap U)=s$, then $x\cap u=x\cap w=x\cap U$.
Since $\dim x>s$, choose a one-dimensional subspace $W\le x$ with $W\cap U=\{0\}$.
Put $y=U+W$. The subspace $W$ is contained in neither of the subspaces $u,w$, and we obtain
\begin{equation}\label{eq:Grassmann-common-neighbor}
\dim y=\kappa,\qquad
\dim(y\cap u)=\dim(y\cap w)=\kappa-1,\qquad
\dim(y\cap x)=s+1.
\end{equation}

If $\dim(x\cap U)=s-1$, take an $(\kappa-2)$-dimensional subspace $W\le U$ containing $x\cap U$, and set $y=W+(x\cap u)+(x\cap w)$. The subspace $(x\cap u)+(x\cap w)$ has dimension $s+1$, and its intersection with $U$ is $x\cap U$. Thus
\begin{align}
\dim y&=(\kappa-2)+(s+1)-(s-1)=\kappa,\\
y\cap u&=W+(x\cap u),\qquad y\cap w=W+(x\cap w),\\
y\cap x&=(x\cap u)+(x\cap w).
\end{align}
These dimension formulas give \eqref{eq:Grassmann-common-neighbor}. In either case, $y\in\Gamma(u)\cap\Gamma(w)\cap\Gamma_{r-1}(x)$.

To prove \ref{item:quadrangle-condition}, suppose $d(u,w)=2$, take $z\in\Gamma(u)\cap\Gamma(w)\cap\Gamma_{r+1}(x)$, and let $U=u\cap w$.
Applying the dimension formula within $z$ gives
\begin{equation}
\dim(U\cap z)\ge\dim(u\cap z)+\dim(w\cap z)-\dim z=\kappa-2=\dim U.
\end{equation}
Hence we have $U\le z$.
Applying the dimension formula within $u$ also gives
\begin{equation}
\dim(x\cap u\cap z)\ge\dim(x\cap u)+\dim(u\cap z)-\dim u=s-1=\dim(x\cap z),
\end{equation}
so $x\cap z\le u$. Applying the same argument with $w$ gives $x\cap z\le w$. Thus $x\cap U=x\cap z$, and $\dim(x\cap U)=s-1$. Put $y=U+(x\cap u)+(x\cap w)$. Then
\begin{align}
\dim y&=(\kappa-2)+(s+1)-(s-1)=\kappa,\\
y\cap u&=U+(x\cap u),\qquad y\cap w=U+(x\cap w),\\
y\cap x&=(x\cap u)+(x\cap w).
\end{align}
The dimension formulas give \eqref{eq:Grassmann-common-neighbor} and hence the required common neighbor.
\end{proof}

By \zcref{prop:wmg-reducible, prop:johnson-wmg, prop:Grassmann-wmg}, Johnson graphs and Grassmann graphs are diagonal.

In the next definition, we use generalized polar spaces and their singular subspaces in the sense of \cite[Section~2.2.3]{CCHO20}.

\begin{defi}[{\cite[Section~2.2.3]{CCHO20}}]
The \emph{dual polar graph} of a nondegenerate generalized polar space has its maximal singular subspaces as vertices.
Two distinct vertices are adjacent if their intersection is a maximal proper singular subspace of both.
\end{defi}

Dual polar graphs are weakly modular \cite[Theorem~5.2]{CCHO20}, so they are diagonal by \zcref{prop:wmg-reducible}.

\begin{defi}[{\cite[Section~2.1.1]{CCHO20}}]\label{def:halved-cube}
For an integer $m\ge1$, the \emph{halved $m$-cube}, or \emph{half cube}, $Q_m^+$ has as vertices the subsets of $[m]$ of even cardinality. Two vertices $u,w$ are adjacent exactly when $\abs{u\triangle w}=2$.
\end{defi}

\begin{defi}[{\cite[Section~10.10]{BvM22}}]\label{def:gosset-schlafli}
The \emph{Gosset graph} has vertex set $\binom{[8]}2\times\{+1,-1\}$. Distinct vertices $(S,\varepsilon),(S',\varepsilon')$ are adjacent if
\begin{equation}
\bigl(\varepsilon=\varepsilon'\ \text{and}\ \abs{S\cap S'}=1\bigr)
\quad\text{or}\quad
\bigl(\varepsilon\ne\varepsilon'\ \text{and}\ S\cap S'=\emptyset\bigr).
\end{equation}
The induced subgraph on the neighborhood of $(\{1,2\},+1)$ in the Gosset graph is the \emph{Schl\"afli graph}.
\end{defi}

Halved cubes, the Schl\"afli graph, and the Gosset graph are weakly modular \cite[Proposition~4.4(1)]{CCHO20}. By \zcref{prop:wmg-reducible}, they are diagonal as well.

\subsubsection{Triangle-square complexes}\label{subsubsec:triangle-square}

The diagonality of weakly modular graphs, proved in \zcref{prop:wmg-reducible}, is also related to the simple connectivity of a complex obtained by filling in triangles and squares.

For a connected graph $G$, the \emph{triangle-square complex} $\mathsf{CW}_{\triangle\square}(G)$ is the CW complex obtained from $G$ by attaching a $2$-cell along the boundary of each triangle and each chordless square \cite[Section~2.1.2]{CCHO20}.
Note that its $1$-skeleton is $G$.

A connected graph $G=(V,E)$ is called \emph{locally weakly modular} if \ref{item:triangle-condition} and \ref{item:quadrangle-condition} hold for $r=2$ and every $x\in V$ \cite[Definition~2.6]{CCHO20}.

\begin{thm}\label{thm:triangle-square}
Let $G$ be a finite connected locally weakly modular graph.
The following conditions are equivalent:
\begin{enumerate}[(i)]
    \labelitem{i}{item:triangle-square-simply-connected} The complex $\mathsf{CW}_{\triangle\square}(G)$ is simply connected.
    \labelitem{ii}{item:triangle-square-wmg} The graph $G$ is weakly modular.
    \labelitem{iii}{item:triangle-square-diagonal} The graph $G$ is diagonal.
    \labelitem{iv}{item:triangle-square-mh2} For every integer $\ell>2$, we have $\MH_{2,\ell}(G)=0$.
    \labelitem{v}{item:triangle-square-koszul} The distance algebra $\mathsf{DA}_{\mathbb{F}}(G)$ is Koszul over every field $\mathbb{F}$.
    \labelitem{vi}{item:triangle-square-quadratic} The distance algebra $\mathsf{DA}_{\mathbb{F}}(G)$ is quadratic over every field $\mathbb{F}$.
\end{enumerate}
\end{thm}

\begin{proof}
The equivalence of \ref{item:triangle-square-simply-connected} and \ref{item:triangle-square-wmg} is due to Chalopin, Chepoi, Hirai, and Osajda \cite[Theorem~3.1]{CCHO20}.
If \ref{item:triangle-square-wmg} holds, then $G$ is reducive by \zcref{prop:wmg-reducible} and hence diagonal by \zcref{thm:reducible}. This gives \ref{item:triangle-square-diagonal}.
By the definition of diagonality, \ref{item:triangle-square-diagonal} implies \ref{item:triangle-square-mh2}.

Now assume \ref{item:triangle-square-mh2}.
By the result of Ivanov and Mukoseev \cite[Corollary~4.22]{IM26}, the complex obtained from $G$ by attaching $2$-cells along all $3$-cycles and $4$-cycles is simply connected.
In $\mathsf{CW}_{\triangle\square}(G)$, the boundary of a chorded square is null-homotopic because it is bounded by two triangular $2$-cells.
Thus attaching $2$-cells along squares with chords does not change the fundamental group. It follows that $\mathsf{CW}_{\triangle\square}(G)$ is simply connected, which proves \ref{item:triangle-square-simply-connected}.

The equivalences of \ref{item:triangle-square-diagonal} with \ref{item:triangle-square-koszul} and of \ref{item:triangle-square-mh2} with \ref{item:triangle-square-quadratic} follow from \cite[Theorem~6.5 and Proposition~6.4]{IM26}, respectively.
\end{proof}

Chalopin, Chepoi, Hirai, and Osajda describe the equivalence of \ref{item:triangle-square-simply-connected} and \ref{item:triangle-square-wmg} as a discrete version of the Cartan--Hadamard theorem \cite[Section~3.1]{CCHO20}.
The result in \zcref{thm:triangle-square} connects these conditions to vanishing in magnitude homology. In particular, for a finite connected locally weakly modular graph, vanishing off the diagonal in homological degree $2$ implies vanishing off the diagonal in all homological degrees.
Equivalently, within this class of graphs, quadraticity of the distance algebra over every field implies Koszulity over every field.

\subsection{Basis graph of matroids and even $\Delta$-matroids}\label{subsec:even-delta}

Matroids capture properties shared by families of spanning trees and families of vector space bases.
\begin{defi}
Let $\Omega$ be a finite set and $\mathcal{F}\subseteq 2^\Omega$ a nonempty family of subsets.
We call $M=(\Omega,\mathcal{F})$ a \emph{matroid} if, for every $x,y\in\mathcal{F}$ and $e\in x\setminus y$, there exists $f\in y\setminus x$ such that $(x\setminus\{e\})\cup\{f\}\in\mathcal{F}$.
The members of $\mathcal{F}$ are the \emph{bases} of $M$.
\end{defi}
For background on matroids, see \cite{Oxley2011}.
For integers $0\le \kappa\le m$, the pair $U_{\kappa,m}\coloneq([m],\binom{[m]}{\kappa})$ satisfies the basis exchange axiom and is therefore a matroid.
The matroid $U_{\kappa,m}$ is called the \emph{uniform matroid}.
In general, all members of $\mathcal{F}$ have the same cardinality \cite{Oxley2011}.

$\Delta$-matroids, also called symmetric matroids, were introduced in \cite{Bouchet1987}. They are also known as a type of Coxeter matroid \cite[Chapter 4]{BGW2003}.
The following definition of a $\Delta$-matroid generalizes the matroid basis exchange axiom.
Recall that for $x,y\subseteq\Omega$, their symmetric difference is denoted by $x\triangle y$.

\begin{defi}\label{def:deltamatroid}
Let $\Omega$ be a finite set and $\mathcal{F}\subseteq 2^\Omega$ a nonempty family of subsets.
We call $M=(\Omega,\mathcal{F})$ a \emph{$\Delta$-matroid} if, for every $x,y\in\mathcal{F}$ and $e\in x\triangle y$, there exists $f\in x\triangle y$ such that $x\triangle\{e,f\}\in\mathcal{F}$.
The members of $\mathcal{F}$ are the \emph{feasible sets} of $M$.

A $\Delta$-matroid is \emph{even} if $\abs{x}$ and $\abs{y}$ have the same parity for every $x,y\in\mathcal{F}$.
\end{defi}

Every matroid is an even $\Delta$-matroid.
Taking $\Omega=[m]$ and letting $\mathcal{F}$ be the vertex set of the halved $m$-cube $Q_m^+$ yields another even $\Delta$-matroid.

Matroid basis graphs are the $1$-skeletons of matroid base polytopes \cite{GGMS1987}. They were defined and characterized in \cite{Maurer1973,Maurer1973-2}.
This characterization was extended to basis graphs of even $\Delta$-matroids in \cite{Chepoi2007}.

\begin{defi}
For an even $\Delta$-matroid $M=(\Omega,\mathcal{F})$, its \emph{basis graph} $B(M)$ has vertex set $\mathcal{F}$. Distinct vertices $x,y\in\mathcal{F}$ are adjacent if and only if $\abs{x\triangle y}=2$.
\end{defi}

The basis graph of the uniform matroid $U_{\kappa,m}$ is isomorphic to the Johnson graph $J(m,\kappa)$.
More generally, the basis graph of a matroid is an induced subgraph of a Johnson graph, and the basis graph of an even $\Delta$-matroid is isomorphic to an induced subgraph of a halved cube.
When the feasible sets have odd cardinality, taking the symmetric difference with a fixed feasible set embeds the basis graph in the halved cube.

\begin{thm}\label{prop:even-delta-reducible}
Basis graphs of even $\Delta$-matroids are reducive and hence diagonal.
\end{thm}

For the rest of this subsection, let $M=(\Omega,\mathcal{F})$ be an even $\Delta$-matroid.

\begin{lemma}[\cite{Chepoi2007}]\label{lem:even-delta-dist}
For all $u,x\in \mathcal{F}$, we have $d(u,x)={\abs{u\triangle x}}/{2}$.
\end{lemma}

\begin{lemma}\label{lem:even-delta-distance}
For $w\in\Gamma(u)$, we have $T_u(w)=\{x\in\mathcal{F}\mid u\triangle w\subseteq u\triangle x\}$.
\end{lemma}

\begin{proof}
Since $\abs{u\triangle w}=2$, the distance formula in \zcref{lem:even-delta-dist} and inclusion-exclusion give $d(w,x)=\abs{(u\triangle w)\triangle(u\triangle x)}/2 = 1 + d(u,x) - \abs{(u\triangle w)\cap(u\triangle x)}$.
We have $x\in T_u(w)$ exactly when $\abs{(u\triangle w)\cap(u\triangle x)}=2$. The conclusion follows from $\abs{u\triangle w}=2$.
\end{proof}

For each $u$ and $\Lambda\subseteq\Omega$, put
\begin{equation}\label{eq:even-delta-N}
    \Gamma_\Lambda(u)\coloneq\{y\in\Gamma(u)\mid(u\triangle y)\cap\Lambda\ne\emptyset\}.
\end{equation}

\begin{lemma}\label{lem:even-delta-target}
For each $u$ and $\Lambda\subseteq\Omega$, we have $T_u(\Gamma_\Lambda(u))=\{x\in\mathcal{F}\mid(u\triangle x)\cap\Lambda\ne\emptyset\}$.
\end{lemma}

\begin{proof}
If $x\in T_u(\Gamma_\Lambda(u))$, then $x\in T_u(y)$ for some $y\in\Gamma_\Lambda(u)$.
By \zcref{lem:even-delta-distance}, we have $u\triangle y\subseteq u\triangle x$, so $(u\triangle x)\cap\Lambda\ne\emptyset$.
Conversely, suppose $(u\triangle x)\cap\Lambda\ne\emptyset$ and take $e\in(u\triangle x)\cap\Lambda$.
By \zcref{def:deltamatroid}, there exists $f\in u\triangle x$ such that $y=u\triangle\{e,f\}\in\mathcal{F}$.
If $f=e$, then the exchange would change the parity of the size of the feasible set. Thus $f\ne e$.
Since $y\triangle u=\{e,f\}\subseteq u\triangle x$, we have $y\in\Gamma_\Lambda(u)$ and $x\in T_u(y)$, proving the other inclusion.
\end{proof}

\begin{proof}[{Proof of \zcref{prop:even-delta-reducible}}]
To prove reducivity by \zcref{lem:basis-reduction}, we construct families of cliques satisfying \ref{item:basis-reduction-clique} and \ref{item:basis-reduction-partition}.
Put $\mathcal{C}_u\coloneq\{\emptyset\}\cup\{\Gamma_{\{e\}}(u)\mid e\in\Omega\}$.
Two distinct vertices of $\Gamma_{\{e\}}(u)$ can be written as $u\triangle\{e,f\}$ and $u\triangle\{e,g\}$. Their symmetric difference is $\{f,g\}$, so each $\Gamma_{\{e\}}(u)$ is a clique.
The identities $\Gamma_\Lambda(u)=\bigcup_{e\in\Lambda}\Gamma_{\{e\}}(u)$ and $\Gamma_\Omega(u)=\Gamma(u)$ give \ref{item:basis-reduction-clique}.

Write $C\in\mathcal{C}_u$ as $C=\Gamma_\Lambda(u)$, where $\Lambda$ is either empty or a singleton, and take any $w\in\Gamma(u)\setminus C$ with $u\triangle w=\{e,f\}$.
Since $w\notin C$, we have $(u\triangle w)\cap \Lambda=\{e,f\}\cap\Lambda=\emptyset$.
The identity
\begin{equation}
    (w\triangle x)\cap(\Lambda\cup\{e,f\})=(u\triangle x\triangle \{e,f\})\cap(\Lambda\cup\{e,f\})
    =\big((u\triangle x)\cap \Lambda\big)\cup\big(\{e,f\}\setminus (u\triangle x)\big)
\end{equation}
shows that $(w\triangle x)\cap(\Lambda\cup\{e,f\})=\emptyset$ holds exactly when $\{e,f\}\subseteq u\triangle x$ and $(u\triangle x)\cap\Lambda=\emptyset$.
Applying \zcref{lem:even-delta-target} at the vertex $w$ gives
\begin{align}
    \mathcal{F}\setminus T_w(\Gamma_{\Lambda\cup\{e,f\}}(w))
    &=\{x\in\mathcal{F}\mid(w\triangle x)\cap(\Lambda\cup\{e,f\})=\emptyset\}\\
    &=\{x\in\mathcal{F}\mid\{e,f\}\subseteq u\triangle x,\ (u\triangle x)\cap\Lambda=\emptyset\}=T_u^C(w).
\end{align}
For the last equality, \zcref{lem:even-delta-distance} shows that $\{e,f\}\subseteq u\triangle x$ is equivalent to $x\in T_u(w)$, and \zcref{lem:even-delta-target} shows that $(u\triangle x)\cap\Lambda=\emptyset$ is equivalent to $x\notin T_u(\Gamma_\Lambda(u))=T_u(C)$.
Put $R\coloneq\Gamma_{\Lambda\cup\{e,f\}}(w)=\bigcup_{g\in \Lambda\cup\{e,f\}}\Gamma_{\{g\}}(w)$, which is a union of members of $\mathcal{C}_w$.
Since $\mathcal{F}=T_u^C(w)\sqcup T_w(R)$, condition \ref{item:basis-reduction-partition} holds.

Thus \zcref{lem:basis-reduction} applies, and $B(M)$ is reducive.
Diagonality follows from \zcref{thm:reducible}.
\end{proof}

\begin{coro}[\cite{Koizumi2026}]\label{cor:matroid-basis-reducive}
Matroid basis graphs are diagonal.
\end{coro}

The basis graph of an even $\Delta$-matroid need not be weakly modular \cite{BC08}.

\subsection{Basis graphs of $q$-matroids}\label{subsubsec:q-mat}

A $q$-matroid is the $q$-analogue of a matroid. This notion was introduced in \cite{JurPel2018}. Its applications include rank-metric codes.

\begin{defi}[\cite{JurPel2018, CeJu2024}]\label{def:q-bases}
Let $U$ be a finite-dimensional $\mathbb{F}_q$-vector space, and let $\mathcal{L}(U)$ be the set of subspaces of $U$.
If a nonempty family $\mathcal{F}\subseteq\mathcal{L}(U)$ of subspaces of $U$ satisfies the following three conditions,
we call $M=(U,\mathcal{F})$ a \emph{$q$-matroid}.
\begin{itemize}
    \labelitem{B1}{item:q-basis-inclusion}
    For all $u,x\in\mathcal{F}$, if $u\leq x$, then $u=x$.
    \labelitem{B2}{item:q-basis-exchange}
    For all $u,x\in\mathcal{F}$ and every hyperplane $h$ of $u$ containing $u\cap x$,
    there exists a line $\nu\leq x$ such that $h+\nu\in\mathcal{F}$.
    Here a line is a one-dimensional subspace, and a hyperplane of $u$ is a subspace of $u$ of codimension $1$.
    \labelitem{B3}{item:q-basis-sum}
    For any $W_0,W_1\leq U$, if $u,x\in\mathcal{F}$ satisfy $\dim(u\cap W_0)=\max_{y\in\mathcal{F}}\dim(y\cap W_0)$ and $\dim(x\cap W_1)=\max_{y\in\mathcal{F}}\dim(y\cap W_1)$, then there exists $z\in\mathcal{F}$ such that $\dim\bigl(z\cap(W_0+W_1)\bigr)=\max_{y\in\mathcal{F}}\dim\bigl(y\cap(W_0+W_1)\bigr)$ and $z\cap(W_0+W_1)\leq(u\cap W_0)+(x\cap W_1)$.
\end{itemize}
The members of $\mathcal{F}$ are the \emph{bases} of $M$, and we write $\mathcal{F}=\mathcal{F}(M)$.
\end{defi}

The bases of a $q$-matroid have a common dimension, just as the bases of an ordinary matroid have a common cardinality.
Indeed, suppose there are $u,x\in\mathcal{F}$ with $\dim u>\dim x$.
Take a hyperplane $h$ of $u$ containing $u\cap x$. Applying \ref{item:q-basis-exchange} to $u,x$ gives $z=h+\nu\in\mathcal{F}$ with $\nu\leq x$.
By \ref{item:q-basis-inclusion}, we have $\nu\nleq h$, so $\dim z=\dim u$ and $\dim(x\cap z)=\dim(u\cap x)+1$.
Repeat the exchange, increasing the dimension of the intersection at each step, until we obtain a basis $z$ with $x\leq z$. The containment $x\leq z$ and the inequality $\dim z=\dim u>\dim x$ contradict \ref{item:q-basis-inclusion}.
We call this common dimension the \emph{rank} of $M$.

The following construction relates a $q$-matroid to an ordinary matroid.
For a rank-$\kappa$ $q$-matroid $M=(U,\mathcal{F})$, put $\Omega\coloneq\{\nu\le U\mid\dim\nu=1\}$ and write $\langle S\rangle\coloneq\sum_{\nu\in S}\nu$ for $S\subseteq\Omega$.
Then
\begin{equation}
    \widehat {\mathcal{F}}\coloneq\{\widehat u\subseteq\Omega\mid\abs{\widehat u}=\kappa,\ \langle\widehat u\rangle\in\mathcal{F}(M)\}
\end{equation}
is the family of bases of a matroid $\widehat{M}=(\Omega,\widehat{\mathcal{F}})$ \cite[Definition~14, Proposition~15]{JPV23}.
We call $\widehat{M}$ the \emph{projectivization matroid} of $M$.

\begin{defi}\label{def:q-basis-graph}
Let $M=(U,\mathcal{F})$ be a $q$-matroid of rank $\kappa$.
Its \emph{basis graph} $B(M)$ has vertex set $\mathcal{F}$,
with distinct vertices $u,w\in\mathcal{F}$ adjacent if and only if $\dim(u\cap w)=\kappa-1$.
\end{defi}

By definition, every $q$-matroid basis graph is an induced subgraph of a Grassmann graph. Conversely, every Grassmann graph is isomorphic to the basis graph of a suitable $q$-matroid.

\begin{prop}\label{lem:q-distance}
Consider the basis graph $B(M)$ of a rank-$\kappa$ $q$-matroid $M=(U,\mathcal{F})$.
For all $u,x\in \mathcal{F}$, we have $d(u,x)=\kappa-\dim(u\cap x)$.
In particular, $B(M)$ is connected.
\end{prop}

\begin{proof}
Extend a linear basis of the intersection to a basis of each subspace.
The lines spanned by the resulting basis vectors give bases $\widehat u,\widehat x\in\widehat {\mathcal{F}}$ of the projectivization matroid $\widehat{M}$ of $M$ such that $\langle\widehat u\rangle=u$, $\langle\widehat x\rangle=x$, and $\abs{\widehat u\cap\widehat x}=\dim(u\cap x)$.
By \zcref{lem:even-delta-dist}, we have $d_{B(\widehat{M})}(\widehat u,\widehat x)=\kappa-\dim(u\cap x)$.
Adjacent vertices of $B(\widehat{M})$ are bases that share $\kappa-1$ lines. The spans of the two bases either coincide or are adjacent in $B(M)$.
Replacing each vertex of a shortest path from $\widehat u$ to $\widehat x$ by its span and removing consecutive repetitions gives a walk from $u$ to $x$.
Thus $B(M)$ is connected and $d(u,x)\le \kappa-\dim(u\cap x)$.
Conversely, moving along one edge increases the dimension of the intersection with $x$ by at most $1$, so $d(u,x)\ge \kappa-\dim(u\cap x)$.
\end{proof}

\begin{thm}\label{thm:q-reducive}
Basis graphs of $q$-matroids are reducive and hence diagonal.
\end{thm}

Let $M=(U,\mathcal{F})$ be a rank-$\kappa$ $q$-matroid.
Fix a nondegenerate symmetric bilinear form on $U$, and write $W^\perp$ for the orthogonal complement of a subspace $W$.
The family $\mathcal{F}^*=\{u^\perp\mid u\in\mathcal{F}\}$ is the family of bases of the dual $q$-matroid $M^*$. The map $u\mapsto u^\perp$ gives an isomorphism between $B(M)$ and $B(M^*)$ \cite[Theorems~42 and~45]{JurPel2018}.

For each $u$ and $\Lambda\subseteq\mathcal{L}(U)$, put
\begin{equation}\label{eq:q-N}
    N_u(\Lambda)\coloneq\{y\in\Gamma(u)\mid\text{there exists } W\in\Lambda\text{ such that } W\le u\cap y\ \text{or}\ u+y\le W\}.
\end{equation}
Taking orthogonal complements gives $N_{u^\perp}(\{W^\perp\})=\{y^\perp\mid y\in N_u(\{W\})\}$.
Since $u\mapsto u^\perp$ is a graph isomorphism, for every $S\subseteq\Gamma(u)$, we have $T_{u^\perp}(\{y^\perp\mid y\in S\})=\{x^\perp\mid x\in T_u(S)\}$.
Here $N_{u^\perp}$ and $T_{u^\perp}$ are defined in $B(M^*)$.

By \zcref{lem:q-distance}, for $w\in\Gamma(u)$ we have $x\in T_u(w)$ exactly when $\dim(x\cap w)=\dim(x\cap u)+1$.
If $x\in T_u(w)$, then $u\cap w$ is a hyperplane of $w$, so $\dim(x\cap u\cap w)\ge\dim(x\cap w)-1=\dim(x\cap u)$. Hence $x\cap u\le w$.

\begin{lemma}\label{lem:q-complement-target}
For every $h\le u\le s\le U$, we have
\begin{align}
    T_u(N_u(\{h\}))&=\{x\in\mathcal{F}\mid h+(x\cap u)\ne u\},\label{eq:q-lower-target}\\
    T_u(N_u(\{s\}))&=\{x\in\mathcal{F}\mid x\cap s\ne x\cap u\}.\label{eq:q-upper-target}
\end{align}
Therefore,
\begin{equation}\label{eq:q-complement-target}
    \mathcal{F}\setminus T_u(N_u(\{h,s\}))=\{x\in\mathcal{F}\mid h+(x\cap u)=u,\ x\cap s=x\cap u\}.
\end{equation}
\end{lemma}

\begin{proof}
Since $h\le u\le s$, \eqref{eq:q-N} gives
\begin{equation}
    N_u(\{h\})=\{y\in\Gamma(u)\mid h\le y\},\qquad N_u(\{s\})=\{y\in\Gamma(u)\mid y\le s\}.
\end{equation}
The set $N_u(\{h\})$ consists of the neighboring bases containing $h$, and $N_u(\{s\})$ consists of the neighboring bases contained in $s$.

We first prove \eqref{eq:q-lower-target}.
If $x\in T_u(N_u(\{h\}))$, there exists $y\in\Gamma(u)$ such that $h\leq u\cap y$ and $x\in T_u(y)$.
Since $x\cap u\le y$, we have $h+(x\cap u)\le u\cap y\ne u$.

Conversely, take $x\in\mathcal{F}$ with $h+(x\cap u)\ne u$.
Choose a hyperplane $h'$ of $u$ containing $h+(x\cap u)$. By \ref{item:q-basis-exchange}, there is a line $\nu\le x$ such that $y=h'+\nu\in\mathcal{F}$.
If $\nu\le u$, then $\nu\le x\cap u\le h'$, which contradicts $\dim y=\kappa$.
Hence $\nu\nleq u$.
Therefore $u\cap y=h'$ and $y\in N_u(\{h\})$.
Also, since $\nu\le x$ and $x\cap h'=x\cap u$, we have
\begin{equation}
    x\cap y=x\cap(h'+\nu)=(x\cap h')+\nu=(x\cap u)+\nu.
\end{equation}
Since $\nu\nleq u$, we have $\dim(x\cap y)=\dim(x\cap u)+1$. By \zcref{lem:q-distance}, we obtain $x\in T_u(y)$ and hence \eqref{eq:q-lower-target}.

Next, since $u\le s$, we have $s\cap(x+u)=u+(x\cap s)$.
Applying \eqref{eq:q-lower-target} to $s^\perp\le u^\perp$ in the dual gives
\begin{align}
    \mathcal{F}\setminus T_u(N_u(\{s\}))
    &=\{x\in\mathcal{F}\mid s^\perp+(x^\perp\cap u^\perp)=u^\perp\}\\
    &=\{x\in\mathcal{F}\mid s\cap(x+u)=u\}\\
    &=\{x\in\mathcal{F}\mid x\cap s=x\cap u\}.
\end{align}
This proves \eqref{eq:q-upper-target}.
Since $N_u(\{h,s\})=N_u(\{h\})\cup N_u(\{s\})$ and target sets preserve unions by \zcref{prop:target-union}, we obtain \eqref{eq:q-complement-target}.
\end{proof}

\begin{lemma}\label{lem:q-target}
Let $w\in\Gamma(u)$ and let $W\le U$ be a subspace such that $W\le u$ or $u\le W$. Assume $w\notin N_u(\{W\})$.
Put $R\coloneq N_w(\{w\cap W,w+W\})$. Then
\begin{equation}
    \mathcal{F}\setminus T_w(R)=T_u^{N_u(\{W\})}(w).
\end{equation}
\end{lemma}

\begin{proof}
Since $(w\cap W)^\perp=w^\perp+W^\perp$ and $(w+W)^\perp=w^\perp\cap W^\perp$, duality reduces the proof to the case $W\le u$.

Since $w\notin N_u(\{W\})$, we have $W\nleq w$.
Since $W\le u$ and $w$ is a hyperplane of $u+w$, we obtain $w+W=u+w$.
Thus $R$ is the union of two sets: the bases adjacent to $w$ that contain $w\cap W$, and the bases adjacent to $w$ that are contained in $u+w$.
By \zcref{lem:q-complement-target},
\begin{equation}\label{eq:q-transfer-conditions}
    \mathcal{F}\setminus T_w(R)=\{x\in\mathcal{F}\mid(w\cap W)+(x\cap w)=w,\ x\cap(u+w)=x\cap w\}.
\end{equation}

We first show $\mathcal{F}\setminus T_w(R)\subseteq T_u(w)$.
Take $x\in\mathcal{F}\setminus T_w(R)$.
The identity $x\cap(u+w)=x\cap w$ in \eqref{eq:q-transfer-conditions} gives $x\cap u\le x\cap w$.
Since $(w\cap W)+(x\cap w)=w$ and $W\le u$, we also have
\begin{equation}
    w=(w\cap W)+(x\cap w)\leq (u\cap w)+(x\cap w)\leq w,
\end{equation}
so $(u\cap w)+(x\cap w)=w$.
Therefore,
\begin{equation}
    \dim(x\cap w)-\dim(x\cap u)=\dim w-\dim(u\cap w)=1,
\end{equation}
and hence $x\in T_u(w)$.

For any $x\in T_u(w)$, the fact that $u$ is a hyperplane of $u+w$ gives
\begin{equation}
    \dim(x\cap(u+w))\le\dim(x\cap u)+1=\dim(x\cap w).
\end{equation}
Since $x\cap w\le x\cap(u+w)$, we obtain $x\cap(u+w)=x\cap w$.

Next, assuming $x\in T_u(w)$, we show that $(w\cap W)+(x\cap w)=w$ and $W+(x\cap u)=u$ are equivalent.
Since $(w\cap W)+(x\cap w)=w$ and $w+W=u+w$, we have $W+(x\cap w)=u+w$.
Since $W\le u$ and $x\cap u\le w$, we get
\begin{equation}
    u=(W+(x\cap w))\cap u=W+(x\cap u\cap w)=W+(x\cap u).
\end{equation}

Conversely, take $x\in T_u(w)$ with $W+(x\cap u)=u$.
Since $x\cap u\le x\cap w$ and $\dim(x\cap w)=\dim(x\cap u)+1$, we have $x\cap w\nleq u$.
Thus
\begin{align}
    W+(x\cap w)&=u+(x\cap w)=u+w,\\
    (W\cap w)+(x\cap w)&=(W+(x\cap w))\cap w=w.
\end{align}
Therefore \eqref{eq:q-lower-target} gives
\begin{align}
    \mathcal{F}\setminus T_w(R)
    &=\{x\in T_u(w)\mid(w\cap W)+(x\cap w)=w\}\\
    &=\{x\in T_u(w)\mid W+(x\cap u)=u\}\\
    &=T_u(w)\setminus T_u(N_u(\{W\}))=T_u^{N_u(\{W\})}(w).\qedhere
\end{align}
\end{proof}

\begin{proof}[{Proof of \zcref{thm:q-reducive}}]
For each $u$, put
\begin{align}
    \mathcal{C}_u\coloneq{}&\{\emptyset\}\cup\{N_u(\{h\})\mid h\le u,\ \dim h=\kappa-1\}\\
    &{}\cup\{N_u(\{s\})\mid u\le s,\ \dim s=\kappa+1\}.
\end{align}
Two distinct bases containing a given hyperplane $h$ intersect in $h$, so each $N_u(\{h\})$ is a clique. By duality, each $N_u(\{s\})$ is also a clique.
Every $y\in\Gamma(u)$ belongs to $N_u(\{u\cap y\})$, so $\Gamma(u)$ is the union of the members of $\mathcal{C}_u$.
Thus \ref{item:basis-reduction-clique} holds.

We next express $N_u(\{W\})$ as a union of members of $\mathcal{C}_u$.
If $W\le u$, then
\begin{equation}
    N_u(\{W\})=\bigcup_{\substack{W\le h\le u\\ \dim h=\kappa-1}}N_u(\{h\}).
\end{equation}
For $y\in N_u(\{W\})$, take $h=u\cap y$. Then $W\le h\le u$, $\dim h=\kappa-1$, and $y\in N_u(\{h\})$.
Conversely, if $W\le h\le u$, then $N_u(\{h\})\subseteq N_u(\{W\})$.
By duality, for $u\le W$ we have
\begin{equation}
    N_u(\{W\})=\bigcup_{\substack{u\le s\le W\\ \dim s=\kappa+1}}N_u(\{s\}).
\end{equation}
If neither of $W$ and $u$ contains the other, then $N_u(\{W\})=\emptyset$.
Hence every $N_u(\Lambda)$ is a union of members of $\mathcal{C}_u$.
Conversely, since $\emptyset=N_u(\{u\})$, every union of members of $\mathcal{C}_u$ can be written as $N_u(\Lambda)$ for a suitable indexing set.

Finally, we check \ref{item:basis-reduction-partition}.
Write any $C\in\mathcal{C}_u$ as $C=N_u(\{W\})$ for some $W$ with $W\le u$ or $u\le W$.
If $C=\emptyset$, take $W=u$.
For $w\in\Gamma(u)\setminus C$, \zcref{lem:q-target} gives
\begin{equation}
    \mathcal{F}=T_u^C(w)\sqcup T_w(R),\qquad R=N_w(\{w\cap W,w+W\}).
\end{equation}
Since $R=N_w(\{w\cap W,w+W\})$ is a union of members of $\mathcal{C}_w$, condition \ref{item:basis-reduction-partition} also holds.

Thus \ref{item:basis-reduction-clique} and \ref{item:basis-reduction-partition} hold. By \zcref{lem:basis-reduction}, the families
\begin{equation}\label{eq:q-family}
    \mathcal{R}_u=\{N_u(\Lambda)\cup S_0\mid\Lambda\subseteq\mathcal{L}(U),\ S_0\subseteq\Gamma(u)\text{ is a clique}\}
\end{equation}
show that $B(M)$ is reducive.
Diagonality follows from \zcref{thm:reducible}.
\end{proof}

\begin{rema}
The conclusion remains valid under \ref{item:q-basis-inclusion} and \ref{item:q-basis-exchange} alone. The projectivization and dual constructions are used above only to simplify the proof; the required distance and target-set identities can also be established directly.
Thus, if a nonempty family $\mathcal{F}$ of subspaces satisfies \ref{item:q-basis-inclusion} and \ref{item:q-basis-exchange}, the same argument shows that the subgraph of the Grassmann graph induced by $\mathcal{F}$ is reducive.
\end{rema}

Put $U=\mathbb{F}_q^n$, take $j\in\mathbb Z_{>0}$, and let $\Xi\in\mathbb{F}_{q^j}^{\kappa\times n}$ be a matrix of rank $\kappa$.
Regard $\Xi$ as a linear map. Write $\Xi|_{W\otimes_{\mathbb{F}_q} \mathbb{F}_{q^j}}$ for the restriction of $\Xi$ to the scalar extension $W\otimes_{\mathbb{F}_q} \mathbb{F}_{q^j}$ of a subspace $W\le U$.
Then
\begin{equation}
    \mathcal{F}(M[\Xi])\coloneq\{W\le U\mid\dim_{\mathbb{F}_{q}} W=\kappa,\ \rank_{\mathbb{F}_{q^j}} \Xi |_{W\otimes_{\mathbb{F}_q}\mathbb{F}_{q^j}}=\kappa\}
\end{equation}
is the family of bases of a $q$-matroid \cite[Section~5]{JurPel2018}.
A $q$-matroid $M$ is \emph{representable} over $\mathbb{F}_{q^j}$ if there exists $\Xi \in\mathbb{F}_{q^j}^{\kappa\times n}$ such that $\mathcal{F}(M)=\mathcal{F}(M[\Xi])$. We call $\Xi$ a \emph{representation} of the given $q$-matroid.

\begin{defi}[{\cite[Section~9]{BCN1989}}]
Let $q$ be a prime power, and let $D, m$ be integers satisfying $m \ge 2D \ge 1$.
The \emph{bilinear forms graph} $H_q(D, m-D)$ has vertex set $\mathbb{F}_q^{(m-D)\times D}$, consisting of all $(m-D)\times D$ matrices over $\mathbb{F}_q$.
There are $q^{D(m-D)}$ vertices.
Two vertices $u$ and $w$ are adjacent if and only if $\rank(u-w)=1$.
\end{defi}

Bilinear forms graphs can also be described using subspaces that have trivial intersection with a fixed subspace~\cite[Section~3.1]{vDKT16}.
Let $M$ be the $q$-matroid represented by $\begin{pmatrix}I_D&O_{D\times (m-D)}\end{pmatrix}\in\mathbb{F}_q^{D\times m}$.
In the direct sum decomposition $\mathbb{F}_q^m=\mathbb{F}_q^D\oplus \mathbb{F}_q^{m-D}$, the matrix $\begin{pmatrix}I_D&O_{D\times (m-D)}\end{pmatrix}$ represents the projection onto $\mathbb{F}_q^D\le \mathbb{F}_q^m$.
Thus the bases of $M$ are the subspaces $W\le \mathbb{F}_q^m$ satisfying $\dim W=D$ and $W\cap \mathbb{F}_q^{m-D}=0$.
Sending each $\gamma\in \mathbb{F}_q^{(m-D)\times D}$ to $W_\gamma\coloneq\{(s,\gamma(s))\mid s\in \mathbb{F}_q^D\}\le \mathbb{F}_q^{m}$ gives a bijection from $\mathbb{F}_q^{(m-D)\times D}$ to the bases of $M$.
The projection $\mathbb{F}_q^m\to \mathbb{F}_q^D$ onto the first component restricts to an isomorphism from $W_\gamma\cap W_{\gamma'}$ to the kernel of $\gamma-\gamma'$, so
\begin{equation}
    \dim(W_\gamma\cap W_{\gamma'})=D-\rank(\gamma-\gamma').
\end{equation}
Thus $\gamma\mapsto W_\gamma$ is an isomorphism from the bilinear forms graph, whose adjacency condition is $\rank(\gamma-\gamma')=1$, to $B(M)$.
Applying \zcref{thm:reducible, thm:q-reducive} to this isomorphism gives the following.

\begin{prop}\label{prop:bilinear-q-basis}
A bilinear forms graph is isomorphic to the basis graph of a $q$-matroid representable over $\mathbb{F}_q$.
In particular, bilinear forms graphs are reducive and hence diagonal.
\end{prop}

\begin{prop}\label{exam:bilinear-not-wmg}
Bilinear forms graphs of diameter at least $3$ are not weakly modular.
In particular, basis graphs of $q$-matroids need not be weakly modular.
\end{prop}

\begin{proof}
The construction in \zcref{prop:bilinear-q-basis} and the distance formula in \zcref{lem:q-distance} give $d(u,w)=\rank(u-w)$ for bilinear forms graphs.
Over $\mathbb{F}_q$, consider the following matrices, padded with zeros if necessary:
\begin{equation}
    x=O,\qquad z=I_3,\qquad
    u=\begin{pmatrix}1&0&0\\0&1&0\\0&0&0\end{pmatrix},\qquad
    w=\begin{pmatrix}1&-1&1\\0&1&0\\0&1&0\end{pmatrix}.
\end{equation}
Computing the ranks of their differences gives
\begin{equation}
    d(x,u)=d(x,w)=d(u,w)=2,\qquad d(x,z)=3,\qquad d(u,z)=d(w,z)=1,
\end{equation}
so the hypotheses of \ref{item:quadrangle-condition} hold.

If there were a vertex $y\in\Gamma(x)\cap\Gamma(u)\cap\Gamma(w)$, then $\rank y=\rank(u-y)=\rank(w-y)=1$.
We write $\operatorname{Col}$ for the column space of a matrix and $\operatorname{Row}$ for its row space.
Since $\rank u=\rank y+\rank(u-y)=2$ and $\rank w=\rank y+\rank(w-y)=2$, we have $\operatorname{Col}(y)\le\operatorname{Col}(u)\cap\operatorname{Col}(w)$ and $\operatorname{Row}(y)\le\operatorname{Row}(u)\cap\operatorname{Row}(w)$.
The intersection $\operatorname{Col}(u)\cap\operatorname{Col}(w)$ consists of column vectors whose entries are all $0$ except possibly in the first coordinate. Similarly, $\operatorname{Row}(u)\cap\operatorname{Row}(w)$ consists of row vectors whose entries are all $0$ except possibly in the second coordinate. Thus only entry $(1,2)$ of $y$ can be nonzero.
But the upper-left $2\times2$ submatrix of $u-y$ has determinant $1$, so $\rank(u-y)\ge2$, a contradiction.
Therefore $\Gamma(u)\cap\Gamma(w)\cap\Gamma(x)=\emptyset$, and \ref{item:quadrangle-condition} fails.
The last claim follows from \zcref{prop:bilinear-q-basis}.
\end{proof}

\subsection{Meshed graphs}\label{subsec:nonred-mesh}
We conclude the section with a class that includes weakly modular graphs and the basis graphs of even $\Delta$-matroids and $q$-matroids.
\begin{defi}[{\cite[Section~2.2.4.1]{CCHO20}}]
A connected graph $G=(V,E)$ is \emph{meshed} if, for all $x,y,z\in V$ with $d(y,z)=2$, there exists $u\in\Gamma(y)\cap\Gamma(z)$ such that
\begin{equation}
2d(x,u)\le d(x,y)+d(x,z).
\end{equation}
\end{defi}

\begin{prop}[{\cite[Lemma~2.22 and Section~2.2.4]{CCHO20}}]\label{prop:wmg-matroid-meshed}
Weakly modular graphs and basis graphs of even $\Delta$-matroids are meshed.
\end{prop}

\begin{prop}\label{prop:q-basis-meshed}
The basis graph of a $q$-matroid is meshed.
\end{prop}

\begin{lemma}\label{lem:q-associated-matroid}
Let $M=(U,\mathcal{F})$ be a $q$-matroid, and let $\widehat{M}=(\Omega,\widehat{\mathcal{F}})$ be its projectivization matroid.
For any $x,y,z\in\mathcal{F}$, there exist $\widehat x,\widehat y,\widehat z\in\widehat {\mathcal{F}}$ such that $\langle\widehat x\rangle=x$, $\langle\widehat y\rangle=y$, $\langle\widehat z\rangle=z$, and $\abs{\widehat u\cap\widehat w}=\dim(u\cap w)$ for $u,w\in\{x,y,z\}$.
\end{lemma}

\begin{proof}
Choose a linear basis of $x\cap y\cap z$ and extend it to linear bases of $x\cap y$, $y\cap z$, and $z\cap x$.
The two bases of the pairwise intersections contained in $x$ have in common only the chosen basis of $x\cap y\cap z$, and their union is linearly independent, so it extends to a linear basis of $x$.
Repeat the construction for $y,z$, and replace each basis vector with the line it spans to obtain $\widehat x,\widehat y,\widehat z$.
We have $\widehat x,\widehat y,\widehat z\in \widehat {\mathcal{F}}$, and the lines common to any two of these sets are the lines obtained from the basis chosen for the corresponding intersection.
\end{proof}

\begin{proof}[{Proof of \zcref{prop:q-basis-meshed}}]
Let $M=(U,\mathcal{F})$ be a $q$-matroid, and take $x,y,z\in \mathcal{F}$ with $d(y,z)=2$.
By \zcref{lem:q-associated-matroid}, in the projectivization $(\Omega,\widehat{\mathcal{F}})$ of $M$ we can choose $\widehat x,\widehat y,\widehat z\in\widehat{\mathcal{F}}$ that span $x,y,z$, respectively, so that the number of shared lines in each pair is the dimension of the corresponding intersection.
The distance formulas for the two basis graphs (\zcref{lem:even-delta-dist, lem:q-distance}) give $d_{B(\widehat{M})}(\widehat{w},\widehat{w'})=d(w,{w'})$ for $w,w'\in\{x,y,z\}$.
By \zcref{prop:wmg-matroid-meshed}, the matroid basis graph is meshed, so $\widehat y,\widehat z$ have a common neighbor $\widehat u$ such that $2d_{B(\widehat{M})}(\widehat x,\widehat u)\le d(x,y)+d(x,z)$.
Set $u\coloneq\langle\widehat u\rangle$. Since $\langle\widehat x\cap\widehat u\rangle\le x\cap u$, we have $d(x,u)\le d_{B(\widehat{M})}(\widehat x,\widehat u)$, and similarly $d(y,u),d(z,u)\le1$.
Since $d(y,z)=2$, we have $u\ne y,z$. Thus $u\in\Gamma(y)\cap\Gamma(z)$ and
\begin{equation}
    2d(x,u)\le2d_{B(\widehat{M})}(\widehat x,\widehat u)\le d(x,y)+d(x,z).
\end{equation}
Hence, $B(M)$ is meshed.
\end{proof}

Thus every graph in the three classes above is both meshed and reducive. The next example shows that a meshed graph need not be diagonal.

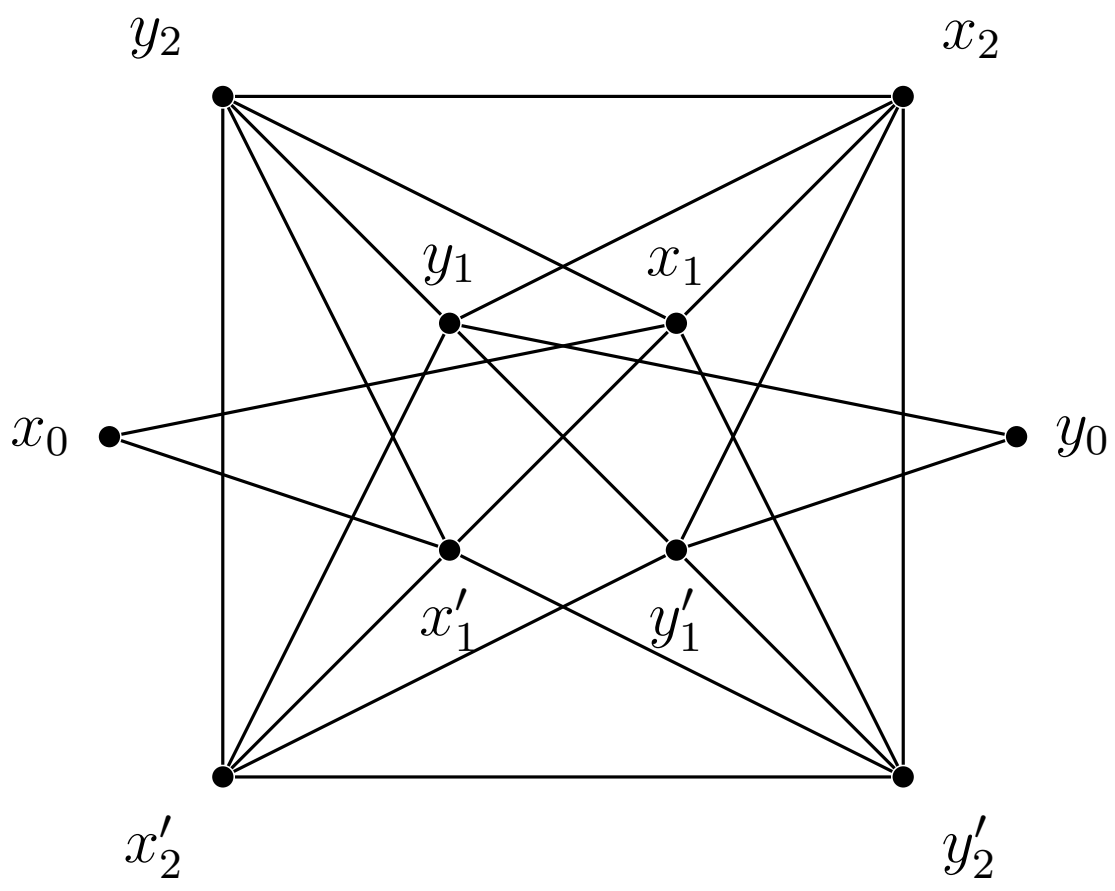
\begin{figure}[htbp]
\centering
\begin{tikzpicture}[
    vertex/.style={circle,fill=black,minimum size=4pt,inner sep=0pt},
    line width=0.6pt,
    scale=1.5
]
    \node[vertex] (x0)  at (-2,0) {};
    \node[vertex] (y0)  at (2,0) {};
    \node[vertex] (x1)  at (0.5,0.5) {};
    \node[vertex] (x1p) at (-0.5,-0.5) {};
    \node[vertex] (y2)  at (-1.5,1.5) {};
    \node[vertex] (y2p) at (1.5,-1.5) {};
    \node[vertex] (x2)  at (1.5,1.5) {};
    \node[vertex] (x2p) at (-1.5,-1.5) {};
    \node[vertex] (y1p)  at (0.5,-0.5) {};
    \node[vertex] (y1) at (-0.5,0.5) {};

    \foreach \from/\to in {
        x0/x1, x0/x1p, x1/x1p, x1/x2, x1p/x2p,
        y2/y1, y2p/y1p, y1/y1p, y0/y1, y0/y1p}{
        \draw (\from)--(\to);
    }
    \foreach \from in {x1,x1p}{
        \foreach \to in {y2,y2p}{
            \draw (\from)--(\to);
        }
    }
    \foreach \from in {x2,x2p}{
        \foreach \to in {y2,y2p,y1,y1p}{
            \draw (\from)--(\to);
        }
    }

    \foreach \point/\position/\labeltext in {
        x0/left/x_0, y0/right/y_0,
        x1/above/x_1, x1p/below/x'_1,
        y2/above left/y_2, y2p/below right/y'_2,
        x2/above right/x_2, x2p/below left/x'_2,
        y1/above/y_1, y1p/below/y'_1
    }{
        \node[\position=3pt]
            at (\point) {$\labeltext$};
    }
\end{tikzpicture}
\caption{A meshed graph with $10$ vertices and $22$ edges.}
\label{fig:meshed-ten}
\end{figure}

\begin{prop}\label{exam:meshed-not-diagonal}
Consider the graph $G=(V,E)$ shown in \zcref{fig:meshed-ten}, defined by
\begin{align}
V&=\{x_0,y_0,x_1,x'_1,x_2,x'_2,y_2,y'_2,y_1,y'_1\},\\
E&=\{x_0x_1,x_0x'_1,x_1x'_1,x_1x_2,x'_1x'_2,y_2y_1,y'_2y'_1,y_1y'_1,y_0y_1,y_0y'_1\}\\
 &\qquad\cup\{uw\mid u\in\{x_1,x'_1\},\ w\in\{y_2,y'_2\}\}\\
 &\qquad\cup\{uw\mid u\in\{x_2,x'_2\},\ w\in\{y_2,y'_2,y_1,y'_1\}\}.
\end{align}
The graph $G$ is meshed, its triangle-square complex $\mathsf{CW}_{\triangle\square}(G)$ is simply connected, and $\MH_{3,4}(G;x_0,y_0)\ne0$.
In particular, $G$ is neither diagonal nor reducive.
\end{prop}

\begin{proof}
To show that the graph is meshed, it suffices to prove the following inequality for every $u,w\in V$ with $d(u,w)=2$ and every $z\in V$:
\begin{equation}
2\min_{z'\in\Gamma(u)\cap\Gamma(w)}d(z,z')
\le d(z,u)+d(z,w).
\end{equation}
The following permutations preserve the edge set and hence are automorphisms of $G$:
\begin{equation}
(x_1\,x'_1)(x_2\,x'_2),\qquad
(y_2\,y'_2)(y_1\,y'_1),\qquad
(x_0\,y_0)(x_1\,y_1)(x'_1\,y'_1)(x_2\,y_2)(x'_2\,y'_2).
\end{equation}
By these automorphisms, it suffices to consider the five pairs at distance $2$ listed in the table below. Each entry is the value of
$d(z,u)+d(z,w)-2\min_{z'\in\Gamma(u)\cap\Gamma(w)}d(z,z')$.
\begin{equation}
\begin{array}{c|rrrrrrrrrr}
\{u,w\}\backslash z
 &x_0&y_0&x_1&x'_1&x_2&x'_2&y_2&y'_2&y_1&y'_1\\ \hline
\{x_0,x_2\} &0&0&2&1&0&0&1&1&0&0\\
\{x_0,y_2\} &0&0&2&2&1&1&0&2&0&1\\
\{x_1,y_1\} &0&0&0&1&2&1&2&1&0&1\\
\{x_1,x'_2\}&1&1&0&2&1&0&2&2&1&1\\
\{x_2,x'_2\}&0&2&1&1&0&0&2&2&2&2
\end{array}
\end{equation}
All entries are nonnegative, so $G$ is meshed.

Next, we show that $\MH_{3,4}(G;x_0,y_0)\ne0$. The sets $\Gamma(x_0)=\{x_1,x'_1\}$ and $\Gamma(y_0)=\{y_1,y'_1\}$ are disjoint, and no edge joins them. Thus there is no path between the endpoints of length less than four. The path $(x_0,x_1,x_2,y_1,y_0)$ has length $4$, which proves $d(x_0,y_0)=4$. Set
\begin{equation}
\xi=(x_0,x_1,y_2,y_0)-(x_0,x_1,y'_2,y_0)-(x_0,x'_1,y_2,y_0)+(x_0,x'_1,y'_2,y_0)
\in\MC_{3,4}(G;x_0,y_0).
\end{equation}
For $u\in\{x_1,x'_1\}$ and $w\in\{y_2,y'_2\}$, we have
\begin{equation}
\partial(x_0,u,w,y_0)=-(x_0,w,y_0)+(x_0,u,y_0)
\end{equation}
and hence $\partial\xi=0$.

We show that $\xi$ is not a boundary. Since $d(x_0,y_0)=4$, the basis of $\MC_{4,4}(G;x_0,y_0)$ consists of the following shortest paths of length $4$ from $x_0$ to $y_0$:
\begin{align}
&(x_0,x_1,x_2,w,y_0),\quad(x_0,x'_1,x'_2,w,y_0)
&& (\text{if }w\in\{y_1,y'_1\}),\\
&(x_0,u,y_2,y_1,y_0),\quad(x_0,u,y'_2,y'_1,y_0)
&& (\text{if }u\in\{x_1,x'_1\}).
\end{align}
Consider the coefficients of
\begin{equation}
p_1=(x_0,x_1,y_2,y_0),\qquad p_2=(x_0,x_1,y_1,y_0),\qquad p_3=(x_0,x_2,y_1,y_0)
\end{equation}
in the boundaries of these paths.
The only basis elements whose boundaries have a nonzero coefficient for at least one of $p_1, p_2, p_3$ are
$(x_0,x_1,x_2,y_1,y_0)$ and $(x_0,x_1,y_2,y_1,y_0)$. The two boundaries are
\begin{align}
\partial(x_0,x_1,y_2,y_1,y_0)
 &=-p_1+p_2-(x_0,y_2,y_1,y_0),\\
\partial(x_0,x_1,x_2,y_1,y_0)
 &=p_2-p_3-(x_0,x_1,x_2,y_0).
\end{align}
Set $P\coloneq P_{3,4}(G;x_0,y_0)\setminus\{p_1,p_2,p_3\}$ and write $\zeta=a_1p_1+a_2p_2+a_3p_3+\sum_{p\in P}a_pp$. The two boundary formulas show that $a_1+a_2+a_3=0$ whenever $\zeta\in \Im(\partial)$.
For $\xi$, the sum of the coefficients of $p_1,p_2,p_3$ is $1$.
Thus, $\xi\notin\Im(\partial)$.
Hence $\MH_{3,4}(G;x_0,y_0)\ne0$. The group $\MH_{3,4}(G;x_0,y_0)$ is a direct summand of $\MH_{3,4}(G)$, so $G$ is not diagonal and, by \zcref{thm:reducible}, is not reducive.

Finally, we show that $\mathsf{CW}_{\triangle\square}(G)$ is simply connected.
Each of $x_0$ and $y_0$ belongs to exactly one triangle and has no other incident edges.
Thus, the clique complexes of $G$ and $G-\{x_0,y_0\}$ are homotopy equivalent.
The clique complex of $G-\{x_0,y_0\}$ consists of two closed disks.
The first disk is formed by joining each edge of the path $(x_2,x_1,x'_1,x'_2)$ to each of $y_2$ and $y'_2$.
Similarly, the second disk is formed by joining each edge of the path $(y_2,y_1,y'_1,y'_2)$ to each of $x_2$ and $x'_2$.
The two disks have the common boundary $(x_2,y_2,x'_2,y'_2,x_2)$ and have no other intersection.
Hence, their union is homeomorphic to $S^2$.
Therefore, the clique complex of $G$ is homotopy equivalent to $S^2$ and hence simply connected.
Since $G$ has no $4$-clique, $\mathsf{CW}_{\triangle\square}(G)$ is obtained from the clique complex of $G$ by attaching $2$-cells along chordless $4$-cycles.
Thus, $\mathsf{CW}_{\triangle\square}(G)$ is also simply connected.
\end{proof}
This shows that the hypothesis of local weak modularity in \zcref{thm:triangle-square} cannot be replaced by meshedness.

\section{Distance-regular graphs}\label{sec:drg}
\subsection{Background and overview}
A graph $G=(V,E)$ is \emph{distance-transitive} if, for any two pairs of vertices $(u,w),(x,y)$ with $d(u,w)=d(x,y)$, there is an automorphism $\sigma:V\to V$ such that $\sigma(u)=x$ and $\sigma(w)=y$.
Cycles, the icosahedral graph, Hamming graphs, Johnson graphs, Grassmann graphs, halved cubes, and bilinear forms graphs are distance-transitive.
So are the classical dual polar graphs defined by nondegenerate alternating, quadratic, or Hermitian forms over finite fields \cite[Sections~2.10 and 3.1]{vDKT16}.
In a distance-transitive graph, the number of vertices at specified distances from two given vertices depends only on the distance between the two vertices.
This counting property gives the following class.

\begin{defi}[{\cite[Section~12.1]{BH2012}}]
A connected graph $G$ of diameter $D$ is \emph{distance-regular} if, whenever $d(u,w)=i$, the numbers
\begin{align}
c_i&=\abs{\Gamma(w)\cap\Gamma_{i-1}(u)} &&(1\le i\le D),\\
a_i&=\abs{\Gamma(w)\cap\Gamma_i(u)} &&(0\le i\le D),\\
b_i&=\abs{\Gamma(w)\cap\Gamma_{i+1}(u)} &&(0\le i<D)
\end{align}
depend only on $i$.
We set $c_0=b_D=0$ and $k=b_0$.
The numbers $a_i,b_i,c_i$ are the \emph{intersection numbers}, and
$\{b_0,\ldots,b_{D-1};c_1,\ldots,c_D\}$ is the \emph{intersection array} of $G$.
\end{defi}

An automorphism preserves distances and neighborhoods, so a distance-transitive graph is distance-regular. By definition, a distance-regular graph is $k$-regular, with
\begin{equation}
a_0=0,\qquad a_i+b_i+c_i=k\quad(0\le i\le D).
\end{equation}
If $D\ge1$, then $c_1=1$.

Two common approaches to classifying distance-regular graphs are to restrict the eigenvalues of the adjacency matrix and to restrict the intersection array.
Examples of the first approach include the classification of distance-regular graphs of diameter $2$ with smallest eigenvalue $-2$ \cite{S1968}, and results on graphs of diameter $3$  whose smallest eigenvalue lies in $[-3,-2)$ \cite{BGK19, KYLCM25}.
Moore graphs, conference graphs, and Taylor graphs are families defined by special relations among their intersection numbers.

We study diagonality and weak modularity in several families of distance-regular graphs.
We first give a sufficient criterion and necessary inequalities for weak modularity in terms of intersection numbers.
For strongly regular graphs with smallest eigenvalue $-2$, we determine which are diagonal and which are weakly modular.
We also classify the weakly modular distance-regular graphs of diameter $3$ with smallest eigenvalue in $[-3,-2)$ and study conference graphs and Taylor graphs.
Finally, we examine what the intersection array determines about magnitude homology.
In particular, we give examples of strongly regular graphs with the same parameters but different $\MH_{2,3}$.

\subsection{Intersection arrays and weak modularity}
\begin{prop}\label{prop:drg-wmg}
Let $G$ be a distance-regular graph of diameter $D\ge2$ with intersection numbers $a_i,b_i,c_i$.
Suppose that
\begin{itemize}
    \labelitem{i}{cond:drg-tc} $c_i>\min\{a_i,b_1\}$ for $2\le i\le D$;
    \labelitem{ii}{cond:drg-qc} $2c_i>c_{i+1}$ for $2\le i\le D-1$.
\end{itemize}
Then $G$ is weakly modular.
\end{prop}

\begin{proof}
We argue by induction on $r$, simultaneously for all vertices $x\in V$. We assume \ref{item:triangle-condition} and \ref{item:quadrangle-condition} for $d(x,u)=d(x,w)=r-1$ and prove them for $d(x,u)=d(x,w)=r$. When $r=1$, the common neighbor is $x$ itself, so assume $r\ge2$.

We first verify \ref{item:triangle-condition}. Let $u,w\in\Gamma_r(x)$ with $uw\in E$. If $c_r>b_1$, then $w\notin\Gamma_{r-1}(x)$, and
$\abs{\Gamma(u)\setminus(\Gamma(w)\cup\{w\})}=b_1$ and $\abs{\Gamma(u)\cap \Gamma_{r-1}(x)}=c_r$ give
\begin{equation}
\abs{\Gamma(u)\cap\Gamma(w)\cap\Gamma_{r-1}(x)}\ge c_r-b_1>0.
\end{equation}

If $c_r\le b_1$, then $c_r>a_r$ by \ref{cond:drg-tc}.
Set
\begin{equation}\label{eq:wmg-parameter}
X=\Gamma(x)\cap\Gamma_{r-1}(u),\qquad
Y=\Gamma(x)\cap\Gamma_{r-1}(w).
\end{equation}
By definition of the intersection numbers, $\abs X=\abs Y=c_r$. Since $uw\in E$, the triangle inequality gives
$X\setminus Y\subseteq\Gamma(x)\cap\Gamma_r(w)$, and hence $\abs{X\cap Y}\ge c_r-a_r>0$.
Choose $x'\in X\cap Y\subseteq \Gamma(x)$. The induction hypothesis at $x'$ gives $y\in\Gamma(u)\cap\Gamma(w)$ with $d(x',y)=r-2$. The triangle inequality then gives
\begin{equation}\label{eq:drg-induction-distance}
r-1=d(x,u)-1\le d(x,y)\le d(x,x')+d(x',y)=r-1,
\end{equation}
and hence $d(x,y)=r-1$.

We next verify \ref{item:quadrangle-condition}. Let $u,w\in\Gamma_r(x)$ with $d(u,w)=2$, and suppose that there is a vertex $z\in\Gamma(u)\cap\Gamma(w)\cap\Gamma_{r+1}(x)$. In this case $r\le D-1$. Define $X,Y$ by \eqref{eq:wmg-parameter}. The triangle inequality gives $X\cup Y\subseteq \Gamma(x)\cap\Gamma_r(z)$.
Since $d(x,z)=r+1$, we have $\abs{\Gamma(x)\cap\Gamma_r(z)} = c_{r+1}$. As $\abs{X}=\abs{Y}=c_r$, condition \ref{cond:drg-qc} gives
\begin{equation}
\abs{X\cap Y} =\abs{X}+\abs{Y}-\abs{X\cup Y} \ge2c_r-c_{r+1}>0.
\end{equation}
Choose $x'\in X\cap Y$. Then $d(x',u)=d(x',w)=r-1$ and $d(x',z)=r$. The induction hypothesis for \ref{item:quadrangle-condition} at $x'$ gives $y\in\Gamma(u)\cap\Gamma(w)$ with $d(x',y)=r-2$. The inequalities in \eqref{eq:drg-induction-distance} give $d(x,y)=r-1$.
\end{proof}

\begin{prop}\label{prop:wmg-drg-necessary}
Let $G$ be a weakly modular distance-regular graph of diameter $D\ge2$.
Then
\begin{equation}\label{eq:wmg-drg-ci-bound}
c_i\le1+(c_2-1)c_{i-1}\qquad(2\le i\le D).
\end{equation}
Moreover, $c_2\ge2$.
In particular, if $D\ge3$ and $c_2=2$, then $c_3\le3$.
\end{prop}

\begin{proof}
First, we prove \eqref{eq:wmg-drg-ci-bound}.
Let $2\le i\le D$ and take two vertices with $d(x,z)=i$.
Set $X=\Gamma(z)\cap\Gamma_{i-1}(x)$. Since $\abs X=c_i>0$, we can choose $u\in X$.
Suppose that $X\setminus u\neq\emptyset$.
For $w\in X\setminus u$, if $uw\in E$, apply \ref{item:triangle-condition}.
If $uw\notin E$, then $d(u,w)=2$, and we apply \ref{item:quadrangle-condition} using the common neighbor $z$.
In either case, $\Gamma(u)\cap\Gamma(w)\cap\Gamma_{i-2}(x)\ne\emptyset$.
Thus for each $w\in X\subseteq \Gamma(z)$ there is a vertex $y\in \Gamma(u)\cap\Gamma_{i-2}(x)$ with $w\in \Gamma(y)$. Hence
\begin{equation}\label{eq:wmg-drg-ci-bound-set}
X\setminus u\subseteq
\bigcup_{y\in\Gamma(u)\cap\Gamma_{i-2}(x)}
\bigl((\Gamma(y)\cap\Gamma(z))\setminus u\bigr).
\end{equation}
Each vertex $y\in \Gamma(u)\cap\Gamma_{i-2}(x)$ satisfies $d(y,z)=2$ by the triangle inequality.
Since $u\in\Gamma(y)\cap\Gamma(z)$, the size of each set in the union is given by $\abs{(\Gamma(y)\cap\Gamma(z))\setminus u}=c_2-1$.
Since $\abs{\Gamma(u)\cap\Gamma_{i-2}(x)}= c_{i-1}$, comparing the sizes of the two sides of \eqref{eq:wmg-drg-ci-bound-set} gives $c_i-1\le (c_2-1)c_{i-1}$, which is \eqref{eq:wmg-drg-ci-bound}.
If $X\setminus u=\emptyset$, then \eqref{eq:wmg-drg-ci-bound-set} still holds and gives \eqref{eq:wmg-drg-ci-bound} in the same way.

Next, we show that $c_2\ge2$. Suppose that $c_2=1$. By \eqref{eq:wmg-drg-ci-bound} and $c_i>0$ we have $c_i=1$ for all $1\le i\le D$. Choose $z\in\Gamma_D(x)$ and write $\Gamma(z)\cap\Gamma_{D-1}(x)=u$. Every $w\in\Gamma(z)\setminus u$ belongs to $\Gamma_D(x)$. The common neighbor in $\Gamma_{D-1}(x)$ supplied by \ref{item:triangle-condition} for the edge $zw$ must be $u$, so $uw\in E$. It follows that
\begin{equation}
\Gamma(z)\setminus u\subseteq\Gamma(u)\setminus z.
\end{equation}
Both sides have size $k-1$, so equality holds.
However, $u$ has a neighbor in $\Gamma_{D-2}(x)$. The neighbor belongs to $\Gamma(u)\setminus z$ and lies outside $\Gamma(z)\setminus u$.
The two sets therefore cannot be equal.

The last assertion follows from \eqref{eq:wmg-drg-ci-bound} with $i=3$.
\end{proof}

\subsection{Strongly regular graphs}
We study diagonality and weak modularity for strongly regular graphs, which are distance-regular graphs of diameter $2$.

\begin{defi}[{\cite[Section~9.1]{BH2012}}]
A connected graph $G$ with $v$ vertices is \emph{strongly regular} with parameters $(v,k,\lambda,\mu)$ if it is neither complete nor edgeless, every vertex has $k$ neighbors, and
\begin{equation}
\abs{\Gamma(u)\cap\Gamma(w)}=
\begin{cases}
\lambda &\text{if }uw\in E,\\
\mu &\text{if }u\ne w\text{ and }uw\notin E.
\end{cases}
\end{equation}
\end{defi}
The Petersen graph, $Q_5^+$, and the Schl\"afli graph are strongly regular with parameters $(10,3,0,1)$, $(16,10,6,6)$, and $(27,16,10,8)$, respectively \cite[Section~1.1.10]{BvM22}.

Connected strongly regular graphs are precisely the distance-regular graphs of diameter $2$. In terms of the parameters of strongly regular graphs, the intersection array is
\begin{equation}\label{eq:srg-intersection-array}
\{k,k-\lambda-1;1,\mu\}.
\end{equation}
Indeed, a connected noncomplete strongly regular graph has a pair of vertices at distance $2$, so $\mu>0$. Thus every pair of nonadjacent vertices has a common neighbor, and the diameter is $2$. The definitions give $a_1=\lambda$, $c_2=\mu$, $b_1=k-\lambda-1$, and $a_2=k-\mu$. Conversely, a distance-regular graph of diameter $2$ is strongly regular with parameters $(v,k,a_1,c_2)$.

\begin{exam}[{\cite[Section~9.2]{BH2012}}]\label{exam:rook-shrikhande}
The \emph{$4\times4$ rook graph} $\mathrm{Ro}$ and the \emph{Shrikhande graph} $\mathrm{Sh}$ have vertex set $(\mathbb Z/4\mathbb Z)^2$.
Two vertices are adjacent in $\mathrm{Ro}$ if their difference lies in
\begin{equation}
\{(0,1),(0,2),(0,3),(1,0),(2,0),(3,0)\},
\end{equation}
and in $\mathrm{Sh}$ if their difference lies in
\begin{equation}
\{\pm(0,1),\pm(1,0),\pm(1,1)\}.
\end{equation}
Both graphs are strongly regular with parameters $(16,6,2,2)$.

Two vertices in $\mathrm{Ro}$ are adjacent if they differ in exactly one coordinate, so $\mathrm{Ro}\cong H(2,4)$ and the graph is diagonal. Gu's computations \cite[Proposition~A.2 and Tables~1--2]{Gu18} give
\begin{equation}
\rank\MH_{3,4}(\mathrm{Ro})=0,\qquad \rank\MH_{3,4}(\mathrm{Sh})=144,
\end{equation}
so $\mathrm{Sh}$ is not diagonal. By \eqref{eq:valency-constant-diam2}, the two graphs have the same magnitude:
\begin{equation}
\Mag(\mathrm{Ro};t)=\Mag(\mathrm{Sh};t)=\frac{16}{1+6t+9t^2}.
\end{equation}
Thus the intersection array determines magnitude, but in general it does not determine magnitude homology or diagonality.
\end{exam}

For diameter $2$, only \ref{cond:drg-tc} is needed in \zcref{prop:drg-wmg}. Together with \zcref{prop:wmg-paw-med}, this gives the following sufficient condition.

\begin{prop}\label{thm:pawful-item1}
A strongly regular graph with parameters $(v,k,\lambda,\mu)$ is pawful if
\begin{equation}
k<\max\{\lambda+\mu+1,2\mu\}.
\end{equation}
\end{prop}

\begin{proof}
We have $D=2$, $a_2=k-\mu$, $b_1=k-\lambda-1$, and $c_2=\mu$.
Thus $c_2>\min\{a_2,b_1\}$ is equivalent to $k<\max\{\lambda+\mu+1,2\mu\}$. The graph is weakly modular by \zcref{prop:drg-wmg} and pawful by \zcref{prop:wmg-paw-med}.
\end{proof}

We next consider the boundary case $k=\lambda+\mu+1=2\mu$.
In this case, \eqref{eq:drg-sphere-sizes} gives
\begin{equation}
k_2=\frac{k(k-\lambda-1)}{\mu}=k,\qquad v=1+k+k_2=2k+1.
\end{equation}
Graphs with these parameters are called conference graphs.
\begin{defi}[{\cite[Section~9.1.6]{BH2012}}]
A \emph{conference graph} is a strongly regular graph with parameters
\begin{equation}
\left(v,\frac{v-1}{2},\frac{v-5}{4},\frac{v-1}{4}\right),\qquad v\ge5.
\end{equation}
\end{defi}

Here is a standard example of a conference graph.
\begin{defi}[{\cite[Section~9.1.2]{BH2012}}]
For a prime power $q\equiv1\pmod4$, the \emph{Paley graph} has vertex set $\mathbb{F}_q$, where two distinct vertices $x$ and $y$ are adjacent if $x-y$ is a nonzero square in $\mathbb{F}_q$.
\end{defi}

When $q\equiv1\pmod4$, $-1$ is a square, so adjacency in the Paley graph is symmetric. The Paley graph has parameters
\begin{equation}
\left(q,\frac{q-1}{2},\frac{q-5}{4},\frac{q-1}{4}\right)
\end{equation}
of strongly regular graphs by \cite[Proposition~9.1.1]{BH2012} and is therefore a conference graph.

\begin{prop}\label{thm:pawful-item2}
Let $G$ be a strongly regular graph with parameters $(v,k,\lambda,\mu)$ satisfying $\lambda>0$ and $k=\lambda+\mu+1=2\mu$, i.e., a conference graph with $v>5$.
Then $G$ is pawful.
\end{prop}
The proof is given in \zcref{proof-conference}.
The unique conference graph on five vertices is the cycle $C_5$, which is also a Moore graph of diameter $2$.
We next consider Moore graphs.

\begin{defi}[{\cite[Section~12.9]{BH2012}}]
A $k$-regular graph of diameter $D\ge1$ is a \emph{Moore graph} if its number of vertices is
\begin{equation}
v=1+k\sum_{i=0}^{D-1}(k-1)^i.
\end{equation}
\end{defi}
A Moore graph is distance-regular with intersection array $\{k,k-1,\ldots,k-1;1,1,\ldots,1\}$ \cite[Section~6.1]{BCN1989}.
Suppose that $D\ge2$. If $k=2$, then the graph is the odd cycle $C_{2D+1}$.
If $k\ge3$, then $D=2$ by \cite[Theorem~12.9.1]{BH2012}, and the intersection
array gives the parameters $(k^2+1,k,0,1)$ of strongly regular graphs.
In both cases, the computation of Asao and Wakatsuki gives the following result.
We use the following consequence.
\begin{prop}[{\cite[Theorem~7.1]{AW24}}]\label{prop:moore-nondiag}
A Moore graph $G$ of diameter $D\ge2$ satisfies
\begin{equation}
\MH_{2,D+1}(G)\cong\mathbb Z^{\,vk(k-1)^D},
\end{equation}
and hence is not diagonal.
\end{prop}

\begin{defi}[{\cite[Section~10.3]{BvM22}}]\label{def:petersen}
The \emph{Petersen graph} has vertex set $\binom{[5]}2$, with two subsets adjacent if they are disjoint.
\end{defi}
The Petersen graph is a Moore graph of diameter $2$.
Its parameter is $(10,3,0,1)$.
In particular, the Petersen graph is not diagonal \cite[Examples~7.6]{AW24}.

\begin{defi}[{\cite[Section~10.11]{BvM22}}]\label{def:chang}
Write $ij=\{i,j\}$ and take the following subsets of the vertex set of $J(8,2)$:
\begin{align}
S_1&=\{12,34,56,78\},\\
S_2&=\{12,23,34,45,56,67,78,18\},\\
S_3&=\{12,23,13,45,56,67,78,48\}.
\end{align}
For each $j\in\{1,2,3\}$, switch adjacency between $S_j$ and its complement, leaving all other adjacencies unchanged. We denote the resulting graph by $\mathrm{Ch}_j$. These are the three \emph{Chang graphs}.
\end{defi}

Each Chang graph is strongly regular with parameters $(28,12,6,4)$ \cite[Section~10.11]{BvM22}.

By \zcref{prop:wmg-paw-med}, weak modularity and pawfulness are equivalent for strongly regular graphs.
Write $\theta_{\min}$ for the smallest eigenvalue of the adjacency matrix of a graph.
The distinct eigenvalues of a distance-regular graph are the eigenvalues of a tridiagonal matrix constructed from the intersection array. In particular, the array determines $\theta_{\min}$ \cite[Section~1.2]{BvM22}.
Among connected, noncomplete strongly regular graphs, $C_5$ has $\theta_{\min}=-(1+\sqrt{5})/2$, whereas every other graph satisfies $\theta_{\min}\le-2$ \cite[Section~1.1.10]{BvM22}. Thus, $\theta_{\min}=-2$ is the largest possible value for this smallest eigenvalue without $C_5$.
Seidel classified all strongly regular graphs in this case. We use the classification to decide which of them are diagonal and which are weakly modular.

\begin{prop}\label{exam:chang-tryy}
All Chang graphs are reducive, and none of them are weakly modular.
\end{prop}
The proof is given in \zcref{proof-of-chang}.

\begin{thm}\label{thm:srg-minus2-wmg}
Let $G$ be a strongly regular graph with smallest eigenvalue $\theta_{\min}=-2$.
The following are equivalent: $G$ is reducive, $G$ is diagonal, and $G$ is isomorphic to neither the Petersen graph nor the Shrikhande graph.
Moreover, $G$ is not weakly modular if and only if $G$ is isomorphic to the Petersen graph, the Shrikhande graph or one of the Chang graphs.
\end{thm}

\begin{proof}
The strongly regular graphs with smallest eigenvalue $-2$ are, up to isomorphism:
\begin{align}
&K_{m\times2}\quad(m\ge2),\qquad H(2,m)\quad(m\ge 2),\qquad J(m,2)\quad(m\ge 4),\\
&\text{Petersen graph},\qquad Q_5^+,\qquad \text{Schl\"afli graph},\qquad \mathrm{Sh},\qquad \mathrm{Ch}_j\quad(j=1,2,3).
\end{align}
Note that $H(2,2)\cong K_{2\times2}$ and $J(4,2)\cong K_{3\times2}$.
This is Seidel's classification, as stated in \cite[Theorem~1.1.1]{BvM22}.
Among these graphs, $Q_5^+$ and the Schl\"afli graph are weakly modular \cite{CCHO20}.
By \zcref{prop:wmg-paw-med, prop:ham-wmg, prop:johnson-wmg}, $K_{m\times 2}$, $H(2,m)$, and $J(m,2)$ are also weakly modular.
By \zcref{prop:wmg-reducible}, these graphs are reducive and diagonal.

The Petersen graph is not diagonal \cite[Examples 7.6]{AW24}.
By the contrapositives of \zcref{thm:reducible, prop:wmg-reducible}, the Petersen graph is neither reducive nor weakly modular.
Similarly, the Shrikhande graph is known not to be diagonal \cite{Gu18}.
Thus $\mathrm{Sh}$ is neither reducive nor weakly modular.

Finally, by \zcref{exam:chang-tryy} the Chang graphs are reducive but not weakly modular.

The diagonality, reducivity, and weak modularity of strongly regular graphs with smallest eigenvalue $-2$ are summarized in \zcref{tab:srg-min-eigenvalue-minus-two}.

\begingroup
\providecommand{\OKref}[1]{
  \shortstack{{\small #1}}}
\providecommand{\NGref}[1]{
  \shortstack{$\times$\\[-2pt]{\scriptsize #1}}}

\begin{table}[htbp]
  \centering
  \small
  \setlength{\tabcolsep}{4pt}
  \renewcommand{\arraystretch}{1.3}
  \resizebox{\linewidth}{!}{
  \begin{tabular}{lcccc}
    \hline
    Graph & Parameters $(v,k,\lambda,\mu)$
    & diagonal & reducive & weakly modular \\
    \hline
    $K_{m\times2}$ ($m\ge2$)
    & $(2m,2m-2,2m-4,2m-2)$
    & \OKref{\cite[Example~7.6]{HW17}}
    & \OKref{\zcref{prop:wmg-reducible}}
    & \OKref{\zcref{prop:wmg-paw-med}} \\
    $H(2,m)$ ($m\ge2$)
    & $(m^2,2(m-1),m-2,2)$
    & \OKref{\cite[Example~2.5, Prop.~7.3]{HW17}}
    & \OKref{\zcref{prop:wmg-reducible}}
    & \OKref{\zcref{prop:ham-wmg}} \\
    $J(m,2)$ ($m\ge4$)
    & $\left(\binom{m}{2},2(m-2),m-2,4\right)$
    & \OKref{\zcref{prop:wmg-reducible}}
    & \OKref{\zcref{prop:wmg-reducible}}
    & \OKref{\zcref{prop:johnson-wmg}} \\
    \hline
    Petersen graph
    & $(10,3,0,1)$
    & \NGref{\cite[Example~7.6]{AW24}}
    & \NGref{\zcref{thm:reducible}}
    & \NGref{\zcref{prop:wmg-reducible}} \\
    $Q_5^+$
    & $(16,10,6,6)$
    & \OKref{\zcref{prop:wmg-reducible}}
    & \OKref{\zcref{prop:wmg-reducible}}
    & \OKref{\cite[Proposition~4.4(1)]{CCHO20}} \\
    Schl\"afli graph
    & $(27,16,10,8)$
    & \OKref{\zcref{prop:wmg-reducible}}
    & \OKref{\zcref{prop:wmg-reducible}}
    & \OKref{\cite[Proposition~4.4(1)]{CCHO20}} \\
    $\mathrm{Sh}$
    & $(16,6,2,2)$
    & \NGref{\cite[Prop.~A.2, Table~2]{Gu18}}
    & \NGref{\zcref{thm:reducible}}
    & \NGref{\zcref{prop:wmg-reducible}} \\
    $\mathrm{Ch}_j\ (j=1,2,3)$
    & $(28,12,6,4)$
    & \OKref{\zcref{thm:reducible}}
    & \OKref{\zcref{exam:chang-tryy}}
    & \NGref{\zcref{exam:chang-tryy}} \\
    \hline
  \end{tabular}
  }
  \caption{Strongly regular graphs with smallest eigenvalue $-2$.}
  \label{tab:srg-min-eigenvalue-minus-two}
  \par\smallskip
\end{table}
\endgroup

\end{proof}

\begin{rema}
    The weak modularity of $Q_5^+$ and the Schl\"afli graph can also be checked from their parameters using \zcref{prop:drg-wmg}.
\end{rema}

\begin{rema}\label{rem:srg-parameters-not-enough}
The parameters of a strongly regular graph do not determine whether it is weakly modular. Indeed, $H(2,4)$ and the Shrikhande graph both have parameters $(16,6,2,2)$, but only the former is weakly modular. Likewise, $J(8,2)$ and the three Chang graphs all have parameters $(28,12,6,4)$, but only $J(8,2)$ is weakly modular. The comparison in \zcref{exam:rook-shrikhande} shows that the parameters do not determine diagonality either.
\end{rema}

\subsection{Taylor graphs and distance-regular graphs of diameter $3$}

\begin{defi}[{\cite[Section~1.2.7]{BvM22}}]\label{def:taylor}
A \emph{Taylor graph} is a distance-regular graph with intersection array $\{k,m,1;1,m,k\}$, where $k\ge2$ and $m\ge1$.
\end{defi}

The icosahedral graph is a Taylor graph with intersection array $\{5,2,1;1,2,5\}$ \cite[Theorem~3.1]{BGK19}.
For the Gosset graph, the intersection array is $\{27,10,1;1,10,27\}$ \cite[Section~10.10]{BvM22}.
The Higman--Sims group and the Conway group $\mathrm{Co}_3$ are sporadic finite simple groups that give rise to Taylor graphs with $352$ and $552$ vertices, respectively \cite[Proposition~2.7(vii), (viii)]{IM24} and \cite[Section~8.10.4]{BvM22}.

\begin{lemma}\label{prop:taylor-wmg}
Let $G$ be a Taylor graph with intersection array $\{k,m,1;1,m,k\}$.
Then $G$ always satisfies \ref{item:triangle-condition}, and \ref{item:quadrangle-condition} is equivalent to $3m>k+1$.
In particular, $G$ is weakly modular if and only if $3m>k+1$.
Moreover, if $3m>k+1$, then $G$ is diagonal.
\end{lemma}

\begin{proof}
The intersection array gives, for every $x\in V$,
$\abs{\Gamma(x)}=\abs{\Gamma_2(x)}=k$ and $\abs{\Gamma_3(x)}=1$. Write $\Gamma_3(x)=\{x^*\}$. Since $c_3=k$, we have $\Gamma(x^*)=\Gamma_2(x)$, and by symmetry of the distance we have $(x^*)^*=x$. We also have
\begin{equation}\label{eq:taylor-antipodes}
d(u,w^*)=3-d(u,w)\qquad(u,w\in V).
\end{equation}
For $d(u,w)\in\{0,3\}$, the formula follows from the definition of $w^*$. If $d(u,w)=2$, then $u\in\Gamma_2(w)=\Gamma(w^*)$. If $d(u,w)=1$, then $u$ cannot be at distance $0,1,3$ from $w^*$, so the distance is $2$.

Fix $x\in V$ and consider $\Gamma_2(x)=\Gamma(x^*)$.
The induced subgraph on $\Gamma_2(x)$ is regular of degree $a_2=a_1=k-m-1$. For distinct $u,w\in \Gamma_2(x)$, set
\begin{equation}
s_0=\abs{\Gamma(u)\cap\Gamma(w)\cap \Gamma_2(x)},\qquad
s_1=\abs{\Gamma(u)\cap\Gamma(w)\cap\Gamma(x)}.
\end{equation}
By \eqref{eq:taylor-antipodes}, the map $(-)^*:V\to V$ given by $w\mapsto w^*$ is an automorphism of $G$ and in particular gives a bijection from $\Gamma_2(x)$ to $\Gamma(x)$.
Since $u^*,w^*\notin \Gamma_2(x)$, the quantity $s_1=\abs{\Gamma(u^*)\cap\Gamma(w^*)\cap\Gamma(x^*)}=\abs{\Gamma_2(u)\cap\Gamma_2(w)\cap\Gamma_2(x)}$ is the number of vertices in $\Gamma_2(x)\setminus\{u,w\}$ adjacent to neither of the vertices $u,w$.

If $uw\in E$, counting neighbors in $\Gamma_2(x)\setminus\{u,w\}$ gives
\begin{align}
    s_1&=\abs{\Gamma_2(x)\setminus(\Gamma(u)\cup\Gamma(w))}\\
    &=\abs{\Gamma_2(x)}-\abs{\Gamma(u)\cap\Gamma(x^*)}-\abs{\Gamma(w)\cap\Gamma(x^*)}+\abs{\Gamma(u)\cap\Gamma(w)\cap\Gamma(x^*)}\\
    &=k-2a_1+s_0.
\end{align}
The common neighbors of $u,w$ consist of $x^*$, $s_0$ vertices in $\Gamma_2(x)$, and $s_1$ vertices in $\Gamma(x)$. Thus $a_1=\abs{\Gamma(u)\cap\Gamma(w)}=1+s_0+s_1$.
Hence $s_1=m/2$.
If $uw\notin E$, then $d(u,w)=2$, and the same count gives
\begin{align}
    s_1&=\abs{\Gamma_2(x)\setminus(\Gamma(u)\cup\Gamma(w)\cup \{u,w\})}\\
    &=\abs{\Gamma_2(x)}-2-\abs{\Gamma(u)\cap\Gamma(x^*)}-\abs{\Gamma(w)\cap\Gamma(x^*)}+\abs{\Gamma(u)\cap\Gamma(w)\cap\Gamma(x^*)}\\
    &=k-2-2a_1+s_0=k-2-2(k-m-1)+s_0=s_0+2m-k,\\
    m&=c_2=\abs{\Gamma(u)\cap\Gamma(w)}=1+s_0+s_1.
\end{align}
Combining the two cases, we obtain
\begin{equation}\label{eq:taylor-common-neighbor-count}
\abs{\Gamma(u)\cap\Gamma(w)\cap\Gamma(x)}=
\begin{cases}
\dfrac m2 &\text{if }uw\in E,\\[7pt]
\dfrac{3m-k-1}{2} &\text{if }uw\notin E,
\end{cases}
\qquad(u,w\in\Gamma_2(x),\ u\ne w).
\end{equation}

For $r=1$, both \ref{item:triangle-condition} and \ref{item:quadrangle-condition} hold with $x$ itself as the common neighbor. Since $\abs{\Gamma_3(x)}=1$, the only remaining case of \ref{item:triangle-condition} is when $u,w\in \Gamma_2(x)$ are adjacent. If such a pair exists, the common-neighbor count in \eqref{eq:taylor-common-neighbor-count} is $m/2>0$, so \ref{item:triangle-condition} always holds.

The remaining case of \ref{item:quadrangle-condition} is a nonadjacent pair $u,w\in \Gamma_2(x)$ with $x^*\in\Gamma_3(x)$ as a common neighbor. This case occurs because, for each vertex of $\Gamma_2(x)$, the number of its nonneighbors in $\Gamma_2(x)$ other than itself is $k-1-a_1=m>0$. By \eqref{eq:taylor-common-neighbor-count}, condition \ref{item:quadrangle-condition} holds if and only if $3m>k+1$. Diagonality follows from \zcref{thm:reducible, prop:wmg-reducible}.
\end{proof}

\begin{prop}[icosahedral graph]\label{exam:icosahedral-tryy}
The icosahedral graph $G$ is reducive.
\end{prop}

\begin{proof}
For each $u\in V$, let $\mathcal{R}_u$ consist of the empty set and all nonempty sets $S\subseteq\Gamma(u)$ that induce a connected subgraph. Since $\Gamma(u)$ induces $C_5$, we have $\Gamma(u)\in\mathcal{R}_u$. For a nonempty set $S\in\mathcal{R}_u$, choose any $w\in S$ if $S=\Gamma(u)$. Otherwise, choose an endpoint of the path induced by $S$ if $S$ has at least two vertices, or the only vertex of $S$ if $S$ is a singleton. Then $S\setminus w\in\mathcal{R}_u$.

Define $R$ by \eqref{eq:tryy-R}. The set $R$ contains $\{u\}\cup(\Gamma(u)\cap\Gamma(w))$, and these three vertices form a path of length $2$ in the $C_5$ induced by $\Gamma(w)$. Every subset of $\Gamma(w)$ containing this path induces a connected subgraph, so $R\in\mathcal{R}_w$ and \ref{item:fc} holds.

By \zcref{prop:taylor-wmg}, $G$ satisfies \ref{item:triangle-condition}, so \ref{item:tc-uws} holds. For \ref{item:qc-uws}, the diameter is $3$, so it remains to check the case $r=2$. Then $d(u,x)=3$, and \eqref{eq:taylor-antipodes} gives $\Gamma(u)=\Gamma_2(x)$. Since $S\setminus w\ne\emptyset$ and $S$ induces a connected subgraph, choose $z\in(S\setminus w)\cap\Gamma(w)$. We have $d(w,x)=d(z,x)=2$, so the triangle condition gives $\Gamma(w)\cap\Gamma(z)\cap\Gamma(x)\ne\emptyset$. Thus \ref{item:qc-uws} holds as well, and \zcref{thm:tryy} shows that $G$ is reducive.
\end{proof}

\begin{thm}\label{cor:taylor-wmg-classification}
Let $G$ be a Taylor graph.
The following are equivalent: $G$ is reducive, $G$ is diagonal, and $G$ is not isomorphic to $C_6$.
Moreover, $G$ is not weakly modular if and only if $G$ is isomorphic to $C_6$ or to the icosahedral graph.
\end{thm}

\begin{proof}
Let $\{k,m,1;1,m,k\}$ be the intersection array of $G$, and suppose that $G$ is not weakly modular. By \zcref{prop:taylor-wmg}, we have $3m\le k+1$.
Fix $x\in V$ and consider $\Gamma_2(x)$.
For $u\in \Gamma_2(x)$, we have $\abs{\Gamma(u)\cap \Gamma(x^*)}=a_1=k-m-1<k$ and $k-1-a_1=m>0$. Thus $u$ has at least $1$ nonneighbor $w$ in $\Gamma_2(x)$ other than itself.
For such a pair $u,w$, the common-neighbor count in \eqref{eq:taylor-common-neighbor-count} is nonnegative: $\abs{\Gamma(u)\cap\Gamma(w)\cap\Gamma(x)}\ge 0$. Together with $3m\le k+1$, this gives $3m= k+1$.
By \eqref{eq:taylor-antipodes}, the left-hand side of \eqref{eq:taylor-common-neighbor-count} is the number of vertices in $\Gamma_2(x)\setminus\{u,w\}$ adjacent to neither of the vertices $u,w$.
Thus $\Gamma_2(x)$ contains no three pairwise nonadjacent vertices.

Since $\Gamma_2(x)=\Gamma(x^*)$ and $(x^*)^*=x$, the neighborhood $\Gamma(v)$ contains no three pairwise nonadjacent vertices for any $v\in V$.
Therefore $G$ contains no induced $K_{1,3}$, and in particular no induced $K_{1,4}$. Of the Taylor graphs in the classification of Bang, Gavrilyuk, and Koolen \cite[Theorem~3.1]{BGK19}, only $C_6$ and the icosahedral graph satisfy $3m=k+1$.

Conversely, the intersection array of $C_6$ is $\{2,1,1;1,1,2\}$ and that of the icosahedral graph is $\{5,2,1;1,2,5\}$. Both graphs have parameters satisfying $3m=k+1$, so neither is weakly modular by \zcref{prop:taylor-wmg}.

Any Taylor graph not isomorphic to $C_6$ is reducive by \zcref{prop:wmg-reducible} and \zcref{exam:icosahedral-tryy},
and hence diagonal by \zcref{thm:reducible}.
The computation for even cycles \cite[Theorem~4.8]{Gu18} gives
\begin{equation}
\MH_{2,3}(C_6)\cong\mathbb Z^6,
\end{equation}
so $C_6$ is not diagonal.
This proves the result.
\end{proof}

We use \zcref{prop:wmg-drg-necessary} to classify the weakly modular distance-regular graphs of diameter $3$ with $-3\le\theta_{\min}<-2$.

\begin{prop}\label{thm:diam3-minus3-wmg}
Let $G$ be a distance-regular graph of diameter $3$ with $-3\le\theta_{\min}<-2$. The following conditions are equivalent.
\begin{itemize}\labelitem{i}{item:diam3-wmg} $G$ is weakly modular.
\labelitem{ii}{item:diam3-isom} $G$ is isomorphic to one of the following:
\begin{align}
&H(3,m)\quad(m\ge2),\qquad J(m,3)\quad(m\ge6),\\
&Q_6^+,\qquad Q_7^+,\qquad \text{Gosset graph}.
\end{align}\end{itemize}
In particular, all these graphs are diagonal.
\end{prop}

\begin{proof}
The implication \ref{item:diam3-isom}$\Rightarrow$\ref{item:diam3-wmg} follows from the weak modularity of the listed graphs, as discussed in \zcref{subsubsec:emg-ex}.

We prove \ref{item:diam3-wmg}$\Rightarrow$\ref{item:diam3-isom}. By \zcref{prop:wmg-drg-necessary}, we have $c_2\ge2$, and if $c_2=2$, then $c_3\le3$. Also, since $b_2>0$, the degree of $G$ satisfies $k=a_2+b_2+c_2\ge3$.

First suppose that $k=3$. The induced subgraph on each neighborhood is $a_1$-regular on $3$ vertices, so $a_1\in\{0,2\}$.
If $a_1=2$, regularity and connectedness give $G=K_4$, contrary to the diameter being $3$. Thus $a_1=0$, and $G$ contains no $3$-cycle.
If two vertices in $\Gamma_i(x)$ were adjacent, \ref{item:triangle-condition} would give a $3$-cycle. Thus $a_i=0$ for all $1\le i\le3$.
Since $c_2\ge2$ and $b_2>0$, the intersection array is $\{3,2,1;1,2,3\}$. Therefore
\begin{equation}
\abs{\Gamma_0(x)}=\abs{\Gamma_3(x)}=1,\qquad
\abs{\Gamma(x)}=\abs{\Gamma_2(x)}=3.
\end{equation}
Hence $G$ is a $3$-regular bipartite graph with parts $\Gamma_0(x)\cup\Gamma_2(x)$ and $\Gamma(x)\cup\Gamma_3(x)$, each of size $4$. The edges of $K_{4,4}$ that are absent from $G$ form a perfect matching, so $G$ is isomorphic to $H(3,2)$. 

From now on assume $k>3$.
In a distance-regular graph of degree $k>0$, a clique $C$ with $\abs C=1-\frac{k}{\theta_{\min}}$ is called a \emph{Delsarte clique}.
A noncomplete distance-regular graph is \emph{geometric} if there is a family $\mathcal{C}$ of Delsarte cliques such that $\abs{\{C\in\mathcal{C}\mid\{u,w\}\subseteq C\}}=1$ for every edge $uw\in E$ \cite[Section~2.8]{vDKT16}.

Suppose that $G$ is geometric.
The size $1-k/\theta_{\min}$ of a Delsarte clique is an integer, so $\theta_{\min}$ is rational.
Since $\theta_{\min}$ is an eigenvalue of an integer matrix, $\theta_{\min}$ is a root of a monic polynomial with integer coefficients. A rational root of such a polynomial is an integer, so the hypothesis gives $\theta_{\min}=-3$.
Let $\mathcal{C}$ be a family as in the definition of a geometric graph. For each vertex $x$,
\begin{equation}
\Gamma(x)=\bigsqcup_{\substack{C\in\mathcal{C}\\x\in C}}(C\setminus\{x\}),\qquad
\abs{C\setminus\{x\}}=\frac{k}{3}.
\end{equation}
Thus exactly three members of $\mathcal{C}$ contain $x$. Among any four vertices in $\Gamma(x)$, at least two lie in the same clique, so $G$ has no induced $K_{1,4}$. This verifies the hypotheses of the classification of Bang, Gavrilyuk, and Koolen \cite[Theorem~4.5]{BGK19}. The only candidates of diameter $3$ with $c_2\ge2$ and $\theta_{\min}=-3$ are $H(3,m)$ $(m\ge3)$ and $J(m,3)$ $(m\ge6)$. Every other candidate has $c_2=1$ or $\theta_{\min}=-2$.

If $G$ is not geometric, we apply the classification of Koolen, Yu, Liang, Choi, and Markowsky \cite[Theorem~1.1]{KYLCM25}. The candidates of diameter $3$ with $c_2\ge2$ are $Q_6^+$, $Q_7^+$, the Gosset graph, $\mathrm{Sh}\square K_4$, and the graphs with intersection arrays in \zcref{tab:diam3-candidates}.
An intersection array need not determine the graph uniquely.
\begin{table}[htbp]
\centering
\begin{tabular}{c|cc}
Intersection array & $c_2$ & $c_3$\\ \hline
$\{5,2,1;1,2,5\}$ & $2$ & $5$\\
$\{10,6,4;1,2,5\}$ & $2$ & $5$\\
$\{7,4,1;1,2,7\}$ & $2$ & $7$\\
$\{9,6,1;1,2,9\}$ & $2$ & $9$\\
$\{15,10,1;1,2,15\}$ & $2$ & $15$\\
$\{18,12,1;1,2,18\}$ & $2$ & $18$
\end{tabular}
\caption{Intersection arrays considered in the proof of \zcref{thm:diam3-minus3-wmg}. Each satisfies $c_2=2$ and $c_3>3$.}
\label{tab:diam3-candidates}
\end{table}
Every graph with one of the listed intersection arrays has $c_2=2$ and $c_3>3$, so no such graph is weakly modular by \zcref{prop:wmg-drg-necessary}.
By \zcref{prop:wmg-reducible, exam:rook-shrikhande}, $\mathrm{Sh}$ is not weakly modular, so $\mathrm{Sh}\square K_4$ is not weakly modular by \zcref{prop:cartesian-wmg}.
The remaining candidates are $Q_6^+$, $Q_7^+$, and the Gosset graph. This proves the equivalence; diagonality follows from \zcref{thm:reducible, prop:wmg-reducible}.
\end{proof}

\subsection{What the intersection array determines}\label{subsec:drg-mag}
In a distance-regular graph, $k_i=\abs{\Gamma_i(u)}$ does not depend on $u$ and is given by
\begin{equation}\label{eq:drg-sphere-sizes}
k_0=1,\qquad k_i=\frac{b_0\cdots b_{i-1}}{c_1\cdots c_i}\quad(1\le i\le D).
\end{equation}
Counting the edges between $\Gamma_i(u)$ and $\Gamma_{i+1}(u)$ in two ways gives
$\abs{\Gamma_i(u)}b_i=\abs{\Gamma_{i+1}(u)}c_{i+1}$. Starting with $\abs{\Gamma_0(u)}=1$, the edge count yields \eqref{eq:drg-sphere-sizes} by induction.
Applying \zcref{prop:valency-constant} gives the following formula.

\begin{prop}[{\cite[Proposition~7.10]{AW24}}]\label{prop:drg-magnitude}
Let $G$ be a distance-regular graph with $v$ vertices and intersection array $\{b_0,\ldots,b_{D-1};c_1,\ldots,c_D\}$.
Then
\begin{equation}\label{eq:drg-magnitude}
\Mag(G;t)=\frac{v}{1+\sum_{i=1}^{D}\frac{b_0\cdots b_{i-1}}{c_1\cdots c_i}t^i}.
\end{equation}
In particular, $\Mag(G;t)$ is determined by the intersection array.
\end{prop}

Indeed, $v=\sum_{i=0}^{D}k_i$ is also determined by the intersection array.
To study the summands with fixed endpoints, we use the relation between distance matrices and the adjacency algebra.
\begin{defi}[adjacency algebra]
Write $\mathbb{C}^{V\times V}$ for the algebra of matrices indexed by $V \times V$ with entries in $\mathbb{C}$.
The \emph{adjacency algebra} $\mathfrak{A}$ of $G$ is the subalgebra of $\mathbb{C}^{V\times V}$ consisting of polynomials in the adjacency matrix $A$ with coefficients in $\mathbb{C}$:
\begin{equation}
\mathfrak{A} = \mathbb{C}[A] = \{\alpha(A) \mid \alpha(\tau)\in \mathbb{C}[\tau] \}.
\end{equation}
\end{defi}

If $G$ is distance-regular with diameter at least $1$, then its distance matrices satisfy
\begin{equation}\label{eq:drg-distance-recurrence}
\begin{aligned}
AA_0&=A_1,\\
AA_i&=b_{i-1}A_{i-1}+a_iA_i+c_{i+1}A_{i+1} &&(1\le i<D),\\
AA_D&=b_{D-1}A_{D-1}+a_DA_D.
\end{aligned}
\end{equation}
These formulas express $(AA_i)_{uw}=\abs{\Gamma(u)\cap\Gamma_i(w)}$ in terms of $d(u,w)$.

\begin{prop}\label{prop:drg-adjacency-algebra}
For a connected graph $G$ of diameter $D$, the following are equivalent:
\begin{itemize}
\labelitem{i}{item:drg-adjacency-algebra-1} $G$ is distance-regular.
\labelitem{ii}{item:drg-adjacency-algebra-2} $Z_G\in\mathfrak{A}[t]$ and $\dim_{\mathbb C}\mathfrak{A}=D+1$.
\labelitem{iii}{item:drg-adjacency-algebra-3} $Z_G^{-1}\in\mathfrak{A}[\![t]\!]$ and $\dim_{\mathbb C}\mathfrak{A}=D+1$.
\end{itemize}
If these conditions hold, then $\{A_0,\ldots,A_D\}$ is a linear basis of $\mathfrak{A}$.
\end{prop}

\begin{proof}
When $D=0$, we have $G=K_1$ and all three conditions hold. Assume $D\ge1$ from now on.
First assume \ref{item:drg-adjacency-algebra-1}.
In \eqref{eq:drg-distance-recurrence} we have $c_{i+1}>0$, so induction gives $A_i\in\mathfrak{A}$.
For the reverse inclusion, the span of the distance matrices contains $I$ and is closed under left multiplication by $A$, by the recurrence. It therefore contains every polynomial in $A$.
Thus $\mathfrak{A}=\operatorname{span}_{\mathbb C}\{A_0,\ldots,A_D\}$.
The distance matrices have disjoint nonempty supports and are therefore linearly independent, so $\dim_{\mathbb C}\mathfrak{A}=D+1$. Now $Z_G=\sum_{i=0}^{D}t^iA_i$ gives \ref{item:drg-adjacency-algebra-2}.

Assume \ref{item:drg-adjacency-algebra-2}. Comparing coefficients gives $A_i\in\mathfrak{A}$, and the dimension hypothesis shows that $A_0,\ldots,A_D$ form a basis of $\mathfrak{A}$. Each entry of $AA_i$ therefore depends only on the distance between the two vertices indexing it. In $(AA_i)_{uw}=\abs{\Gamma(u)\cap\Gamma_i(w)}$, taking the values $i=d(u,w)-1,d(u,w),d(u,w)+1$ in the range $0\le i\le D$ shows that the intersection numbers exist. This proves \ref{item:drg-adjacency-algebra-1}.

The coefficients of the inverse of a formal power series with constant term $I$ are signed sums of products of the original coefficients. Thus the inverse of such a series over a unital matrix algebra has coefficients in the same algebra. Applying this to $Z_G$ and $Z_G^{-1}$ proves the equivalence of \ref{item:drg-adjacency-algebra-2} and \ref{item:drg-adjacency-algebra-3}.
\end{proof}

The condition $Z_G\in\mathfrak{A}[t]$ is equivalent to requiring all distance matrices $A_i$ to belong to $\mathfrak{A}$. This condition does not imply distance-regularity. For example, the triangular prism $K_3\square K_2$ is regular with diameter $2$, and the eigenvalues of its adjacency matrix are $3,1,0,-2$. Let $J$ be the matrix with every entry equal to $1$. Then
\begin{equation}
J=\frac15A(A-I)(A+2I),\qquad A_2=J-I-A,
\end{equation}
so $Z_G\in\mathfrak{A}[t]$ and $Z_G^{-1}\in\mathfrak{A}[\![t]\!]$. However, $\dim_{\mathbb{C}} \mathfrak{A}=4$.
Moreover, the endpoints of an edge within a $3$-cycle have one common neighbor, while the endpoints of an edge joining the two $3$-cycles have none. The graph is therefore not distance-regular, so the dimension hypothesis in \zcref{prop:drg-adjacency-algebra} cannot be dropped.

\begin{prop}\label{prop:drg-zeta-inverse}
Let $G$ be distance-regular with diameter $D\ge1$.
Define polynomials $\alpha_i(\tau)\in\mathbb Q[\tau]$ $(0\le i\le D)$ by
\begin{align}
\alpha_0(\tau)&=1,\qquad \alpha_1(\tau)=\tau,\\
c_{i+1}\alpha_{i+1}(\tau)&=(\tau-a_i)\alpha_i(\tau)-b_{i-1}\alpha_{i-1}(\tau)\quad(1\le i<D).
\end{align}
Let $\widetilde{A}$ be the tridiagonal matrix indexed by $0,\ldots,D$ with
\begin{equation}
\widetilde{A}_{ii}=a_i \quad (0\le i\le D),\qquad
\widetilde{A}_{i,i-1}=c_i \quad (1\le i\le D),\qquad
\widetilde{A}_{i,i+1}=b_i \quad (0\le i<D),
\end{equation}
and all other entries zero.
For $u,w\in V$ with $d(u,w)=r$, we have
\begin{equation}\label{eq:drg-zeta-inverse}
(Z_G^{-1})_{uw}
=\left[\left(\sum_{i=0}^{D}t^i\alpha_i(\widetilde{A})\right)^{-1}\right]_{r,0}.
\end{equation}
In particular, $(Z_G^{-1})_{uw}$ and $\abs{P_{n,\ell}(G;u,w)}$ for all $n,\ell\ge0$ are determined by $r$ and the intersection array.
\end{prop}

\begin{proof}
By \eqref{eq:drg-distance-recurrence}, we have $A_i=\alpha_i(A)$.
By \zcref{prop:drg-adjacency-algebra}, $(A_0,\ldots,A_D)$ is a basis of $\mathfrak{A}$. In this basis, left multiplication by $A$ has matrix $\widetilde{A}$, so left multiplication by $Z_G$ has matrix $\sum_{i=0}^{D}t^i\alpha_i(\widetilde{A})$.
The matrix for left multiplication by $Z_G$ has constant term $I$ and is therefore invertible as a formal power series.

By \zcref{prop:drg-adjacency-algebra}, $Z_G^{-1}$ is a linear combination of the distance matrices with formal power series coefficients.
Writing $Z_GZ_G^{-1}=A_0$ in this basis, we see that the coefficient vector of $Z_G^{-1}$ is the column of $\left(\sum_{i=0}^{D}t^i\alpha_i(\widetilde{A})\right)^{-1}$ corresponding to $A_0$.
At a pair of vertices at distance $r$, $(Z_G^{-1})_{uw}$ equals the coefficient of $A_r$. This proves \eqref{eq:drg-zeta-inverse}.

Applying the same matrix representation to $(Z_G-I)^n$ and using the count preceding \eqref{eq:entries-mobius} gives
\begin{equation}
\abs{P_{n,\ell}(G;u,w)}
=[t^\ell]\left[\left(\sum_{i=1}^{D}t^i\alpha_i(\widetilde{A})\right)^n\right]_{r,0}.
\end{equation}
Thus the chain counts are also determined by $r$ and the intersection array.
\end{proof}


Hence if $G$ is diagonal and distance-regular, \eqref{eq:rank-of-diagonal} gives
\begin{equation}
\rank\MH_{\ell,\ell}(G;u,w)=(-1)^\ell[t^\ell](Z_G^{-1})_{uw},
\end{equation}
and $\MH_{\ell,\ell}(G;u,w)$ is free abelian by \zcref{prop:freeness-diagonal-part} and \eqref{eq:chain-hom-direct-sum}. Thus $d(u,w)$ and the intersection array determine $\MH_{\ell,\ell}(G;u,w)$ up to isomorphism of abelian groups.

\begin{rema}
For a diagonal distance-regular graph $G$, the intersection array and the distance between two vertices $u,w$ determine the dimensions of the graded Tor groups between the corresponding right and left simple modules over $\mathsf{DA}_{\mathbb{F}}(G)$ for every field $\mathbb{F}$, by \cite[Theorem~3.3]{IM26}.
Thus they determine the corresponding graded Betti numbers.
\end{rema}

Equations \eqref{eq:drg-zeta-inverse} and \eqref{eq:rank-of-diagonal} determine $\rank\MH_{\ell,\ell}(G;u,w)$. The argument based on \zcref{prop:valency-constant} recovers the entire homology from the magnitude.
For a distance-regular graph, the low-degree groups with fixed endpoints in \zcref{prop:low-degree-magnitude-homology} are likewise determined by the endpoint distance and the intersection array.
By \eqref{eq:drg-sphere-sizes}, the intersection array gives the layer sizes $k_r$ and hence the number $\abs V=\sum_{r=0}^D k_r$ of vertices. There are $\abs V k_r$ ordered pairs at distance $r$. Taking the direct sum in \eqref{eq:chain-hom-direct-sum} over these pairs for all distances shows that the intersection array determines the entire magnitude homology of $G$ in length degrees at most $2$.
\zcref{sec:diagonal-homology-tables} lists diagonal ranks for several diagonal distance-regular graphs.
Without the diagonality hypothesis, the intersection array need not determine magnitude homology even in length degree $3$.

Distance-regular graphs with the same intersection array can have different ranks of $\MH_{2,3}$.
To give an example, we use Ihringer's publicly available data \cite[Section~3.12]{Ihr22}.

\begin{exam}\label{exam:96-19-2-4}
We use the file \texttt{srgs\_96\_19\_2\_4.gz} from \cite{IhrData}. The file contains $178,040$ graphs.
Write $G_1,G_2,G_3$ for the graphs on lines $1,2,3$ of the file, respectively.
All three have parameters $(96,19,2,4)$.
However, $\rank \MH_{2,3}(G_1)=0$, $\rank \MH_{2,3}(G_2)=6$, and $\rank \MH_{2,3}(G_3)=12$, so the parameters and magnitude of a strongly regular graph do not determine the isomorphism type of $\MH_{2,3}$.
Detailed computations are given in \zcref{append:96-19-2-4}.
\end{exam}
\begin{rema}
    Since the second magnitude homology is torsion-free \cite[Corollary~4.7]{IM26} and $\MH_{2,4}(G_j)=0$ for $j=1,2,3$, \cite[Proposition~6.4]{IM26} implies that the distance algebra of $G_1$ is quadratic over every field, whereas those of $G_2$ and $G_3$ are not quadratic over any field.
    Nevertheless, all three distance algebras have the same Hilbert series $96+1824t+7296t^2$.
    Thus quadraticity of a distance algebra is not determined by its Hilbert series.
\end{rema}

\bibliographystyle{alpha}
\bibliography{reference}

\appendix

\section{\zcref{lem:basis-reduction} and \ref{item:triangle-condition}}\label{sec:appendix-house}

By \zcref{lem:basis-reduction}, suitable families of cliques suffice to establish reducivity. However, the existence of such families is a strictly stronger condition than reducivity itself. In particular, it implies the triangle condition \ref{item:triangle-condition}.

\begin{prop}
For a graph $G=(V,E)$, suppose that for each $u\in V$ we are given a family $\mathcal{C}_u\subseteq 2^{\Gamma(u)}$ such that $\emptyset\in\mathcal{C}_u$.
If the collection $(\mathcal{C}_u)_{u\in V}$ satisfies \ref{item:basis-reduction-partition}, then $G$ satisfies the triangle condition in \ref{item:triangle-condition}.
\end{prop}

\begin{proof}
Let $x,u,w\in V$ satisfy $uw\in E$ and
$d(x,u)=d(x,w)=r$.
Since $u$ and $w$ are distinct, $r\geq 1$.

Since $\emptyset\in\mathcal{C}_u$ and $w\notin\emptyset$,
applying \ref{item:basis-reduction-partition} to $C=\emptyset$ gives
a set $R\subseteq\Gamma(w)$ such that
\begin{equation}
  V=T_u(w)\sqcup T_w(R).
\end{equation}
Since $d(x,u)=d(x,w)$, we have $x\notin T_u(w)$, and hence
$x\in T_w(R)$. Thus there is a vertex $y\in R$ such that
\begin{equation}
  d(y,x)=d(w,x)-1=r-1.
\end{equation}
Furthermore, $y\neq u$ because $d(y,x)=r-1$ while $d(u,x)=r$.

Since $y\in R$, we have $y\in T_w(R)$.
The partition $V=T_u(w)\sqcup T_w(R)$ gives $y\notin T_u(w)$. Since $uw,wy\in E$, we have $d(u,y)\leq 2$.
If $d(u,y)=2$, then
\begin{equation}
  d(u,y)=1+d(w,y)
\end{equation}
would give $y\in T_u(w)$, a contradiction. Therefore
$d(u,y)\leq 1$. Since $y\neq u$, we have $y\in\Gamma(u)$.

Thus,
\begin{equation}
  y\in\Gamma(u)\cap\Gamma(w)\cap\Gamma_{r-1}(x),
\end{equation}
and \ref{item:triangle-condition} follows.
\end{proof}

Thus a reducive graph that fails \ref{item:triangle-condition} cannot have the clique families required by \zcref{lem:basis-reduction}.
The House graph is one such example. \zcref{rema:house-not-wmg} explains why the House graph fails \ref{item:triangle-condition}.

\section{Proof of \zcref{thm:pawful-item2}}\label{proof-conference}
Let $x,y,z\in V$ satisfy $d(x,y)=1$ and $d(x,z)=d(y,z)=2$.
Since $\mu=\lambda+1$, we have $\abs{\Gamma(x)\cap\Gamma(y)}=\lambda>0$ and
$\abs{\Gamma(x)\cap\Gamma_2(y)}=k-\lambda-1=\mu>0$.
Thus, we can choose $u\in\Gamma(x)\cap\Gamma(y)$ and
$w\in\Gamma(x)\cap\Gamma_2(y)$.
The decompositions
\begin{align}
    \Gamma(x)
    &=\{y\}\sqcup(\Gamma(x)\cap\Gamma(y))
      \sqcup(\Gamma(x)\cap\Gamma_2(y)),\\
    \Gamma(y)
    &=\{x\}\sqcup(\Gamma(x)\cap\Gamma(y))
      \sqcup(\Gamma_2(x)\cap\Gamma(y)),
\end{align}
give
\begin{align}
    \lambda=\abs{\Gamma(x)\cap\Gamma(u)}
    &=1+\abs{\Gamma(x)\cap\Gamma(y)\cap\Gamma(u)}
      +\abs{\Gamma(x)\cap\Gamma_2(y)\cap\Gamma(u)},\\
    \lambda=\abs{\Gamma(y)\cap\Gamma(u)}
    &=1+\abs{\Gamma(x)\cap\Gamma(y)\cap\Gamma(u)}
      +\abs{\Gamma_2(x)\cap\Gamma(y)\cap\Gamma(u)}.
\end{align}
Taking the difference yields
\begin{equation}\label{eq:lem-pawful1}
    \abs{\Gamma(x)\cap\Gamma_2(y)\cap\Gamma(u)}
    =\abs{\Gamma_2(x)\cap\Gamma(y)\cap\Gamma(u)}.
\end{equation}
Similarly,
\begin{align}
    \lambda=\abs{\Gamma(x)\cap\Gamma(w)}
    &=\abs{\Gamma(x)\cap\Gamma(y)\cap\Gamma(w)}
      +\abs{\Gamma(x)\cap\Gamma_2(y)\cap\Gamma(w)},\\
    \mu=\abs{\Gamma(y)\cap\Gamma(w)}
    &=1+\abs{\Gamma(x)\cap\Gamma(y)\cap\Gamma(w)}
      +\abs{\Gamma_2(x)\cap\Gamma(y)\cap\Gamma(w)}.
\end{align}
Since $\mu=\lambda+1$, we obtain
\begin{equation}\label{eq:lem-pawful2}
    \abs{\Gamma(x)\cap\Gamma_2(y)\cap\Gamma(w)}
    =\abs{\Gamma_2(x)\cap\Gamma(y)\cap\Gamma(w)}.
\end{equation}

Observe the following decompositions:
\begin{align}
    \Gamma(x)\cap\Gamma(z)
    &=(\Gamma(x)\cap\Gamma(y)\cap\Gamma(z))
      \sqcup(\Gamma(x)\cap\Gamma_2(y)\cap\Gamma(z)),\\
    \Gamma(y)\cap\Gamma(z)
    &=(\Gamma(x)\cap\Gamma(y)\cap\Gamma(z))
      \sqcup(\Gamma_2(x)\cap\Gamma(y)\cap\Gamma(z)),\\
    \Gamma(x)\cap\Gamma_2(y)
    &=(\Gamma(x)\cap\Gamma_2(y)\cap\Gamma(z))
      \sqcup(\Gamma(x)\cap\Gamma_2(y)\cap\Gamma_2(z)),\\
    \Gamma_2(x)\cap\Gamma(y)
    &=(\Gamma_2(x)\cap\Gamma(y)\cap\Gamma(z))
      \sqcup(\Gamma_2(x)\cap\Gamma(y)\cap\Gamma_2(z)).
\end{align}
Since each set on the left-hand side has cardinality $\mu$, we have
\begin{align}
    \abs{\Gamma(x)\cap\Gamma(y)\cap\Gamma(z)}
    &=\abs{\Gamma(x)\cap\Gamma_2(y)\cap\Gamma_2(z)}
    =\abs{\Gamma_2(x)\cap\Gamma(y)\cap\Gamma_2(z)},\\
    \abs{\Gamma(x)\cap\Gamma_2(y)\cap\Gamma(z)}
    &=\abs{\Gamma_2(x)\cap\Gamma(y)\cap\Gamma(z)}
    =\mu-\abs{\Gamma(x)\cap\Gamma(y)\cap\Gamma(z)}.
\end{align}

Suppose for contradiction that
$\Gamma(x)\cap\Gamma(y)\cap\Gamma(z)=\emptyset$.
Under this assumption, $\Gamma(x)\cap\Gamma(y)\subseteq\Gamma_2(z)$ and $\Gamma(x)\cap\Gamma_2(y)\subseteq\Gamma(z)$. Hence $d(z,u)=2$ and $d(z,w)=1$.
Thus,
\begin{equation}\label{eq:lem-pawful3}
    \abs{\Gamma(z)\cap\Gamma(u)}
    =\mu=\lambda+1
    =\abs{\Gamma(z)\cap\Gamma(w)}+1,
\end{equation}
so the two common-neighbor counts have different parity.
Moreover,
\begin{align}
    k=\abs{\Gamma(z)}
    &=\abs{\Gamma(x)\cap\Gamma(z)}
      +\abs{\Gamma(y)\cap\Gamma(z)}
      +\abs{\Gamma_2(x)\cap\Gamma_2(y)\cap\Gamma(z)}\\
    &=2\mu+\abs{\Gamma_2(x)\cap\Gamma_2(y)\cap\Gamma(z)}.
\end{align}
Since $k=2\mu$, we obtain
\begin{gather}
    \Gamma_2(x)\cap\Gamma_2(y)\cap\Gamma(z)=\emptyset,\\
    \Gamma(x)\cap\Gamma_2(y)\cap\Gamma_2(z)=\emptyset,\qquad
    \Gamma_2(x)\cap\Gamma(y)\cap\Gamma_2(z)=\emptyset.
\end{gather}
Consequently,
\begin{gather}
    \Gamma(z)
    =(\Gamma(x)\cap\Gamma(z))\sqcup(\Gamma(y)\cap\Gamma(z)),\\
    \Gamma(x)\cap\Gamma_2(y)=\Gamma(x)\cap\Gamma(z),\qquad
    \Gamma_2(x)\cap\Gamma(y)=\Gamma(y)\cap\Gamma(z).
\end{gather}
Together with \eqref{eq:lem-pawful1} and \eqref{eq:lem-pawful2},
these identities yield
\begin{align}
    \abs{\Gamma(z)\cap\Gamma(u)}
    &=\abs{\Gamma(x)\cap\Gamma(z)\cap\Gamma(u)}
      +\abs{\Gamma(y)\cap\Gamma(z)\cap\Gamma(u)}\\
    &=\abs{\Gamma(x)\cap\Gamma_2(y)\cap\Gamma(u)}
      +\abs{\Gamma_2(x)\cap\Gamma(y)\cap\Gamma(u)}\\
    &=2\abs{\Gamma(x)\cap\Gamma_2(y)\cap\Gamma(u)},\\
    \abs{\Gamma(z)\cap\Gamma(w)}
    &=\abs{\Gamma(x)\cap\Gamma(z)\cap\Gamma(w)}
      +\abs{\Gamma(y)\cap\Gamma(z)\cap\Gamma(w)}\\
    &=\abs{\Gamma(x)\cap\Gamma_2(y)\cap\Gamma(w)}
      +\abs{\Gamma_2(x)\cap\Gamma(y)\cap\Gamma(w)}\\
    &=2\abs{\Gamma(x)\cap\Gamma_2(y)\cap\Gamma(w)}.
\end{align}
Both cardinalities are even, contradicting \eqref{eq:lem-pawful3}.
Therefore $x,y,z$ have a common neighbor, and $G$ is pawful.

\section{Proof of \zcref{exam:chang-tryy}}\label{proof-of-chang}

We first show that all Chang graphs are reducive.
Let $G\in\{\mathrm{Ch}_1,\mathrm{Ch}_2,\mathrm{Ch}_3\}$ and use the vertex notation of \zcref{def:chang}.
For each $x\in V$, set
\begin{equation}
\mathcal{B}_x\coloneq
\bigl\{\{u',w\}\in E\mid
u',w\in\Gamma_2(x),\
\Gamma(x)\cap\Gamma(u')\cap\Gamma(w)=\emptyset
\bigr\}.
\end{equation}

A permutation of $\{1,\ldots,8\}$ that preserves the set $S_j$ used to define $\mathrm{Ch}_j$ induces an automorphism of $\mathrm{Ch}_j$.
Viewed as edge sets on $[8]$, $S_1$, $S_2$, and $S_3$ form a perfect matching, an $8$-cycle, and the disjoint union of a triangle and a $5$-cycle, respectively. The symmetries of these graphs yield the orbit representatives listed in the table below.

For each listed vertex $x$, the adjacency rule gives the entries in the table.
The rightmost column counts edges in $\bigcup_{y\in V}\mathcal{B}_y$ that are incident with $x$.
\begin{equation}
\begin{array}{c|c|l|c}
G&x&\mathcal{B}_x&\text{Number of edges}\\ \hline
\mathrm{Ch}_1&12&
\{\{13,14\},\{15,16\},\{17,18\},
  \{23,24\},\{25,26\},\{27,28\}\}&0\\
\mathrm{Ch}_1&13&\emptyset&2\\
\mathrm{Ch}_2&12&\{\{13,14\},\{27,28\}\}&4\\
\mathrm{Ch}_2&13&\{\{12,23\},\{25,27\}\}&4\\
\mathrm{Ch}_2&14&\{\{12,35\},\{28,34\}\}&4\\
\mathrm{Ch}_2&15&\emptyset&0\\
\mathrm{Ch}_3&12&\emptyset&0\\
\mathrm{Ch}_3&14&\emptyset&0\\
\mathrm{Ch}_3&45&\{\{46,47\},\{57,58\},\{67,78\}\}&3\\
\mathrm{Ch}_3&46&\{\{45,56\},\{48,57\},\{58,67\}\}&3
\end{array}
\end{equation}
The table shows that, for each $x\in V$, no two distinct edges in $\mathcal{B}_x$ share an endpoint.
Each $u\in V$ is incident with $12$ edges, at most $4$ of which belong to $\bigcup_{x\in V}\mathcal{B}_x$.
Thus there are at least $8$ vertices $u'\in\Gamma(u)$ such that $\{u,u'\}\notin\mathcal{B}_x$ for every $x\in V$.

For each $u\in V$, let $\mathcal{R}_u$ consist of the empty set and
all nonempty sets $S\subseteq\Gamma(u)$ that induce a connected subgraph and satisfy the following condition:
\begin{equation}
\text{There is a vertex }u'\in S\text{ such that for every }x\in V
\text{ we have }\{u,u'\}\notin\mathcal{B}_x.
\end{equation}
The induced subgraph on $\Gamma(u)$ is $6$-regular on $12$ vertices.
The neighborhoods of two nonadjacent vertices in this subgraph are $6$-element subsets of the other $10$ vertices, so they intersect.
Thus the induced subgraph has diameter at most $2$.
The neighborhood $\Gamma(u)$ is connected and contains a vertex $u'$ such that $\{u,u'\}\notin\mathcal{B}_x$ for every $x\in V$. Hence $\Gamma(u)\in\mathcal{R}_u$, so \ref{item:canonical-condition} holds.

Take a nonempty set $S\in\mathcal{R}_u$ and fix $u'\in S$ as in the definition of the family.
If $S=\{u'\}$, set $w=u'$.
Otherwise, take a spanning tree of the subgraph induced by $S$ and choose a leaf $w$ distinct from $u'$.
In the latter case $S\setminus w$ induces a connected subgraph and contains $u'$. Thus in either case,
$S\setminus w\in\mathcal{R}_u$.

The set $R$ in \eqref{eq:tryy-R} contains the $7$-element set $\{u\}\cup(\Gamma(u)\cap\Gamma(w))$.
Since the induced subgraph on $\Gamma(w)$ has diameter at most $2$, every vertex in $R\setminus\{u\}$ is either adjacent to $u$ or joined to $u$ via a vertex in $\Gamma(u)\cap\Gamma(w)$.
Thus $R$ induces a connected subgraph.
Also, at most $4$ edges in $\bigcup_{x\in V}\mathcal{B}_x$ are incident with $w$.
Hence $R$ contains a vertex $y$ such that $\{w,y\}\notin\mathcal{B}_x$ for every $x\in V$.
Thus $R\in\mathcal{R}_w$, and \ref{item:fc} holds.

For \ref{item:tc-uws}, we only need to check $r=2$.
If the common neighbor required by the triangle condition does not exist, then $\{u,w\}\in\mathcal{B}_x$.
By the choice of $u'$, we have $w\ne u'$, so $S$ contains at least two vertices.
Let $z\in(S\setminus w)\cap\Gamma(w)$ be a neighbor of $w$ in the chosen spanning tree.
Since $\Gamma(u)\cap\Gamma(w)\cap\Gamma(x)=\emptyset$, we have $z\in\Gamma_2(x)$.
No two distinct edges in $\mathcal{B}_x$ share an endpoint, so $\{w,z\}\notin\mathcal{B}_x$ and hence
\begin{equation}
\Gamma(w)\cap\Gamma(z)\cap\Gamma(x)\ne\emptyset.
\end{equation}
Since $z\in S\setminus w$, this proves \ref{item:tc-uws}.
Since $G$ has diameter $2$, the hypotheses of \ref{item:qc-uws} can hold only when $r=1$.
When $r=1$, take $x$ as the common neighbor.
By \zcref{thm:tryy}, this proves that $G$ is reducive.

Finally, we show that none of the Chang graphs is weakly modular.

For each Chang graph, we choose three vertices satisfying $u,w\in\Gamma_2(x)$ and $uw\in E$ but with
\begin{equation}
\Gamma(x)\cap\Gamma(u)\cap\Gamma(w)=\emptyset
\end{equation}
which shows that \ref{item:triangle-condition} fails.

In the notation of \zcref{def:chang}, take the three vertices in each row of the table.
The adjacency rule gives $u,w\in\Gamma_2(x)$ and $uw\in E$, and the set in the last column is disjoint from $\Gamma(x)$.
\begin{center}
\begin{tabular}{c|ccc|l}
$G$ & $x$ & $u$ & $w$ & $\Gamma(u)\cap\Gamma(w)$\\ \hline
$\mathrm{Ch}_1$ & $12$ & $13$ & $14$ & $\{15,16,17,18,56,78\}$\\
$\mathrm{Ch}_2$ & $12$ & $13$ & $14$ & $\{15,16,17,56,67,78\}$\\
$\mathrm{Ch}_3$ & $45$ & $46$ & $47$ & $\{12,13,23,14,24,34\}$
\end{tabular}
\end{center}
Thus none of the Chang graphs is weakly modular.

\section{Details of \zcref{exam:96-19-2-4}}\label{append:96-19-2-4}
The dataset \texttt{srgs\_96\_19\_2\_4.gz} from \cite{IhrData}
contains $178{,}040$ graphs.
Let $G_1, G_2, G_3$ be the first three graphs in this dataset, in order; their \texttt{graph6} representations are provided at the end of this section.

By \zcref{prop:drg-magnitude}, we have
\begin{equation}
\Mag(G_j;t)=\frac{96}{1+19t+76t^2}\qquad(j=1,2,3).
\end{equation}
The ranks of $\MH_{n,\ell}(G_j)$ for $0\le n,\ell\le5$
are given in the following tables.
For $n=\ell$, we computed the ranks of $\MH_{\ell,\ell}(G_j)$ as the dimensions of the kernels of the corresponding boundary maps over $\mathbb Q$.

By applying the formula of Ivanov and Mukoseev \cite[Theorem~3.5]{IM26}
to graphs of diameter $2$, we obtain
\begin{equation}
\MH_{1,2}(G_j)=\MH_{2,4}(G_j)=\MH_{3,5}(G_j)=0.
\end{equation}
It follows from the definition of the chain complex and the above
equalities that, for $0\le n,\ell\le5$, the only off-diagonal groups
that can be nonzero are $\MH_{2,3},\MH_{3,4},\MH_{4,5}$.
By the formula for the Euler characteristic in \zcref{subsec:maghom},
we have
\begin{equation}
\rank\MH_{\ell-1,\ell}(G_j)
=\rank\MH_{\ell,\ell}(G_j)
-(-1)^\ell[t^\ell]\Mag(G_j;t)
\qquad(\ell=3,4,5).
\end{equation}

\begin{equation}
G_1:\qquad
\begin{array}{c|rrrrrr}
\ell\backslash n & 0 & 1 & 2 & 3 & 4 & 5\\ \hline
0 & 96 & 0 & 0 & 0 & 0 & 0\\
1 & 0 & 1824 & 0 & 0 & 0 & 0\\
2 & 0 & 0 & 27360 & 0 & 0 & 0\\
3 & 0 & 0 & 0 & 381216 & 0 & 0\\
4 & 0 & 0 & 0 & 0 & 5163744 & 0\\
5 & 0 & 0 & 0 & 0 & 0 & 69138720
\end{array}
\end{equation}

\begin{equation}
G_2:\qquad
\begin{array}{c|rrrrrr}
\ell\backslash n & 0 & 1 & 2 & 3 & 4 & 5\\ \hline
0 & 96 & 0 & 0 & 0 & 0 & 0\\
1 & 0 & 1824 & 0 & 0 & 0 & 0\\
2 & 0 & 0 & 27360 & 0 & 0 & 0\\
3 & 0 & 0 & 6 & 381222 & 0 & 0\\
4 & 0 & 0 & 0 & 234 & 5163978 & 0\\
5 & 0 & 0 & 0 & 0 & 5814 & 69144534
\end{array}
\end{equation}

\begin{equation}
G_3:\qquad
\begin{array}{c|rrrrrr}
\ell\backslash n & 0 & 1 & 2 & 3 & 4 & 5\\ \hline
0 & 96 & 0 & 0 & 0 & 0 & 0\\
1 & 0 & 1824 & 0 & 0 & 0 & 0\\
2 & 0 & 0 & 27360 & 0 & 0 & 0\\
3 & 0 & 0 & 12 & 381228 & 0 & 0\\
4 & 0 & 0 & 0 & 468 & 5164212 & 0\\
5 & 0 & 0 & 0 & 0 & 11628 & 69150348
\end{array}
\end{equation}

The following three blocks give $G_1,G_2,G_3$, respectively, in \texttt{graph6} format.
\begingroup

\begin{itemize}
\item[$G_1$:] {\footnotesize
\begin{verbatim}
~?@_?????DQCC?CO?@?OA?C@@OgQ@g_@?AoPO_O??o_A@?CGiA?b@A@G?oPG??`??M_?
F?KAAAOX@GO??P@CcCQ?CDOAQAO?PA?IOOoBQ_O_?CC@BO_O?[oGoC?DCoQCWA?cACOg
?a?g@Q??P_gc@S??_]C?GaG__OwIGw??@H?CAOH?gqE_?D?O?DIDE?O_SO?O_ODS?__W
[`K?@_`SU_@?@aGGAWa??D?_AgfOO?G_cOgOiQ?E?T?OceOFG_CD@?@GeDO?G\CcG?IK
DW@IC?AFS?OIASO???B{????B}@B_C?CAGacOGAQ_?DA@_K?AG`LG@KG?Q?GO?HJ_gC@
AO?M_GO_SAQ?A_dGa?AgB@_OA@IS?PG?aA__ccc?`???eo_O??D?ICIARK_OaDG???F`
H?O___Q???PCP_?W_ga??Gd`OW??OG@I?OCKoP`?I_?O?CGxOod??A?A?@??????????
F~~c?E?IAECS?OWOi@@GG@COa?CL?OKEOBd?G?C@DS?aDE?PAo??SH@CES??A@@D`IOC
Gs_?C?E@F?DC_P?L?F??KOCbAC?B_C?Bo?Go@B?K@cY?c??G?AFkO@G?lACH_oG_?oA?
O[GGc?o_OS?CE?OAPOGM_A@O@?h?O?DEE@?[@_G?IA?C@CG`{?CW?OAAC?G?FB??Bw@A
A?M??G???CWKO?EMg?K_?o??C_po?OIk?PO@C??ESq@?BBQ?AC?G??????F~@`g?@_?C
_???????????@~~~
\end{verbatim}}\vspace{10pt}
\item[$G_2$:] {\footnotesize
\begin{verbatim}
~?@_?????P_cA??GKG?G_?Q?D?O@cI?PHA@??DOC_C??CE_P`?GHO?AOoBQa_??KAL?`
?KA_GA?CY@?@G??SD??KCAP?__EW?AOW@P?T?@??Hs?WKP?G?[@GwGc_G?i?KC?caOCC
@O@IOBG?aG_IO_GAC__COQSCS?_SY[??_U?A?aO@goQA?@`_?DFA?ob?@?PUG`@_?o@A
HSoCOG_AJOc?Do?k?_?A?ASG?l_?Y?COSB@Eh?BcCc@C?GPK@`BkB?C?b?aGGOtC?GKO
?xJ@C_?BL??IOoO?_?B{????B}@C_G@OPH?e?HAA?@G?oWoDa?D?OC_i??CLg?OIEO??
AG?SKOOSGAB@@O?@?WaGHC?ICSG_SGCCQ?qGgO??OH`OYc_g????CtCgG@?AicIO?_AA
EO`GPG?a?_AIDOWBI?C?C?Wt?OCH_D????dA_WS_AK??C?@dQOsAO?O???O?????????
F~~AcQoAoA?G??MH?CPB_?ODOb@@G@_D?GKHS?QoGKI??WCA_og?DSC_Aco?WGOgB@?`
S??Bc@AAWOL?G@?P?B__LGG?Po?a?o?@}?GQ??go@PSeG?G?O?Dac?Cya@C`?W?_??HO
?eh_A?GOO@DEO@?CDACB@@Ag?`g?@?a@jP??d??_J????OxOGWHa??B?A??G^B???W@A
E?EA?G???CcDg?AWw?Q_?W??PCSG?QFc?SO?o??BIYO?BAq??c?O??????F~@`g?@_?C
_???????????@~~~
\end{verbatim}}\vspace{10pt}
\item[$G_3$:] {\footnotesize
\begin{verbatim}
~?@_?????P_cA??GKG?G_?Q?D?O@cIG?Q?PG`?_?@S?c?_?`GC??_Y?bA?CCgO?aa?Cb
@_?OOOACA@??a??SE_?KC@KO__EL?AOW?cgT?@??Hs?WEA?G?[?Tw?oO@cH?OCHAHC?q
?cGGDGOE?QO_COPDCSDG?i?G?C_SZJ????IMA@`K?sAI_?@A?_GdPO`_e?K_g?A??o@_
bSo?JC_C_I?CCAGCWcGA?ASG@J_?i?COS?pEQ?bcCaGS?O?k@`@rB?I?B?aGA?tGGGKO
?[R@H?_BL??A_oaC??B{????B}@C_G?SPAPI?JC?C__ASW??Ou__?goO??P?DOQACO_O
@hGC`???DKB@@O?@?ocGGK?IADG_g_CCa?qGY???CHa_wGG??_E@BCKgG@@?igGo?_AA
EOOWOQ@A@?@GTGMAS?A?G?Wt?g?H_E????dA_kC_CI??G?AfAGDBG?_???O?????????
F~~AcQW?oC?c?AEH?CO_w?`P?a@`G@_AgGKPG?P`GKI??ECB?TG?DSC_AeW?WO?WA@_`
S??@w@CASOX?G@?P??o_XoG?oo?a?o??}?Wo??go@PSeC?G?_?Dcc?Cya@AH?O@G?A@O
?eh_@?GO_@DKO@?CDAE?p?C`?bG?@?cDJG?AY??_T????GjQF?AD??AGA??G^B???W@A
E?EA?G???APDg?AWw?Q_?W??Q`CG?KFc?SO?o??BDeO?Caq??c?O??????F~@`g?@_?C
_???????????@~~~
\end{verbatim}}
\end{itemize}

\endgroup

\section{Magnitude homology of diagonal strongly regular and distance-regular graphs}
\label{sec:diagonal-homology-tables}

By \eqref{eq:rank-of-diagonal} and \zcref{prop:drg-zeta-inverse}, the endpoint distance and the intersection array determine the ranks of the diagonal groups with fixed endpoints in a diagonal distance-regular graph.
The diagonal groups are free abelian by \zcref{prop:freeness-diagonal-part}. The intersection array also gives the number of ordered pairs at each distance by \eqref{eq:drg-sphere-sizes}, so it determines the full magnitude homology by the direct sum in \eqref{eq:chain-hom-direct-sum}.
The off-diagonal groups are $0$, so we list only the diagonal ranks.
The row for distance $r$ gives the sum over ordered pairs of endpoints:
\begin{equation}
\sum_{\substack{(u,u')\in V^2\\d(u,u')=r}}
\rank\MH_{i,i}(G;u,u').
\end{equation}
For fixed endpoints at distance $r$, divide each entry in the row for $r$ by the number $\abs V k_r$ of ordered pairs at that distance to obtain the rank in the corresponding column (see \eqref{eq:drg-sphere-sizes}).

The $4\times4$ rook graph $\mathrm{Ro}$ in \zcref{exam:rook-shrikhande}
has intersection array $\{6,3;1,2\}$.
The distance layer sizes are $(k_0,k_1,k_2)=(1,6,9)$. The table gives the ranks of its magnitude homology groups.
\begin{equation}
\begin{array}{c|rrrrrr}
i & 0 & 1 & 2 & 3 & 4 & 5\\ \hline
d(u,u')=0 & 16 & 0 & 96 & 192 & 816 & 2496\\
d(u,u')=1 & 0 & 96 & 192 & 960 & 3072 & 11040\\
d(u,u')=2 & 0 & 0 & 144 & 576 & 2592 & 9792\\ \hline
\rank\MH_{i,i}(\mathrm{Ro})
& 16 & 96 & 432 & 1728 & 6480 & 23328
\end{array}
\end{equation}
Since $\Mag(\mathrm{Ro};t)=16/(1+3t)^2$, for every $i\ge0$ we have
\begin{equation}
\rank\MH_{i,i}(\mathrm{Ro})=16(i+1)3^i.
\end{equation}
The following graphs are also diagonal by
\zcref{thm:srg-minus2-wmg,thm:diam3-minus3-wmg}.

The Johnson graph $G=J(5,2)$
has intersection array $\{6,2;1,4\}$.
The distance layer sizes are $(k_0,k_1,k_2)=(1,6,3)$.
\begin{equation}
\begin{array}{c|rrrrrr}
i & 0 & 1 & 2 & 3 & 4 & 5\\ \hline
d(u,u')=0 & 10 & 0 & 60 & 180 & 1050 & 5400\\
d(u,u')=1 & 0 & 60 & 180 & 1140 & 5880 & 32220\\
d(u,u')=2 & 0 & 0 & 90 & 480 & 2880 & 15840\\ \hline
\rank\MH_{i,i}(G)
& 10 & 60 & 330 & 1800 & 9810 & 53460
\end{array}
\end{equation}

The halved cube $G=Q_5^+$
has intersection array $\{10,3;1,6\}$.
The distance layer sizes are $(k_0,k_1,k_2)=(1,10,5)$.
\begin{equation}
\begin{array}{c|rrrrrr}
i & 0 & 1 & 2 & 3 & 4 & 5\\ \hline
d(u,u')=0 & 16 & 0 & 160 & 960 & 8880 & 81600\\
d(u,u')=1 & 0 & 160 & 960 & 9280 & 85760 & 809440\\
d(u,u')=2 & 0 & 0 & 400 & 4160 & 41760 & 400960\\ \hline
\rank\MH_{i,i}(G)
& 16 & 160 & 1520 & 14400 & 136400 & 1292000
\end{array}
\end{equation}

The Schl\"afli graph $G$
has intersection array $\{16,5;1,8\}$.
The distance layer sizes are $(k_0,k_1,k_2)=(1,16,10)$.
\begin{equation}
\begin{array}{c|rrrrrr}
i & 0 & 1 & 2 & 3 & 4 & 5\\ \hline
d(u,u')=0 & 27 & 0 & 432 & 4320 & 61182 & 905040\\
d(u,u')=1 & 0 & 432 & 4320 & 63072 & 939600 & 14295312\\
d(u,u')=2 & 0 & 0 & 1890 & 34560 & 564030 & 8817120\\ \hline
\rank\MH_{i,i}(G)
& 27 & 432 & 6642 & 101952 & 1564812 & 24017472
\end{array}
\end{equation}

The Hamming graph $G=H(3,3)$
has intersection array $\{6,4,2;1,2,3\}$.
The distance layer sizes are $(k_0,k_1,k_2,k_3)=(1,6,12,8)$.
\begin{equation}
\begin{array}{c|rrrrrr}
i & 0 & 1 & 2 & 3 & 4 & 5\\ \hline
d(u,u')=0 & 27 & 0 & 162 & 162 & 810 & 1458\\
d(u,u')=1 & 0 & 162 & 162 & 1134 & 2106 & 6966\\
d(u,u')=2 & 0 & 0 & 324 & 648 & 2916 & 7128\\
d(u,u')=3 & 0 & 0 & 0 & 216 & 648 & 2592\\ \hline
\rank\MH_{i,i}(G)
& 27 & 162 & 648 & 2160 & 6480 & 18144
\end{array}
\end{equation}

The Gosset graph $G$
has intersection array $\{27,10,1;1,10,27\}$.
The distance layer sizes are $(k_0,k_1,k_2,k_3)=(1,27,27,1)$.
\begin{equation}
\begin{array}{c|rrrrrr}
i & 0 & 1 & 2 & 3 & 4 & 5\\ \hline
d(u,u')=0 & 56 & 0 & 1512 & 24192 & 535248 & 12894336\\
d(u,u')=1 & 0 & 1512 & 24192 & 548856 & 13329792 & 337172976\\
d(u,u')=2 & 0 & 0 & 13608 & 435456 & 12221496 & 326180736\\
d(u,u')=3 & 0 & 0 & 0 & 12152 & 411264 & 11674152\\ \hline
\rank\MH_{i,i}(G)
& 56 & 1512 & 39312 & 1020656 & 26497800 & 687922200
\end{array}
\end{equation}
\end{document}